\documentclass{eectOF}
\usepackage{amsmath}
\usepackage{amssymb}
\usepackage{paralist}
\usepackage{graphics} 
\usepackage{epsfig} 
\usepackage{graphicx}  \usepackage{epstopdf}
\usepackage[colorlinks=true]{hyperref}
\hypersetup{urlcolor=blue, citecolor=red}

\def\currentvolume{X}
\def\currentissue{X}
\def\currentyear{200X}
\def\currentmonth{XX}
\def\ppages{X--XX}
\def\DOI{10.3934/eect.2021044}

\usepackage{bigints}
\usepackage{empheq}
\allowdisplaybreaks
\usepackage{dsfont}
\usepackage{amsfonts}
\usepackage{graphicx}
\usepackage{shadow}
\usepackage{color}
\usepackage{tcolorbox}
\usepackage[all]{xy}
\usepackage{tikz}
\usepackage{graphicx}
\usepackage{enumitem}
\usepackage{cases}

\newtheorem{theorem}{\textcolor{blue}{Theorem}}[section]
\newtheorem{corollary}{\textcolor{blue}{Corollary}}

\newtheorem{lemma}[theorem]{\textcolor{blue}{Lemma}}
\newtheorem{proposition}{\textcolor{blue}{Proposition}}

\theoremstyle{definition}
\newtheorem{definition}[theorem]{\textcolor{blue}{Definition}}
\newtheorem{remark}{\textcolor{blue}{Remark}}
\newtheorem{example}{\textcolor{blue}{Example}}
\newtheorem*{notation}{\textcolor{blue}{Notation}}

\def \d {\mathrm{~d}}
\def \B {\mathbf{B}}
\def \n {\mathbf{n}}

\title[The VANISHING
VISCOSITY LIMIT]
{On uniform null-controllability of tangential transport-Diffusion equations with vanishing
	viscosity limit}

\author[F. ET-TAHRI, J. A. B\'arcena-Petisco, I. Boutaayamou  AND L. Maniar ]{}

\subjclass{35B25, 35K05, 93B05, 93C20.}
\keywords{Carleman estimates; Uniform controllability; Transport equation; Singular limits; Cost control; Agmon inequality.}

\email{fouad.et-tahri@edu.uiz.ac.ma}
\email{jonasier.barcena@ehu.eus}
\email{dsboutaayamou@gmail.com}
\email{maniar@uca.ma}

\thanks{J.A.B.P was supported by the Grant PID2021-126813NB-I00 
	funded by MCIN/AEI/10.13039/501100011033 and by “ERDF A way of making 
	Europe” and by the grant~IT1615-22 funded the Basque Government.}

\begin{document}
\maketitle

\centerline{\scshape Fouad Et-tahri}

\medskip
{\footnotesize
	
	\centerline{Lab-SIV, Facult\'e des Sciences-Agadir,}
	\centerline{ Universit\'e Ibnou Zohr,  B.P. 8106, Agadir, Morocco }
} 

\bigskip
\centerline{\scshape Jon Asier B\'arcena-Petisco}

\medskip
{\footnotesize
	\centerline{Department of Mathematics, University 
		of the Basque Country UPV/EHU, }
	\centerline{Barrio Sarriena s/n, 48940, Leioa, 
		Spain}
}

\bigskip

\centerline{\scshape Idriss Boutaayamou}

\medskip
{\footnotesize
	\centerline{Lab-SIV, Facult\'e Polydisciplinaire de Ouarzazate,}
	\centerline{ Universit\'e Ibnou Zohr, BP 638, Ouarzazate 45000, Morocco}
	
} 

\bigskip
\centerline{\scshape Lahcen Maniar}

\medskip
{\footnotesize
	\centerline{Cadi Ayyad University, Faculty of Sciences Semlalia}
	\centerline{	LMDP, UMMISCO (IRD-UPMC) B.P. 2390, Marrakesh, Morocco}
} 

\bigskip

	\begin{abstract}
		This paper aims to address an interesting open problem posed in the paper ``Singular Optimal Control for a
		Transport-Diffusion Equation" of Sergio Guerrero and Gilles Lebeau in 2007. The problem involves studying the null-controllability cost of a transport-diffusion equation with Neumann conditions, where the diffusivity coefficient is denoted by $\varepsilon>0$ and the velocity by $\B(x,t)$. Our objective is twofold. Firstly, we investigate the scenario where each trajectory of the tangential velocity $\B$ originating from $\overline{\Omega}$ 
		enters the control region in a shorter time at a fixed entry time.
		By employing Agmon inequalities and Carleman estimates, we establish that the control cost remains bounded for sufficiently small $\varepsilon$ and large control time. Secondly, we explore the case where at least one trajectory fails to enter the control region. In this scenario, we prove that the control cost explodes exponentially when the diffusivity approaches zero and the control time is sufficiently small.
	\end{abstract}
	\maketitle
	\numberwithin{equation}{section}
	\allowdisplaybreaks
	\section{Introduction and main results.}
	Transport-diffusion equations with vanishing diffusivity are widely used to model various physical and biological phenomena. They play a significant role in fluid dynamics by describing the movement of particles while accounting for transport and diffusion effects.
	\par 
	Let $\Omega \subset \mathbb{R}^{d}$, $d\geq 1$ be a bounded open set,  $\Gamma$ denote the boundary of $\Omega$, $\mathbf{n}$ represent the outward unit normal field on $\Gamma$ and $\omega \subset \Omega$ be a nonempty open subset. 
	Throughout this paper, the following notations will be consistently employed:
	\begin{eqnarray*}
		\Omega_{T}=\Omega\times (0,T),\quad \omega_{T}=\omega\times (0,T),\quad \text{and} \quad \Gamma_{T}=\Gamma\times (0,T),
	\end{eqnarray*}
	where $T>0$ is the control time. The main goal of this work is to study the cost of controllability of the following parabolic–transport equation with Neumann boundary conditions:
	\begin{equation}\label{s1}
		\left\{
		\begin{aligned}
			\partial_t y-\varepsilon\Delta y+\mathbf{B}(x,t)\cdot\nabla y &=u(x,t)\mathbf{1}_{\omega} & & \text {in}\;\Omega_{T}, \\
			\partial_{\mathbf{n}} y &=0 & & \text {on}\;\Gamma_{T}, \\
			y(x,0)&=y_{0}(x) & & \text { in } \Omega,  \\	
		\end{aligned}
		\right.
	\end{equation}
	where $\varepsilon>0$ is the viscosity (diffusion coefficient) and $\mathbf{B}$ is the velocity (speed), that satisfies
	\begin{equation}
		\mathbf{B}\in W^{1,\infty}(\mathbb{R}^{d}\times (0,\infty))^{d},
	\end{equation}
	as considered in \cite{guerrero2007singular} in the case of uniform controllability in $\varepsilon$ for Dirichlet conditions case. The reason for this spatial extension (a regular open strictly containing $\Omega$ is sufficient) is to define geometric conditions for the trajectories of the vector field $\B$ and time extension, to have norms of $\B$ that do not depend on $T$. In the case of the controllability cost explosion, we only assume that 
	\begin{equation}
		\mathbf{B}\in L^{\infty}(0,T; W^{1,\infty}(\Omega)^{d}).
	\end{equation}
	In system \eqref{s1}, $y=y(x,t)$ represents the state, $u\in L^{2}(\omega_{T})$ is the control function that only affects the system through $\omega_{T}$, $\mathbf{1}_{\omega}$ denotes the characteristic function of $\omega$ and $y_{0}\in L^{2}(\Omega)$ is the initial state. Let us start with the definition of null-controllability:
	\begin{definition}
		We say that system \eqref{s1} is null controllable at time $T$ if for every state $y_{0}\in L^{2}(\Omega)$, there exists a control $u\in L^{2}(\omega_{T})$ such that the solution of \eqref{s1} satisfies $y(\cdot,T)=0$.
	\end{definition}
	\par
	It is well known that the system \eqref{s1} is null controllable for any time $T$ and any control region $\omega$ when $\mathbf{B}\in L^{\infty}(\Omega_{T})$. Specifically, in the case where $\varepsilon=1$, we can refer to \cite{fernandez2006exact} and \cite{fursikov1996controllability}	 for further details. Additionally, by rescaling the time variable using the transformation $t\mapsto\frac{t}{\varepsilon}$, we can effectively reduce the problem to this specific case, as shown in \cite[Proposition 3.1]{ettahri:hal-04131920}. It is a well-known fact that the controls depend continuously on the initial data. In other words, there exists a constant $C:=C(T,\varepsilon)>0$, such that
	\begin{equation}
		\|u\|_{L^{2}(\omega_{T})}\leqslant C\|y_{0}\|_{L^{2}(\Omega)} \label{controllability}
	\end{equation}
	\par
	and the null-controllability 
	is equivalent to the following observability inequality:
	\begin{eqnarray}
		\exists\,\mathcal{C}:=\mathcal{C}(T,\varepsilon)>0, \; \forall \varphi_{T}\in L^{2}(\Omega),\; \|\varphi(\cdot, 0)\|_{L^{2}(\Omega)}\leq \mathcal{C} \|\varphi\|_{L^{2}(\omega_{T})}, \label{observability}
	\end{eqnarray}
	where $\varphi$ is the solution of the adjoint system of \eqref{s1}:
	\begin{equation}
		\left\{
		\begin{aligned}
			\partial_t \varphi+\varepsilon\Delta \varphi+\nabla \cdot\left(\varphi \mathbf{B}(x,t)\right)  &=0 & & \text { in } \Omega_{T}, \\
			\left(\varepsilon \nabla\varphi +\varphi \mathbf{B}(t,x)\right)\cdot\mathbf{n}(x) &=0 & & \text { on } \Gamma_{T}, \\
			\varphi(x,T)&=\varphi_{T}(x) & & \text { in } \Omega.
			\label{s2}
		\end{aligned}
		\right.
	\end{equation}
	\par
	By employing the Hilbert Uniqueness Method \cite{russell1978controllability, lions1988controlabilite}, it can be shown that the optimal constants satisfying \eqref{controllability} and \eqref{observability} are equal that is:
	\begin{equation}
		\sup_{y_{0}\in L^{2}(\Omega)\setminus\{0\}}\inf_{u\in \mathbb{A}(y_{0})}\frac{\lvert\lvert u\lvert\lvert_{L^{2}(\omega_{T})}}{\lvert\lvert y_{0}\lvert\lvert_{L^{2}(\Omega)}}=\sup_{\varphi_{T}\in L^{2}(\Omega)\setminus\{0\}}\frac{\lvert\lvert \varphi(\cdot, 0)\lvert\lvert_{L^{2}(\Omega)}}{\lvert\lvert \varphi\lvert\lvert_{L^{2}(\omega_{T})}},
	\end{equation}
	where $\varphi$ is the solution of the adjoint system \eqref{s2} and $$\mathbb{A}(y_{0}):=\{u\in L^{2}(\omega_{T}): \;\text{the solution of}\;\eqref{s1}\;\text{satifies}\;y(\cdot,T)=0\}.$$
	\par
	In the sequel, we will adopt the following definition:
	\begin{definition} 
		We define the cost of null-controllability of system \eqref{s1} by the following quantity: 
		\begin{eqnarray}
			\mathcal{K}(\varepsilon,T,\Omega,\omega):=\sup_{\varphi_{T}\in L^{2}(\Omega)\setminus\{0\}}\frac{\lvert\lvert \varphi(\cdot, 0)\lvert\lvert_{L^{2}(\Omega)}}{\lvert\lvert \varphi\lvert\lvert_{L^{2}(\omega_{T})}}. \label{cost of control}
		\end{eqnarray}
	\end{definition}
	The main objective of this paper is to investigate the asymptotic behavior of the null-controllability cost for system \eqref{s1} when the viscosity is small enough. To elucidate the main findings of this paper, let us examine the trajectories of the vector field $\mathbf{B}$ given by the mapping $t\mapsto \varPhi(t,t_{0},x_{0})$:
	\begin{equation} \label{OD}
		\left\{
		\begin{aligned}
			\frac{d}{dt}\varPhi(t,t_{0},x_{0}) &=\mathbf{B}(\varPhi(t,t_{0},x_{0}),t) & & t\in (0,T), \\
			\varPhi(t_{0},t_{0},x_{0}) &=x_{0},
		\end{aligned}
		\right.
	\end{equation}
	for each $(x_{0},t_{0})\in \mathbb{R}^{d}\times [0,T]$. The solutions of the ordinary differential equation \eqref{OD} encompass all the relevant information regarding the trajectories of a particle moving with the velocity $\B$. The following definition serves a specific purpose that is essential to prove our first main result.
	\begin{definition} \label{definition}
		Let $T_{0}\in (0,T)$, $r_{0}>0$ and  $\mathcal{O}\subset\mathbb{R}^{d}$ non-empty open. We say that $(T,T_{0},r_{0},\mathbf{B},\Omega)$ satisfies the Flushing condition \eqref{Flushing Condition} for $\mathcal{O}$ 
		if 
		\begin{equation}
			\forall x_{0}\in\overline{\Omega},\;\forall t_{0}\in [T_{0},T],\; \exists\, t\in (t_{0}-T_{0},t_{0}),\;\forall x\in\overline{B}(x_{0},r_{0}), \varPhi(t,t_{0},x_{0})\in \mathcal{O}. \tag{$\mathcal{F}\mathcal{C}$} \label{Flushing Condition}
		\end{equation}
	\end{definition}
	\begin{remark}
		In the case \eqref{OD} is autonomous, i.e.,  $\B=\B(x)$, we can characterize \eqref{Flushing Condition} by any backward trajectories of $\B$ originating from $\overline{\Omega}$ at time $0$ and entering the open set $\mathcal{O}$. This characterization is discussed in Proposition \ref{P5} in Section \ref{Section 2}.
	\end{remark}
	We will show that if every backward trajectory of $\mathbf{B}(x,t)$ starting from $\overline{\Omega}$ enters the control region within a time that does not surpass a fixed time barrier, then the cost of null-controllability of \eqref{s1} remains uniformly bounded with respect to $\varepsilon$ when it is small enough and the control time is sufficiently large. To be more precise, we will establish the following theorem:
	\begin{theorem}\label{m1}
		Under the following conditions:
		\begin{enumerate}[label=\textcolor{red}{(\arabic*)}]
			\item $\Omega\subset \mathbb{R}^{d}$ is a $\mathcal{C}^{2}$ domain, $d\geq 2$ and $\omega\subset\subset\Omega$ is an open subset, 
			\item $\mathbf{B}\in W^{1,\infty}(\mathbb{R}^{d}\times (0,\infty))^{d}$,
			\item \label{cc3} there exist $T_{0}\in (0,T)$ and $r_{0}>0$ such that $(T,T_{0},r_{0},\mathbf{B},\Omega)$ satisfies \eqref{Flushing Condition} for $\omega$, 
			\item $\forall x\in\Gamma,\forall t\in [0,T],\;\mathbf{B}(x,t)\cdot\mathbf{n}(x)=0$.
		\end{enumerate}
		There exists a constant $\rho_{0}\geq 1$ depending only on $T_{0}$, $r_{0}$ and $\|\mathbf{B}\|_{W^{1,\infty}(\mathbb{R}^{d}\times (0,\infty))}$ such that if $T\geq \rho_{0} T_{0}$, there exists a 
		constant $C>0$ independent of $\varepsilon$ that satisfies the following estimate:
		\begin{equation} \label{cost1}
			\mathcal{K}(\varepsilon,T,\Omega,\omega)\leqslant C,\quad\text{for}\;\varepsilon\;\text{small enough,}
		\end{equation}
		where $\mathcal{K}$ the cost of the null controllability of \eqref{s1}.
	\end{theorem}
	Before giving an example that illustrates the condition \ref{C3m2} of Theorem \ref{m1}, we introduce some notations that will be useful in the sequel:
	\begin{notation}
		\begin{enumerate}
			\item The canonical Euclidean scalar product of $\mathbb{R}^{d}$ is denoted by $\cdot$ and $\left|\cdot\right|$ stands for the
			associated canonical Euclidean norm.
			\item For all $x,y\in\mathbb{R}^{d}$ and $r>0$, $B(x,r)$  and $\overline{B}(x,r)$ are the open and closed balls of center $x$ and radius $r$, respectively.
			\item For $A$ and $B$ two parts of $\mathbb{R}^{d}$,  $\mbox{dist}(x,A)$ and $\mbox{dist}(A,B)$ designate the distance from $x$ to $A$ and the distance between $A$ and $B$, respectively.
		\end{enumerate} 
	\end{notation}
	In 1-D, the conditions of Theorem \ref{m1} are not satisfied because the extremities of $\Omega$ are points of equilibria of \eqref{OD}. In 2-D, the following example illustrates the conditions of Theorem \ref{m1}:
	\begin{example}
		Let
		$\begin{cases}
			&\Omega=B(0,R),\; \omega=B(0,r),\;\mbox{for}\; 0<r<R,\\
			& \B(x)=e^{-|x|^{2}}\left[(-x_{2},x_{1})+(R^{2}-\left|x\right|^{2})x\right],\;\mbox{for}\; x=(x_{1},x_{2})\in\mathbb{R}^{2}.
		\end{cases}$
			Indeed, in this case $\n(x)=\frac{x}{R}$, then $\n(x)\cdot \B(x)=0$ for all $\left|x\right|=R$ and $0$ is an asymptotically stable equilibrium point of the system:
			\begin{equation} 
				\left\{
				\begin{aligned}
					\frac{d}{dt}\Psi(t,0,x_{0}) &=-\mathbf{B}(\Psi(t,0,x_{0})) & & t\geq 0, \\
					\Psi(0,0,x_{0}) &=x_{0}.
				\end{aligned}
				\right.
			\end{equation}
			Hence, there exists $R>0$ such that for all $r<R$, we have
			
			$$\forall x_{0}\in \overline{B}(0,R),\;\;\exists t\in (-\infty,0),\;\; \Phi(t,0,x_{0})=\Psi(-t,0,x_{0})\in B(0,r).$$
			Then, the condition \ref{cc3} of Theorem \ref{m1} is satisfied (see above remark).
			This property is illustrated in the following figure:
			\par
			\begin{figure}[!h]
				\centering
				\includegraphics[scale=0.5]{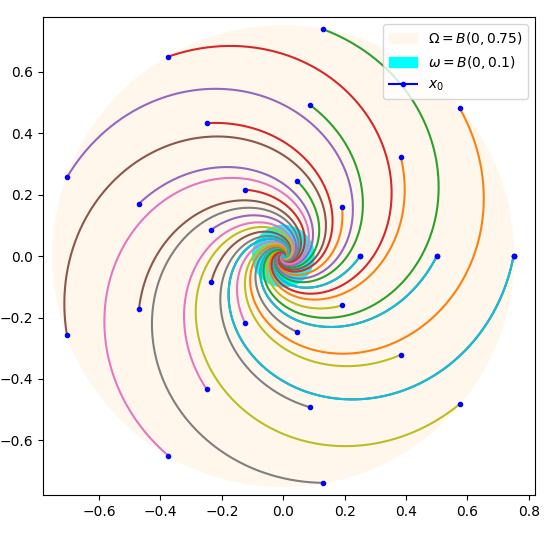}
				\caption{Backward trajectories $t\mapsto\varPhi(t,0,x_{0})$ starts in $\overline{\Omega}$ enter $\omega$.}
			\end{figure}
		\end{example}
	\par
	The second main result of this paper is to show that if there is a backward trajectory of $\B$ that starts in $\overline{\Omega}$ at time $T$ stays in $\overline{\Omega}$ and does not enter the control region during the time $[0,T]$, as illustrated in the following figure, then for a small time control, the control cost explodes exponentially 
	when the viscosity  vanishes. To be more precise, we provide a proof of the following theorem.
	\begin{theorem} \label{m2}
		We assume that:
		\begin{enumerate}[label=\textcolor{red}{(\arabic*)}]
			\item $\Omega\subset\mathbb{R}^{d}$ is a domain with Lipschitz boundary $\Gamma$, $d\geq 1$ and $\omega\subset\Omega$ an open subset, 
			\item $\mathbf{B}\in L^{\infty}(0,T; W^{1,\infty}(\Omega)^{d})$,
			\item \label{C3m2} $\exists x_{0}\in\Omega$ such that for all $t\in [0,T]$, $\varPhi(t,T,x_{0})\in \Omega\setminus\overline{\omega}$.
		\end{enumerate}
		Then there exists a 
		constant $C>0$ independent of $\varepsilon$ such that, we have the following estimate:
		\begin{equation} \label{cost2}
			\mathcal{K}(\varepsilon,T,\Omega,\omega)\geq\exp\left(\frac{C}{\varepsilon}\right),\quad\text{for}\;\varepsilon\;\text{small enough,}
		\end{equation}
		where $\mathcal{K}$ the cost of the null controllability of \eqref{s1}.
	\end{theorem}
Generally, the condition \ref{C3m2} of Theorem \ref{m2} holds true provided that $T$ is sufficiently small, using a continuity argument. For instance:
	\begin{example} 
		 Let $x_{0}\in \Omega\setminus\overline{\omega}$, $r=\mbox{dist}(x_{0},\omega)$ and $b=\|\B\|_{L^{\infty}(\Omega\times(0,\infty))}$ for $\B\in L^{\infty}(\Omega\times (0,\infty))^{d}$. one has
		\begin{eqnarray*}
			\left|\varPhi(t,T,x_{0})-x_{0}\right| &=&\left|\int_{T}^{t}\B(\varPhi(s,T,x_{0}),s)\d s\right|\leq Tb,\; \mbox{for all}\; t\in [0,T].
		\end{eqnarray*}
		Taking $T>0$ such that $0< Tb <r$, we obtain $\varPhi(t,T,x_{0})\in B(x_{0},r)\subset \Omega\setminus\overline{\omega}$.
		\par 
		We can have this property for all $T>0$, as the following example shows: In 2-D, let $\mathbf{B}(x_{1},x_{2})=(x_{2},-x_{1})$. In this case, the matrix associated with the equation \eqref{OD} is skew-symmetric. As a result, we have $\left|\varPhi(t,T,x_{0})\right|=\left|x_{0}\right|$ for any $x_{0}\in\mathbb{R}^{2}$ and any $t\in \mathbb{R}$, as shown in the following figure:
	\begin{figure}[!h] \label{fig1}
		\centering
		\includegraphics[scale=0.5]{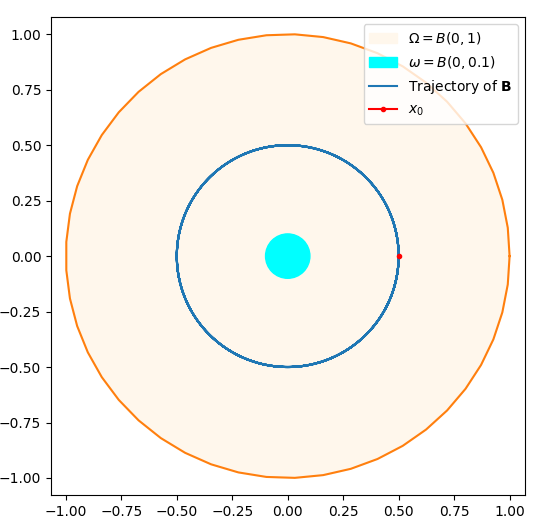}
		\caption{A trajectory of $\B$ in $\Omega\setminus\overline{\omega}$.} 
	\end{figure}
\end{example}
\begin{remark}
	Theorem \ref{m2} is a generalization of \cite[Theorem 2.8]{barcena2021cost} for a small time control. Indeed, the conditions specified in Theorem \cite[Theorem 2.8]{barcena2021cost} lead directly to the realization of the condition \ref{C3m2} of Theorem \ref{m2} for $T<h$ and $x_{0}\in (p_{l}+T,p_{l}+h)$.
\end{remark}
		\par 
		 In the context of the problem under study, Sergio Guerrero and Gilles Lebeau established analog results in their work, specifically in \cite[Theorem 1 and Theorem 3]{guerrero2007singular}, with a transport flow belonging to $W^{1,\infty}(\mathbb{R}^{d}\times (0,+\infty))$ and employing Dirichlet conditions. Jon Asier B\'arcena-Petisco, in \cite[Theorem 2.7 and Theorem 2.8]{barcena2021cost}, demonstrated the same outcomes for the case of the first vector of the canonical basis of $\mathbb{R}^{d}$ as a velocity  and autonomous Robin (or Fourier) conditions. The generalization of these results to a velocity field expressed as a gradient belonging to $W^{1,\infty}(\Omega)^{d}$ was accomplished by \cite{ettahri:hal-04131920}, extending the findings from \cite[Theorem 2.7 and Theorem 2.8]{barcena2021cost}. Additionally, in \cite{laurent2021uniform} Camille Laurent and Matthieu Léautaud investigated uniform controllability and the corresponding optimal time for homogeneous Dirichlet conditions on a smooth, connected, compact manifold, while in \cite{laurent2022uniform} they considered the analogous scenario for 1-D systems.
		 \par 
		 Our contribution is to answer the open questions presented in \cite[Remark 3]{guerrero2007singular} and \cite{barcena2021cost}. The objective then is to study the null-controllability cost of a transport-diffusion equation with Neumann conditions and velocity $\B(x,t)\in W^{1,\infty}(\mathbb{R}^{d}\times (0,\infty))^{d}$ in the case of uniform controllability.
		 The main difficulty encountered initially was to establish an Agmon's inequality. However, this obstacle has been overcome with the help of estimate \eqref{Trace estimate} while still using the tools presented in \cite{guerrero2007singular}, which allowed for the problem's resolution. Subsequently, the second challenge was to establish a new Carleman estimate, yielding an observability constant of the form $e^{C/\varepsilon}$. 
		 Additionally, it should be noted that the vector $\B$ depends on both $x$ and $t$ and does not take the form of a gradient, preventing the use of a state change to eliminate the transport term, as is possible in cases \cite{barcena2021cost, ettahri:hal-04131920}. In the proof of Theorem \ref{m1}, we use a new decomposition (see the system introduced in Section \ref{Section 3}) of the adjoint system \eqref{s2} to prove estimates inside the control region and outside, using the uniform Agmon inequality in $\varepsilon$ and the Carleman estimate.  In the proof of Theorem \ref{m2}, we construct a solution of the adjoint system \eqref{s2}, considering a smooth initial data and exploiting Agmon's inequality.
		\par 
		From a historical standpoint, the problem under investigation has its roots in the field of controllability problems within singular limits, which were initially introduced by Jacques-Louis Lions. The articles \cite{lopez2000null} of Antonio L\'opez, Xu Zhang, and Enrique Zuazua and \cite{phung2002null} of Kim Dang Phung provide an illustration of how the null controllability of the heat equation emerges as a singular limit from the exact controllability of dissipative wave equations. Afterwards, the specific focus of the study revolves around the evanescent viscosity limit, which was first introduced by Jean-Michel Coron and Sergio Guerrero in their work \cite{coron2005singular}. Initially explored in the context of 1-D transport equations, this problem was subsequently extended to higher dimensions by Sergio Guerrero and Gilles Lebeau in their work cited as \cite{guerrero2007singular}. To explore other references related to optimal time of the controllability with
		vanishing viscosity limit of the heat equation see \cite{glass2010complex} of Olivier Glass and \cite{lissy2012link, lissy2014application, lissy2015explicit} of Pierre Lissy and the reference therein.  For the motivation behind the problem and other applications, we invite you to refer to the introductions of the mentioned works \cite{laurent2021uniform, barcena2021cost, ettahri:hal-04131920}.  A last interesting paper involving parabolic equations, in this case fourth order equation, is \cite{carreno2016cost} of Nicolas Carreño and Patricio Gúzman.
		\par 
		Our paper is organized as follows. In Section \ref{Section 2}, we present several results related to property \eqref{Flushing Condition}. Moving on to Section \ref{Section 3}, we delve into the examination of the existence and uniqueness of both strong and weak solutions for a parabolic system that includes the adjoint system \eqref{s2}. In Section \ref{Section 4}, we will present some technical results and we will prove a new Agmon inequalities.
			Transitioning to Section \ref{Section 5}, our focus shifts to establishing significant dissipation results. In Section \ref{Section 6}, we prove a new  Carleman estimate for the adjoint system \eqref{s2}.
			Section \ref{Section 7} is devoted to the proof of our main results, Theorem \ref{m1} and Theorem \ref{m2}. In the last section, we cite some open questions and problems.
		\section{Some preliminary results of condition \eqref{Flushing Condition}} \label{Section 2}
		In this section, we present two relevant results concerning property \eqref{Flushing Condition}. The first result offers a characterization of this property specifically in the autonomous case, while the second result introduces a refinement of the regions associated with this property, which will be used later in our analyses.\\
		\par 
		The following lemma ensures the existence and uniqueness of differentiable solutions of the ordinary differential equation \eqref{OD}.
		\begin{lemma}  Let $\mathbf{B}\in L^{\infty}(0,T;W^{1,\infty}(\mathbb{R}^{d})^{d})$.
			For all $(x_{0},t_{0})\in \mathbb{R}^{d}\times [0,T]$, the ordinary differential equation \eqref{OD} admits a unique global differentiable solution $\varPhi$. Moreover, for all $t\in [0,T]$ and $(x_{0},t_{0}), (y_{0},s_{0})\in\mathbb{R}^{d}\times [0,T]$
			\begin{eqnarray}
				&&|\varPhi(t,t_{0},x_{0})-\varPhi(t,s_{0},y_{0})|\nonumber\\
				&& \leqslant \exp(\|\nabla\mathbf{B}\|_{L^{\infty}(\mathbb{R}^{d}\times (0,T))}T)(\|\mathbf{B}\|_{L^{\infty}(\mathbb{R}^{d}\times (0,T))}|t_{0}-s_{0}|+|x_{0}-y_{0}|). \label{f1}
			\end{eqnarray} 
		\end{lemma}
		\begin{proof}
			For all $t\in (0,T)$, $\mathbf{B}(\cdot,t)\in W^{1,\infty}(\mathbb{R}^{d})^{d}$ and $\mathbb{R}^{d}$ (is convex),  from \cite{brezis1983analyse} for scalar-valued functions, we deduce that
			\begin{eqnarray}
				|\mathbf{B}(x,t)-\mathbf{B}(y,t)|&\leqslant & \|\nabla\mathbf{B}(\cdot,t)\|_{L^{\infty}(\mathbb{R}^{d})}|x-y| \nonumber\\
				&\leqslant & \|\nabla\mathbf{B}\|_{L^{\infty}(\mathbb{R}^{d}\times (0,T))}|x-y|.\label{odf1}
			\end{eqnarray}
			where $\B:=(\B_{1},\cdots,\B_{d})$ and $\|\nabla\mathbf{B}(\cdot,t)\|^{2}_{L^{\infty}(\mathbb{R}^{d})}:=\displaystyle\sum_{1\leq i,j\leq d}\|\partial_{x_j}\mathbf{B}_{i}(\cdot,t)\|^{2}_{L^{\infty}(\mathbb{R}^{d})}$.
			The Cauchy-Lipschitz theorem affirms that, for all $(x_{0},t_{0})\in \mathbb{R}^{d}\times [0,T]$, \eqref{OD} has a unique global solution given by 
			\begin{eqnarray}
				\varPhi(t,t_{0},x_{0})=x_{0}+\int_{t_{0}}^{t}\mathbf{B}(\varPhi(\tau,t_{0},x_{0}),r)\mathrm{~d}\tau,\quad 0\leqslant t\leqslant T. \label{odf2}
			\end{eqnarray}
			Let $(x_{0},t_{0}), (y_{0},s_{0})\in \mathbb{R}^{d}\times [0,T]$ with $t_{0}\leqslant s_{0}$, from \eqref{odf1} and \eqref{odf2}, we obtain
			\begin{eqnarray}
				&&\left|\varPhi(t,t_{0},x_{0})-\varPhi(t,s_{0},y_{0})\right| \leqslant |x_{0}-y_{0}|+ \|\mathbf{B}\|_{L^{\infty}(\mathbb{R}^{d}\times (0,T))}|s_{0}-t_{0}| \nonumber\\
				&&\quad + \|\nabla\mathbf{B}\|_{L^{\infty}(\mathbb{R}^{d}\times (0,T))}\int_{\min(s_{0},t)}^{\max(s_{0},t)}\left|\varPhi(\tau,t_{0},x_{0})-\varPhi(\tau,s_{0},y_{0})\right|\mathrm{~d}\tau.\nonumber
			\end{eqnarray}
			Applying Gronwall's lemma to this last inequality, we obtain \eqref{f1}.
		\end{proof}
		\par 
		In the autonomous case, we can characterize the condition \eqref{Flushing Condition} as follows: 
		\begin{proposition} \label{P5}
			Let $\mathcal{O}$ a non-empty open of $\mathbb{R}^{d}$ and $\mathbf{B}\in W^{1,\infty}(\mathbb{R}^{d})^{d}$.\\
			Assume that 
			$$\forall x_{0}\in\overline{\Omega},\exists\, t\in (-\infty,0),\; \varPhi(t,0,x_{0})\in\mathcal{O}.$$
			Then there exist $T_{0}>0$ and $r_{0}>0$ such that for all $T>T_{0}$, $(T,T_{0},r_{0},\mathbf{B},\Omega)$ satisfies  the condition \eqref{Flushing Condition} for $\mathcal{O}$.
		\end{proposition}
		\begin{proof}
			For all $x_{0}\in\overline{\Omega}$ there exists $t:=t(x_{0})\in (-\infty,0)$ such that $ \varPhi(t(x_{0}),0,x_{0})\in\mathcal{O}.$ By the continuity of the flow, there exists $r:=r(x_{0})>0$ such that 
			\begin{eqnarray}
				|x-x_{0}|< 2r(x_{0})\Longrightarrow \varPhi(t(x_{0}),0,x)\in\mathcal{O}. \label{od4}
			\end{eqnarray}
			By compactness of $\overline{\Omega}$, there exist $x^{1}_{0},\cdots,x^{d}_{0}\in \overline{\Omega}$ such that $\overline{\Omega}\subset\displaystyle\bigcup_{i=1,\cdots, d}B(x^{i}_{0},r(x^{i}_{0}))$. We put $T_{0}:=-\displaystyle\min_{i=1,\cdots, d}t(x^{i}_{0})>0$ and $r_{0}:= \displaystyle\min_{i=1,\cdots, d}r(x^{i}_{0})> 0$.\\ From \eqref{od4}, we obtain
			\begin{eqnarray}
				\forall x_{0}\in\overline{\Omega},\exists\, t\in (-T_{0},0),\;|x-x_{0}|< r_{0}\Longrightarrow \varPhi(t,0,x)\in\mathcal{O}. \label{od5}
			\end{eqnarray}
			Let $T>T_{0}$, $x_{0}\in\overline{\Omega}$ and $t_{0}\in [T_{0},T]$. Since $\mathbf{B}$ is independent of $t$, then
			\begin{eqnarray}
				\varPhi(t,t_{0},x)=\varPhi(t-t_{0},0,x),\; t\in\mathbb{R},\;x\in \mathbb{R}^{d}. \label{od6}
			\end{eqnarray}
			From \eqref{od5} and \eqref{od6}, we conclude that $(T,T_{0},r_{0},\mathbf{B},\Omega)$ satisfies the condition \eqref{Flushing Condition} for $\mathcal{O}$.
		\end{proof}
		The following proposition guarantees that the condition \eqref{Flushing Condition} remains true for small regions.
		\begin{proposition} \label{P1}
			Let $\mathcal{O}\subset\mathbb{R}^{d}$ an open non-empty with bounded boundary and assume that $(T,T_{0},r_{0},\mathbf{B},\Omega)$ satisfies the condition \eqref{Flushing Condition} for $\mathcal{O}$. Then there exists $\mathcal{O}_{0}\subset \subset \mathcal{O}$ an open such that $\left(T,T_{0},\frac{r_{0}}{2},\B,\Omega\right)$ satisfies the Condition \eqref{Flushing Condition} for $\mathcal{O}_{0}$.
		\end{proposition}
		\begin{proof}
			For all $(x_{0},t_{0})\in \overline{\Omega}\times [T_{0},T]$, there exists $t:=t(x_{0},t_{0})\in (t_{0}-T_{0},t_{0})$, such that for all $x\in B(x_{0},r_{0})$ we have $\varPhi(t(x_{0},t_{0}),t_{0},x)\in\mathcal{O}$. We put
			\begin{eqnarray}
				d(x_{0},t_{0}):=\text{dist}\left(\left\{\varPhi(t(x_{0},t_{0}),t_{0},x):\;x\in \overline{B}\left(x_{0},\frac{r_{0}}{2}\right)\right\},\partial\mathcal{O}\right)>0, \label{od7}
			\end{eqnarray}
			since $\left\{\varPhi(t(x_{0},t_{0}),t_{0},x):\;x\in \overline{B}\left(x_{0},\frac{r_{0}}{2}\right)\right\}$ is a closed thanks to \eqref{f1} and $\partial\mathcal{O}$ is a compact.  The continuity of the flow in \eqref{f1}, asserts that it exists $r:=r(x_{0},t_{0})>0$, such that, for all $|s-t_{0}|<r(x_{0},t_{0})$, we have
			\begin{eqnarray}
				\forall x\in \overline{B}(x_{0},r_{0}),\; |\varPhi(t(x_{0},t_{0}),s,x)-\varPhi(t(x_{0},t_{0}),t_{0},x)|\leq \frac{d(x_{0},t_{0})}{2}. \label{od8}
			\end{eqnarray}
			We then consider $\mathcal{U}_{x_{0},t_{0}}$ the set of couples $(x,t)$ satisfying 
			$$x\in B\left(x_{0},\frac{r_{0}}{2}\right),\quad t-T_{0}<t(x_{0},t_{0})<t\quad \text{and}\quad t\in (t_{0}-r(x_{0},t_{0}),t_{0}+r(x_{0},t_{0})).$$
			Since $\mathcal{U}_{x_{0},t_{0}}$ is on open containing $(x_{0},t_{0})$ and $\overline{\Omega}\times [T_{0},T]$ is compact, then it admits a finite covering by $\mathcal{U}_{x^{i}_{0},t^{i}_{0}}$, $i=1,\cdots,I$. Taking $d_{0}=\displaystyle\min_{i=1,\cdots,I}d(x^{i}_{0},t^{i}_{0})$ and using \eqref{od7} and \eqref{od8}, we can show that  $\left(T,T_{0},\frac{r_{0}}{2},\B, \Omega\right)$ satisfies  the condition \eqref{Flushing Condition} for all open $\mathcal{O}_{0}$ such that $$\left\lbrace x\in\mathcal{O}:\;\text{dist}(x,\partial\mathcal{O})\geq \frac{d_{0}}{2}\right\rbrace \subset \mathcal{O}_{0}.$$
		\end{proof}
		\section{Wellposedness and results of a parabolic equation including the adjoint system \eqref{s2}.} \label{Section 3}
		In this section, we will establish the well-posedness and regularity properties of solutions for the following backward, inhomogeneous linear transport-diffusion equation, accompanied by mixed boundary conditions Dirichlet and non-autonomous Robin conditions:
		\begin{equation}  \label{s3}
			\left\{
			\begin{aligned}
				\partial_t \varphi+\varepsilon\Delta \varphi+\nabla \cdot\left(\varphi \mathbf{B}(x,t)\right) &=f(x,t) & & \text { in } \mathcal{U}\times (t_{1},t_{2}), \\
				\left(\varepsilon \nabla\varphi +\varphi \mathbf{B}(x,t)\right)\cdot \mathbf{n}(x)\mathbf{1}_{\Gamma}(x)+\varphi\mathbf{1}_{\Gamma_{0}}(x) &=0 & & \text { on } \partial\mathcal{U}\times (t_{1},t_{2}), \\
				\varphi(x,t_{2}) &=g(x) & & \text { in } \mathcal{U},\\
			\end{aligned} 
			\right.
		\end{equation}
		where $0\leqslant t_{1}<t_{2}\leqslant T$ and $\Omega_{0}\subset\subset\Omega$ a regular open, $\mathcal{U}:=\Omega\setminus\overline{\Omega_{0}}$, $\Gamma=\partial\Omega$, $\Gamma_{0}:=\partial\Omega_{0}$, $f\in L^{2}(\mathcal{U}\times (t_{1},t_{2}))$ and $g\in L^{2}(\mathcal{U})$.\\ 
		\begin{notation}
			In the following, $\mathcal{S}(\Omega_{0},t_{1},t_{2},f,g,\varepsilon,\mathbf{B})$ will refer to system \eqref{s3}.
		\end{notation}
		Note that in the case $\Omega_{0}=\varnothing$ the system $\mathcal{S}(\varnothing,t_{1},t_{2},f,g,\varepsilon,\mathbf{B})$ is the adjoint system \eqref{s2}.
		\subsection{Notations and function spaces}
		Let $\Omega\subset\mathbb{R}^{d}$ a domain with Lipschitz boundary $\Gamma$, $L^{2}(\Omega)$ and $L^{2}(\Gamma)$ are the classical Hilbert spaces over $\mathbb{R}$ with respect to the Lebesgue measure $\mathrm{~d} x$ on $\Omega$ and the $(d-1)$-dimensional Hausdorff measure $\mathrm{~d}\sigma$ on $\Gamma$ and $(\cdot,\cdot)$ is the canonical scalar product of $L^{2}(\Omega)$. We consider $H^{1}(\Omega)$ and $W^{1,\infty}(\Omega)$ are the usual $L^2$-based and $L^{\infty}$-based Sobolev spaces over $\Omega$, respectively. For any $\Omega_{0}\subset\subset\Omega$ a regular open, we set 
		\begin{equation*}
			\mathcal{U}=\Omega\setminus\overline{\Omega_{0}}\quad\mbox{and}\quad \Gamma_{0}=\partial\Omega_{0}.
		\end{equation*}
		We introduce  $H^{1}_{\Gamma_{0}}(\mathcal{U})$ the space of all those functions in $H^{1}(\mathcal{U})$ whose trace vanishes on $\Gamma_{0}$:
		\begin{eqnarray}
			H^{1}_{\Gamma_{0}}(\mathcal{U}):= \begin{cases}
				\{u\in H^{1}(\mathcal{U}):\;\;u=0\;\text{on}\;\Gamma_{0}\},\quad &\text{if}\; \Omega_{0}\neq \varnothing,\\
				H^{1}(\Omega) ,\quad &\text{if}\; \Omega_{0}= \varnothing.
			\end{cases}
		\end{eqnarray}
		We will keep this space the induced norm of $H^{1}(\mathcal{U})$.
		We note by $\left(H_{\Gamma_{0}}^{1}(\mathcal{U})\right)^{'}$ the dual of $H_{\Gamma_{0}}^{1}(\mathcal{U})$ and the product duality is denoted by $\langle\;\cdot\;,\;\cdot\rangle_{\left(H_{\Gamma_{0}}^{1}(\mathcal{U})\right)^{'}, H_{\Gamma_{0}}^{1}(\mathcal{U})}$. Clearly $H^{1}_{\Gamma_{0}}(\mathcal{U})$ is dense in $L^{2}(\mathcal{U})$, as usual we can identify $L^{2}(\Omega)$ with a
		dense subspace of $\left(H^{1}_{\Gamma_{0}}(\mathcal{U})\right)^{'}$.
		We will also note $D(\Omega)$ the space of the test functions on $\Omega$. We recall that there exists a unique linear bounded operator $\gamma_{0} : H^{1}(\Omega)\longrightarrow L^{2}(\Gamma)$ such that $\gamma_{0}(u)=\left.u\right|_{\Gamma}$ if $u\in H^{1}(\Omega)\cap \mathcal{C}(\overline{\Omega})$, see \cite{arendt2011dirichlet}. $\gamma_{0}(u)$ is called the trace of $u$ and one can also use the notation $\left.u\right|_{\Gamma}$ for $u\in H^{1}(\Omega)$ (to simplify, we note $u$ instead of $\left.u\right|_{\Gamma}$).
		\par
		In the sequel, we will employ the following $H^{1}(\Omega)-$trace estimate
		\begin{eqnarray}
			\int_{\Gamma}|u|^{2}d\sigma\leqslant C \|u\|_{H^{1}(\Omega)}\|u\|_{L^{2}(\Omega)}, \label{Trace estimate}
		\end{eqnarray}
		where $C>0$ depending only on $\Omega$.
		For the proof of the inequality \eqref{Trace estimate}, we refer to \cite[Theorem. 1.5.1.10]{grisvard1985elliptic}.
		\par
		Here, we use the following weak definitions of Laplacien and the normal derivative. Let $u \in H^1(\Omega)$, we say that $\Delta u\in L^{2}(\Omega)$ if there exists a function $g \in L^2(\Omega)$ such that
		\begin{equation}
			\int_{\Omega} \nabla u \cdot \nabla v \mathrm{~d}x=-\int_{\Omega} g\,v \mathrm{~d}x, \mbox{ for all } v \in \mathcal{D}(\Omega). \label{lw}
		\end{equation}
		In this case, the function $g \in L^2(\Omega)$ verifying \eqref{lw} is unique, we denote $g$ by $\Delta u$.\\
		Let $u \in H^1(\Omega)$ which satisfy $\Delta u \in L^2(\Omega)$, we say that $u$ has a weak normal derivative if there exists a function $h \in L^2(\Gamma)$ such that
		\begin{equation}
			\int_{\Omega} \Delta u\,v \mathrm{~d} x+\int_{\Omega} \nabla u \cdot \nabla v \mathrm{~d}x=\int_{\Gamma} h\,v \mathrm{~d}\sigma, \mbox{for all } v \in H^1(\Omega). \label{dn}
		\end{equation}
		In this case, the function $h \in L^2(\Gamma)$ verifying \eqref{dn} is unique, we denote $h$ by $\partial_{\n} u$.\\
		This means that we define the normal derivative $\partial_{\n} u$ of $u$ by the validity of Green's formula.
		\subsection{Weak solutions of system $\mathcal{S}(\Omega_{0},t_{1},t_{2},f,g,\varepsilon,\mathbf{B})$.}
		The Lions' theorem \cite{Lions, showalter2013monotone} provides a significant framework for establishing the existence and uniqueness of weak solutions for  $\mathcal{S}(\Omega_{0},t_{1},t_{2},f,g,\varepsilon,\mathbf{B})$.
		\par
		Considering the data 
		\begin{eqnarray}
			\hat{f}(\cdot,t)=-f(\cdot ,t_{2}-t),\;\; \hat{\mathbf{B}}(\cdot,t)=-\mathbf{B}(\cdot,t_{2}-t)\;\;\text{and}\;\; \tau=t_{2}-t_{1}. \label{data}
		\end{eqnarray}
		One can pass from system
		$\mathcal{S}(\Omega_{0},t_{1},t_{2},f,g,\varepsilon,\mathbf{B})$ to system $\mathcal{S}(\Omega_{0},0,\tau,\hat{f},g,-\varepsilon,\hat{\mathbf{B}})$, and vice versa, by means of the transformation $t\mapsto t_{2}-t$:
		\begin{equation} \label{s4}
			\left\{
			\begin{aligned}
				\partial_t \varphi-\varepsilon\Delta \varphi - \nabla\cdot \left(\varphi \hat{\mathbf{B}}(x,t)\right) &=\hat{f}(x,t)\quad &&\text { in } \mathcal{U}\times (0,\tau), \\
				\left(\varepsilon \nabla\varphi +\varphi \hat{\mathbf{B}}(x,t)\right)\cdot \mathbf{n}(x)\mathbf{1}_{\Gamma}(x)+\varphi\mathbf{1}_{\Gamma_{0}}(x) &=0 \quad && \text { on } \partial\mathcal{U}\times (0,\tau), \\
				\varphi(x,0) &=g(x) \quad &&\text { in } \mathcal{U}.
			\end{aligned}
			\right.
		\end{equation}
		\par 
		We consider the bilinear form defined on $[0,\tau]\times H^{1}_{\Gamma_{0}}(\mathcal{U})\times H^{1}_{\Gamma_{0}}(\mathcal{U})$ by
		\begin{eqnarray}
			\mathfrak{a}_{w}(t,u,v):=\varepsilon\int_{\mathcal{U}}\nabla u\cdot\nabla v\mathrm{~d}x+\int_{\mathcal{U}}u\hat{\mathbf{B}}(x,t)\cdot\nabla v\mathrm{~d}x.
		\end{eqnarray}
		\begin{definition}
			Let $g\in L^{2}(\mathcal{U})$ and $f\in L^{2}(t_{1},t_{2};\left(H_{\Gamma_{0}}^{1}(\Omega)\right)^{'})$. A weak solution of \eqref{s3} is a function $u\in L^{2}(t_{1},t_{2};H_{\Gamma_{0}}^{1}(\mathcal{U}))\cap H^{1}\left(t_{1},t_{2}; \left(H_{\Gamma_{0}}^{1}(\mathcal{U})\right)^{'}\right)$ such that 
			\begin{eqnarray} \label{weak solution}
				&&-\int_{t_{1}}^{t_{2}}(u(t),v^{'}(t)) \d t-\int_{t_{1}}^{t_{2}}\mathfrak{a}_{w}(t_{2}-t,u(t),v(t))\d t\nonumber\\
				&& \quad\quad =\int_{t_{1}}^{t_{2}}\langle f(\cdot,t),v(t) \rangle_{(H_{\Gamma_{0}}^{1}(\mathcal{U}))^{'},H_{\Gamma_{0}}^{1}(\mathcal{U})}\d t + (g,v(t_{2})), 
			\end{eqnarray}
			for all $v\in H^{1}(t_{1},t_{2};L^{2}(\mathcal{U}))\cap L^{2}(t_{1},t_{2};H^{1}_{\Gamma_{0}}(\mathcal{U}))$ and $v(t_{1})=0$.
		\end{definition}
		\begin{proposition} \label{pw} Let $\mathbf{B}\in L^{\infty}(\mathcal{U}\times (0,T))^{d}$, then for all $0\leqslant t_{1}< t_{2}\leqslant T$, the system $\mathcal{S}(\Omega_{0},0,\tau,\hat{f},g,-\varepsilon,\hat{\mathbf{B}})$, and hence the system $\mathcal{S}(\Omega_{0},t_{1},t_{2},f,g,\varepsilon,\mathbf{B})$ has a unique weak solution. Moreover there exists a constant $C>0$ independent of $\varepsilon$ such that, the weak solution of $\mathcal{S}(\Omega_{0},t_{1},t_{2},f,g,\varepsilon,\mathbf{B})$ verifies
			\begin{eqnarray}\label{dissw}
				&&\|\varphi\|_{\mathcal{C}([t_{1},t_{2}];L^{2}(\mathcal{U}))}	+\sqrt{\varepsilon}\|\varphi\|_{L^{2}(t_{1},t_{2};H^{1}(\mathcal{U}))}\\&&\quad \leqslant C\exp\left(C(t_{2}-t_{1})\left(\frac{\|\mathbf{B}\|^{2}_{L^{\infty}(\Omega_{T})}}{\varepsilon}+\varepsilon+1\right)\right)\left(\|f\|_{L^{2}(t_{1},t_{2};L^{2}(\mathcal{U}))}+\|g\|_{L^{2}(\mathcal{U})}\right).\nonumber
			\end{eqnarray}
		\end{proposition}
		\begin{proof}
			To prove the existence and uniqueness of a weak solution of $\mathcal{S}(\Omega_{0},0,\tau,\hat{f},g,-\varepsilon,\hat{\mathbf{B}})$, we apply Lions' theorem, so it suffices to check that:
			\begin{itemize}
				\item $t\longmapsto \mathfrak{a}_{w}(t,u,v)$ is measurable for all $u,v\in H_{\Gamma_{0}}^{1}(\Omega)$;
				\item $\mathfrak{a}_{w}$ is $H_{\Gamma_{0}}^{1}(\mathcal{U})$-bounded;
				\item $\mathfrak{a}_{w}$ is quasi-coercive; i.e., there exist $\alpha> 0$ and $\omega\geq 0$ such that
				\begin{eqnarray}
					\mathfrak{a}_{w}(t,u,u) + \omega \|u\|_{L^{2}(\mathcal{U})}\geq  \alpha\|u\|^{2}_{H^{1}(\mathcal{U})}\;\mbox{for all}\; u\in H^{1}_{\Gamma_{0}}(\mathcal{U}), t\in [0,\tau]. \label{quasi-coercive}
				\end{eqnarray}
			\end{itemize}
			Using the boundedness of $\mathbf{B}$, we obtain $(x,t)\longmapsto \varepsilon\nabla u(x)\cdot\nabla v(x)+u\hat{\B}(x,t)\cdot\nabla v (x)$ is integrable on $\Omega\times [0,T]$, for all $u,v\in H^{1}(\mathcal{U})$, then, in particular from Fubini's theorem, we obtain 
			$t\longmapsto \mathfrak{a}_{w}(t,u,v)$ is measurable for all $u,v\in H_{\Gamma_{0}}^{1}(\Omega)$.
			On the other hand 
			\begin{eqnarray}
				|\mathfrak{a}_{w}(t,u,v)| &\leqslant & \varepsilon \|\nabla u\|_{L^{2}(\mathcal{U})}\|\nabla v\|_{L^{2}(\mathcal{U})}+\|\mathbf{B}\|_{L^\infty(\mathcal{U}\times (0,T))}\| u\|_{L^{2}(\mathcal{U})}\|\nabla v\|_{L^{2}(\mathcal{U})} \nonumber\\
				&\leqslant & \left(\varepsilon+\|\mathbf{B}\|_{L^\infty(\mathcal{U})\times (0,T)}\right)\| u\|_{H^{1}(\mathcal{U})}\| v\|_{H^{1}(\mathcal{U})}. \nonumber
			\end{eqnarray}
			Hence, the form $\mathfrak{a}_{w}$ is $H_{\Gamma_{0}}^{1}(\mathcal{U})$-bounded. We claim that $\mathfrak{a}_{w}$ is quasi-coercive. By Hölder's inequality, we get
			\begin{eqnarray*}
				\left|\int_{\mathcal{U}}u\hat{\mathbf{B}}(x,t)\cdot\nabla u\mathrm{~d}x\right| &\leqslant & \|\mathbf{B}\|_{L^\infty(\mathcal{U}\times (0,T))}\| u\|_{L^{2}(\mathcal{U})}\|\nabla u\|_{L^{2}(\mathcal{U})}\\
				&\leqslant & \frac{\varepsilon}{2}\|\nabla u\|^{2}_{L^{2}(\mathcal{U})}+\frac{\|\mathbf{B}\|^{2}_{L^\infty(\mathcal{U}\times (0,T))}}{2\varepsilon}\| u\|^{2}_{L^{2}(\mathcal{U})}.
			\end{eqnarray*}
			Then 
			\begin{eqnarray*}
				\mathfrak{a}_{w}(t,u,u) &\geq &  \frac{\varepsilon}{2}\| u\|^{2}_{H^{1}(\mathcal{U})}-\left(\frac{\varepsilon}{2}
				+\frac{\|\mathbf{B}\|^{2}_{L^\infty(\mathcal{U}\times (0,T))}}{2\varepsilon}\right)\| u\|^{2}_{L^{2}(\mathcal{U})}.
			\end{eqnarray*}
			Lions’ theorem and \cite[Proposition III.2.1]{showalter2013monotone} yields that  $\mathcal{S}(\Omega_{0},0,\tau,\hat{f},g,-\varepsilon,\hat{\mathbf{B}})$ has a unique weak solution.
			\par 
			Let $\varphi$ a weak solution of \eqref{s3}. 
			From \cite[Proposition III.1.2]{showalter2013monotone}, $\|\varphi(\cdot)\|^{2}_{L^{2}(\mathcal{U})}$ is absolutely continuous on $[0,T]$ and we have the standard energy identity:
			\begin{eqnarray}
				\frac{1}{2}\frac{\d}{\d t}\|\varphi\|^{2}_{L^{2}(\Omega)}=\langle \partial_{t}\varphi,\varphi \rangle_{(H_{\Gamma_{0}}^{1}(\mathcal{U}))^{'},H_{\Gamma_{0}}^{1}(\mathcal{U})}\;\;a.e\; t\in [0,T]. \label{energy identity}
			\end{eqnarray}
			Using \eqref{energy identity}
			and integrating by parts, we obtain
			\begin{eqnarray*}
				-\frac{1}{2}\frac{d}{\mathrm{~d}t}\int_{\mathcal{U}}|\varphi|^{2}\mathrm{~d}x	+\varepsilon\int_{\mathcal{U}}|\nabla\varphi|^{2}\mathrm{~d}x &=&- \int_{\mathcal{U}}\varphi \mathbf{B}(x,t)\cdot\nabla\varphi\mathrm{~d}x+\int_{\mathcal{U}}f\varphi\mathrm{~d}x\\
				&\leqslant & \|\mathbf{B}\|_{L^{\infty}(\Omega_{T})}\|\varphi\|_{L^{2}(\mathcal{U})}\|\nabla\varphi\|_{L^{2}(\mathcal{U})}+\int_{\mathcal{U}}f\varphi\mathrm{~d}x.
			\end{eqnarray*}
			By Young's inequality, we get
			\begin{eqnarray*}
				-\frac{1}{2}\frac{d}{\mathrm{~d}t}\int_{\mathcal{U}}|\varphi|^{2}\mathrm{~d}x	+\varepsilon\int_{\mathcal{U}}|\nabla\varphi|^{2}\mathrm{~d}x&\leqslant& \frac{\|\mathbf{B}\|^{2}_{L^{\infty}(\Omega_{T})}\|\varphi\|^{2}_{L^{2}(\mathcal{U})}}{2\varepsilon}+\frac{\varepsilon}{2}\|\nabla\varphi\|^{2}_{L^{2}(\mathcal{U})}\\
				&&+\int_{\mathcal{U}}f\varphi\mathrm{~d}x.
			\end{eqnarray*}
			Adding $\frac{\varepsilon}{2}\|\varphi\|_{L^{2}(\Omega)}$ in both side and by Young's inequality, we deduce that 
			\begin{eqnarray*}
				-\frac{1}{2}\frac{d}{\mathrm{~d}t}\int_{\mathcal{U}}|\varphi|^{2}\mathrm{~d}x	+\frac{\varepsilon}{2}\|\varphi\|^{2}_{H^{1}(\mathcal{U})}&\leqslant&\left(\frac{\|\mathbf{B}\|^{2}_{L^{\infty}(\Omega_{T})}}{2\varepsilon}+\frac{\varepsilon}{2}+\frac{1}{2}\right)\|\varphi\|^{2}_{L^{2}(\mathcal{U})}+\frac{1}{2}\int_{\mathcal{U}}|f|^{2}\mathrm{~d}x.
			\end{eqnarray*}
			Integrating this inequality on $[t,t_{2}]$, we obtain
			\begin{eqnarray*}
				\int_{\mathcal{U}}|\varphi(x,t)|^{2}\mathrm{~d}x	+\varepsilon\|\varphi\|^{2}_{L^{2}(t,t_{2};H^{1}(\mathcal{U}))}&\leqslant&\left(\frac{\|\mathbf{B}\|^{2}_{L^{\infty}(\Omega_{T})}}{\varepsilon}+\varepsilon+1\right)\int_{t}^{t_{2}}\int_{\mathcal{U}}|\varphi(x,s )|^{2}\mathrm{~d}x\mathrm{~d}s\\&&+\int_{t}^{t_{2}}\int_{\mathcal{U}}|f|^{2}\mathrm{~d}x\mathrm{~d}s+\int_{\mathcal{U}}|g|^{2}\mathrm{~d}x\\
				&\leqslant&\left(\frac{\|\mathbf{B}\|^{2}_{L^{\infty}(\Omega_{T})}}{\varepsilon}+\varepsilon+1\right)\int_{t}^{t_{2}}\int_{\mathcal{U}}|\varphi(x,s )|^{2}\mathrm{~d}x\mathrm{~d}s\\&&+\left(\|f\|^{2}_{L^{2}(t_{1},t_{2};L^{2}(\mathcal{U}))}+\|g\|^{2}_{L^{2}(\mathcal{U})}\right).
			\end{eqnarray*}
			The Gronwall's lemma gives the desired result.
		\end{proof}
		The following result gives an important estimate of the solutions of $\mathcal{S}(\Omega_{0},t_{1},t_{2},f,0,\varepsilon,\mathbf{B})$ for a particular source $f$.
		\begin{proposition}\label{P0}
			Let $f_{0}\in L^{2}(0,T;L^{2}(\mathcal{U}))$,  $f_{1},\cdots,f_{d}\in L^{2}(0,T;H^{1}_{0}(\mathcal{U}))$ and\\ $f=f_{0}+\varepsilon\displaystyle\sum_{i=1}^{d}\partial_{x_{i}}f_{i}$ and ssume that $\mathbf{B}(x,t)\cdot\mathbf{n}(x)=0$ on $\Gamma_{T}$. There exists $C>0$ depending on $\mathbf{B}$, $T$, $d$ and $\Omega$ such that for any $0\leqslant t_1 \leqslant t_2 \leqslant T$ and any $\varepsilon\in (0,1)$, the
			solution $\varphi$ of $\mathcal{S}(\Omega_{0},t_{1},t_{2},f,0,\varepsilon,\mathbf{B})$:
			\begin{equation}
				\left\{
				\begin{aligned}
					\partial_t \varphi+\varepsilon\Delta \varphi+\nabla \cdot\left(\varphi \mathbf{B}(x,t)\right) &=f(x,t) & & \text { in } \mathcal{U}\times (t_{1},t_{2}), \\
					\varepsilon\partial_{\mathbf{n}}\varphi\mathbf{1}_{\Gamma}(x)+\varphi\mathbf{1}_{\Gamma_{0}}(x)&=0  & & \text{ on } \partial\mathcal{U}\times (t_{1},t_{2}), \\
					\varphi(x,t_{2}) &=0 & & \text { in } \mathcal{U},
					\nonumber
				\end{aligned}
				\right.
			\end{equation}
			satisfies
			\begin{eqnarray}
				\|\varphi(\cdot,t_{1})\|^{2}_{L^{2}(\mathcal{U})}\leqslant C\sum_{i=0}^{d}\|f_{i}\|^{2}_{L^{2}(t_{1},t_{2};L^{2}(\mathcal{U}))}.
			\end{eqnarray}
		\end{proposition}
		\begin{proof}
			Using the energie identity \eqref{energy identity}, $\mathbf{B}(x,t)\cdot\mathbf{n}(x)=0$ on $\Gamma_{T}$, $\partial_{\mathbf{n}}\varphi=0$ on $\Gamma$, $\varphi=0$ on $\Gamma_{0}$ and integration by parts, we obtain 
			\begin{eqnarray}
				&&\frac{1}{2}\frac{d}{dt}\|\varphi(\cdot,t)\|^{2}_{L^{2}(\mathcal{U})}=\varepsilon\int_{\mathcal{U}}|\nabla\varphi(x,t)|^{2}\mathrm{~d}x-\frac{1}{2}\int_{\mathcal{U}}\nabla\cdot\mathbf{B}(x,t)|\varphi(x,t)|^{2}\mathrm{~d}x \nonumber\\
				&& \quad+ \int_{\mathcal{U}}f_{0}(x,t)\varphi(x,t)\mathrm{~d}x+\varepsilon\sum_{i=0}^{d}\int_{\mathcal{U}}\partial_{x_{i}}f_{i}(x,t)\varphi(x,t)\mathrm{~d}x. \label{i1}
			\end{eqnarray}
			On the other hand, since $f_{i}(\cdot,t)\in H^{1}_{0}(\mathcal{U})$ for $i=1,\cdots,d$, by integration by parts, we have
			\begin{eqnarray*}
				\int_{\mathcal{U}}\partial_{x_{i}}f_{i}(x,t)\varphi(x,t)\mathrm{~d}x=-\int_{\mathcal{U}}f_{i}(x,t)\partial_{x_{i}}\varphi(x,t)\mathrm{~d}x.
			\end{eqnarray*}
			Thus, by Cauchy Schwartz and Hölder inequalities, $\varepsilon\in (0,1)$, we get 
			\begin{eqnarray}
				\varepsilon\left|\int_{\mathcal{U}}\partial_{x_{i}}f_{i}(x,t)\varphi(x,t)\mathrm{~d}x\right|&\leqslant & \varepsilon \|f_{i}(\cdot,t)\|_{L^{2}(\mathcal{U})}\|\nabla\varphi (\cdot,t)\|_{L^{2}(\mathcal{U})} \nonumber\\
				&\leqslant & \left(\sqrt{\frac{\varepsilon}{d}}\|\nabla\varphi (\cdot,t)\|_{L^{2}(\mathcal{U})}\right)\left(\sqrt{d}\|f_{i}(\cdot,t)\|_{L^{2}(\mathcal{U})}\right) \nonumber\\
				&\leqslant & \frac{\varepsilon}{2d}\|\nabla\varphi (\cdot,t)\|^{2}_{L^{2}(\mathcal{U})}+\frac{d}{2}\|f_{i} (\cdot,t)\|^{2}_{L^{2}(\mathcal{U})} \label{i2}
			\end{eqnarray}
			and 
			\begin{eqnarray}
				\left|\int_{\mathcal{U}}f_{0}(x,t)\varphi(x,t)\mathrm{~d}x\right| 
				&\leqslant & \frac{1}{2}\|f_{0} (\cdot,t)\|^{2}_{L^{2}(\mathcal{U})}+\frac{1}{2}\|\varphi (\cdot,t)\|^{2}_{L^{2}(\mathcal{U})} \nonumber\\
				&\leqslant & \frac{d}{2}\|f_{0} (\cdot,t)\|^{2}_{L^{2}(\mathcal{U})}+\frac{1}{2}\|\varphi (\cdot,t)\|^{2}_{L^{2}(\mathcal{U})}, \label{i3}
			\end{eqnarray}
			From \eqref{i1}, \eqref{i2} and \eqref{i3}, we obtain 
			\begin{eqnarray}
				\frac{1}{2}\frac{d}{dt}\|\varphi(\cdot,t)\|^{2}_{L^{2}(\mathcal{U})} &\geq& \frac{\varepsilon}{2}\int_{\mathcal{U}}|\nabla\varphi(x,t)|^{2}\mathrm{~d}x-\frac{d}{2}\sum_{i=0}^{d}\|f_{i}(\cdot,t)\|^{2}_{L^{2}(\mathcal{U})} \nonumber\\
				&& \quad -C\|\varphi(\cdot,t)\|^{2}_{L^{2}(\mathcal{U})}\nonumber\\
				&&\geq-\frac{d}{2}\sum_{i=0}^{d}\|f_{i}(\cdot,t)\|^{2}_{L^{2}(\mathcal{U})} -C\|\varphi(\cdot,t)\|^{2}_{L^{2}(\mathcal{U})}, \label{i4}
			\end{eqnarray}
			where $C:=\|\nabla\cdot\mathbf{B}\|_{L^{\infty}(\Omega\times (0,T))}+1$.\\ Integrating \eqref{i4} on $(\tau,t_{2})$ for $t_{1}\leqslant \tau < t_{2}$, we have
			\begin{eqnarray}
				\|\varphi(\cdot,\tau)\|^{2}_{L^{2}(\mathcal{U})}&\leqslant&  d\sum_{i=0}^{d}\int_{\tau}^{t_{2}}\|f_{i}(\cdot,t)\|^{2}_{L^{2}(\mathcal{U})}\d t +2C\int_{\tau}^{t_{2}}\|\varphi(\cdot,t)\|^{2}_{L^{2}(\mathcal{U})}\mathrm{~d}t. \nonumber  
			\end{eqnarray}
			Applying Gronwall's lemma, we obtain the desired result.
		\end{proof}
		\subsection{Strong solutions of system $\mathcal{S}(\Omega_{0},t_{1},t_{2},f,g,\varepsilon,\mathbf{B})$.}
		The existence and uniqueness of strong solutions for system \eqref{s3} is derived mainly from the reference \cite{arendt2014maximal}. In this section, we will assume that $\B\in W^{1,\infty}(\Omega_{T})^{d}$, allowing us to write \eqref{s3} as
		\begin{equation} \label{s5}
			\left\{
			\begin{aligned}
				\partial_t \varphi-\varepsilon\Delta \varphi - \B (x,t)\cdot\nabla\varphi-\left(\nabla\cdot\B (x,t)\right)\varphi &=f(x,t)\quad &&\text { in } \mathcal{U}\times (t_{1},t_{2}), \\
				\left(\varepsilon \nabla\varphi +\varphi \B (x,t)\right)\cdot \mathbf{n}(x)\mathbf{1}_{\Gamma}(x)+\varphi\mathbf{1}_{\Gamma_{0}}(x) &=0 \quad && \text { on } \partial\mathcal{U}\times (t_{1},t_{2}), \\
				\varphi(x,0) &=g(x) \quad &&\text { in } \mathcal{U}.
			\end{aligned}
			\right.
		\end{equation}
		This system is equivalent to the following
		\begin{equation} \label{s6}
			\left\{
			\begin{aligned}
				\partial_t \varphi-\varepsilon\Delta \varphi - \hat{\mathbf{B}}(x,t)\cdot\nabla\varphi-\left(\nabla\cdot\hat{\mathbf{B}}(x,t)\right)\varphi &=\hat{f}(x,t)\quad &&\text { in } \mathcal{U}\times (0,\tau), \\
				\left(\varepsilon \nabla\varphi +\varphi \hat{\mathbf{B}}(x,t)\right)\cdot \mathbf{n}(x)\mathbf{1}_{\Gamma}(x)+\varphi\mathbf{1}_{\Gamma_{0}}(x) &=0 \quad && \text { on } \partial\mathcal{U}\times (0,\tau), \\
				\varphi(x,0) &=g(x) \quad &&\text { in } \mathcal{U},
			\end{aligned}
			\right.
		\end{equation}
		where $\hat{\mathbf{B}}, \hat{\mathbf{f}}$ and $\tau$ are defined in \eqref{data}.
		\par 
		We consider the bilinear form defined on $[0,\tau]\times H_{\Gamma_{0}}^{1}(\mathcal{U})\times H_{\Gamma_{0}}^{1}(\mathcal{U})$ by
		\begin{eqnarray*}
			&&\mathfrak{a}(t,u,v):=\varepsilon\int_{\mathcal{U}}\nabla u\cdot\nabla v\mathrm{~d}x+\int_{\Gamma}\left(\hat{\mathbf{B}}(x,t)\cdot\n(x)\right)u\; v\mathrm{~d}\sigma\\
			&& \quad\quad -
			\int_{\mathcal{U}}\left(\hat{\mathbf{B}}(x,t)\cdot\nabla u\right)v\mathrm{~d}x
			-\int_{\mathcal{U}}\left(\nabla\cdot\hat{\mathbf{B}}(x,t)\right)u\,v\mathrm{~d}x.
		\end{eqnarray*}
		and the following maximal regularity space
		$$MR_{\mathfrak{a}}(t_{1},t_{2}):=\{u\in H^{1}(t_{1},t_{2};L^{2}(\mathcal{U}))\cap L^{2}(t_{1},t_{2};H_{\Gamma_{0}}^{1}(\mathcal{U})):\;\; \mathcal{A}(\cdot)u(\cdot)\in L^{2}(t_{1},t_{2};L^{2}(\mathcal{U}))\},$$
		where $\mathcal{A}(t)\in \mathcal{L}\left(H_{\Gamma_{0}}^{1}(\mathcal{U}), \left(H_{\Gamma_{0}}^{1}(\mathcal{U})\right)^{'}\right)$ is the operator associated with $\mathfrak{a}(t,\cdot,\cdot)$ and defined by 
		$$\langle \mathcal{A}(t)u,v\rangle_{H_{\Gamma_{0}}^{1}(\mathcal{U})^{'},H_{\Gamma_{0}}^{1}(\mathcal{U})} :=\mathfrak{a}(t,u,v).$$
		It is a Hilbert space for the norm $\|\cdot\|_{MR_{\mathfrak{a}}(t_{1},t_{2})}$ defined by
		$$\|u\|^{2}_{MR_{\mathfrak{a}}}=\|u\|^{2}_{L^{2}(t_{1},t_{2};H^{1}(\mathcal{U}))}+\|\partial_{t}u \|^{2}_{L^{2}(t_{1},t_{2};L^{2}(\mathcal{U}))}+\|\mathcal{A}(\cdot)u(\cdot)\|^{2}_{L^{2}(t_{1},t_{2};L^{2}(\mathcal{U}))}.$$
		We have the follawing important result:
		\begin{proposition} \label{embedded}
			We have the follawing injections:
			\begin{enumerate}[label=\textcolor{red}{(\arabic*)}]
				\item The space $MR_{\mathfrak{a}}(t_{1},t_{2})$ embeds continuously into $\mathcal{C}([t_{1},t_{2}]; H_{\Gamma_{0}}^{1}(\mathcal{U}))$.
				\item The space $MR_{\mathfrak{a}}(t_{1},t_{2})$  embeds compactly into $L^{2}(t_{1},t_{2};H_{\Gamma_{0}}^{1}(\mathcal{U}))$.
			\end{enumerate}
		\end{proposition}
		\begin{proof}
			For more details, we refer to
			\cite[Corollary 3.3]{arendt2014maximal}.
		\end{proof}
		For all $t\in [0,\tau]$, we define the operators $A_{1}(t)$ and $A_{2}(t)$ by
		\begin{eqnarray*}
			&&D(A_{1}(t)):=\{u\in H_{\Gamma_{0}}^{1}(\mathcal{U}):\;\Delta u\in L^{2}(\mathcal{U}),\; \varepsilon\partial_{\n} u +\left(\hat{\mathbf{B}}(x,t)\cdot\n(x)\right)u=0\},\\
			&&D(A_{2}(t)):= H_{\Gamma_{0}}^{1}(\mathcal{U}),
		\end{eqnarray*} 
		for all $u\in D(A_{1}(t))$ and $v\in H_{\Gamma_{0}}^{1}(\mathcal{U})$ the operators $A_{1}(t)$ and $A_{2}(t)$ are given by $A_{1}(t)u:=-\varepsilon\Delta u$ and $A_{2}(t) v:=-\hat{\mathbf{B}}(x,t)\cdot\nabla v-(\nabla\cdot\hat{\mathbf{B}}(x,t))v$.\\
		The system \eqref{s5} can be written equivalently as a Cauchy initial valued problem
		\begin{equation}\label{c1}
			\left\{
			\begin{aligned}
				Y^{'}+A(t)Y &=F(t) & & t\in [0,\tau], \\
				Y(0)&=g, 
			\end{aligned}
			\right.
		\end{equation}
		where $A(t)=A_{1}(t)+A_{2}(t)$, $D(A(t))=D(A_{1}(t))$, $Y(t)=\varphi(\cdot,t)$ and $F(t)=\hat{f}(\cdot, t)$.\\
		We start with the definition of a strong solution of \eqref{s5}.
		\begin{definition}
			Let $g\in L^{2}(\mathcal{U})$ and $f\in L^{2}(t_{1},t_{2};L^{2}(\mathcal{U}))$. A strong solution of \eqref{s3} is a function $u\in MR_{\mathfrak{a}}(t_{1},t_{2})$ fulfilling $\eqref{s5}_{1}$ in $L^{2}(t_{1},t_{2};L^{2}(\mathcal{U}))$, $\eqref{s5}_{2}$ in $L^{2}(t_{1},t_{2};L^{2}(\partial\mathcal{U}))$ and $\eqref{s5}_{3}$, where $\eqref{s5}_{j}$ is the j-th equation in system \eqref{s5}.
		\end{definition}
		No we are in position to establish the following existence, uniqueness and regularity results.
		\begin{proposition} \label{ps}
			Let $\B\in W^{1,\infty}(\Omega_{T})^{d}$, $g\in H_{\Gamma_{0}}^{1}(\mathcal{U})$ and $f\in L^{2}(t_{1},t_{2};L^{2}(\mathcal{U}))$. Then the Cauchy problem \eqref{c1}, and hence the system \eqref{s3} has a unique strong 
			solution $u\in MR_{\mathfrak{a}}(t_{1},t_{2})$.\\ Moreover, $u\in \mathcal{C}([t_{1},t_{2}];L^{2}(\mathcal{U}))$ and there exists a constant $C:=C(T,\varepsilon)\geq 0$ such that
			\begin{equation}
				\|u\|_{MR_{\mathfrak{a}}(t_{1},t_{2})}\leqslant C\left(\|g\|_{L^{2}(\mathcal{U})}+\|f\|_{L^{2}(t_{1},t_{2};L^{2}(\mathcal{U}))}\right).
			\end{equation}
		\end{proposition}
		\begin{proof} 
			To prove the existence and uniqueness of a strong solution of \eqref{s1}, we apply \cite[Theorem 4.2 and Remark 4.6]{arendt2014maximal}, so 
			we consider the bilinear forms defined on $[0,\tau]\times H_{\Gamma_{0}}^{1}(\mathcal{U})\times H_{\Gamma_{0}}^{1}(\mathcal{U})$ by 
			\begin{eqnarray}
				&&\mathfrak{a}_{1}(t,u,v):=\varepsilon\int_{\mathcal{U}}\nabla u\cdot\nabla v\mathrm{~d}x+\int_{\Gamma}\left(
				\hat{\mathbf{B}}(x,t)\cdot\n(x)\right)u\;v\mathrm{~d}\sigma, \\
				&&\mathfrak{a}_{2}(t,u,v):= -\int_{\mathcal{U}}\left(\hat{\mathbf{B}}(x,t)\cdot\nabla u\right)v\mathrm{~d}x-\int_{\mathcal{U}}\left(\nabla\cdot\hat{\mathbf{B}}(x,t)\right)u\;v\mathrm{~d}x.
			\end{eqnarray}
			Clearly, we have 
			\begin{eqnarray*}
				\mathfrak{a}_{1}(t,u,v)+\mathfrak{a}_{2}(t,u,v)=\mathfrak{a}(t,u,v)
			\end{eqnarray*}
			and we claim that, $\mathfrak{a}_{1}$ and $\mathfrak{a}_{2}$ satisfies the conditions:
			\begin{itemize}
				\item $|\mathfrak{a}_{1}(t,u,v)|\leq M_{1}\|u\|_{H^{1}(\mathcal{U})}\|v\|_{H^{1}(\mathcal{U})}\;\mbox{for all}\; u,v\in H^{1}_{\Gamma_{0}}(\mathcal{U})$, $t\in [0,\tau]$;
				\item $\mathfrak{a}_{1}$ is quasi-coercive, see \eqref{quasi-coercive};
				\item  $\mathfrak{a}_{1}$ satisfies the square root property; i.e., $R\left(A_{1}(t)^{-\frac{1}{2}}\right)=H_{\Gamma_{0}}^{1}(\mathcal{U})$;
				\item $\mathfrak{a}_{1}$ is Lipschitz-continuous; i.e., there exists a constant $C_{1}\geq 0$ such that, for all $u,v\in H^{1}_{\Gamma_{0}}(\mathcal{U})$ and $s,t\in [0,\tau]$,
				\begin{eqnarray*}
					|\mathfrak{a}_{1}(t,u,v)-\mathfrak{a}_{1}(s,u,v)|\leq C_{1}|t-s|\|u\|_{H^{1}(\mathcal{U})}\|v\|_{H^{1}(\mathcal{U})};
				\end{eqnarray*}
				\item $|\mathfrak{a}_{2}(t,u,v)|\leq M_{2}\|u\|_{H^{1}(\mathcal{U})}\|v\|_{L^{2}(\mathcal{U})}\;\mbox{for all}\; u\in H^{1}_{\Gamma_{0}}(\mathcal{U}), v\in L^{2}(\mathcal{U})$ $t\in [0,\tau]$;
				\item $t\longmapsto \mathfrak{a}_{2}(t,u,v)$ is measurable for all $u, v\in H^{1}_{\Gamma_{0}}(\mathcal{U})$.
			\end{itemize}
			By the boundedness of $\B$ and the continuity of the trace operator, the form $\mathfrak{a}_{1}$ is $H^{1}(\mathcal{U})$-bounded. Since $\mathfrak{a}_{1}$ is symetric, then it satisfies the square root property, see \cite{kato2013perturbation}. Using the classical inequality \eqref{Trace estimate} and Young's inequality, we obtain
			\begin{eqnarray*}
				\left|\int_{\Gamma}\left(\hat{\B}(x,t)\cdot\n (x)\right)|u|^{2}\mathrm{~d}\sigma\right|\leq \frac{\varepsilon}{2}\|\nabla u\|^{2}_{L^{2}(\mathcal{U})}+ \frac{C}{\varepsilon}\| u\|^{2}_{L^{2}(\mathcal{U})}.
			\end{eqnarray*}
			Hence $\mathfrak{a}_{1}$ is quasi-coercive. By the Lipschitz-continuous of $\B$, the form $\mathfrak{a}_{1}$ is also Lipschitz-continuous. 
			The boundedness of $\B$ implies the form $\mathfrak{a}_{2}: H_{\Gamma_{0}}^{1}(\mathcal{U})\times L^{2}(\mathcal{U})\longrightarrow\mathbb{R}$ is bounded for all fixed $t\in [0,\tau]$. We also have $t\longmapsto \mathfrak{a}_{2}(t,u,v)$ is measurable for all $u,v\in H_{\Gamma_{0}}^{1}(\mathcal{U})$ as for the form $\mathfrak{a}_{w}$ above.\\
			By \cite[Theorem 4.2 and Remark 4.6]{arendt2014maximal}, the Cauchy problem
			\begin{equation}\label{ccal1}
				\left\{
				\begin{aligned}
					Y^{'}+\mathcal{A}(t)Y &=F(t) & & t\in [0,\tau], \\
					Y(0)&=y_{0}, 
				\end{aligned}
				\right.
			\end{equation}
			has a unique strong solution $Y\in H^{1}(0,\tau;L^{2}(\mathcal{U}))\cap L^{2}(0,\tau;H_{\Gamma_{0}}^{1}(\mathcal{U}))$.\\
			Let us then show that \eqref{c1} has a unique strong solution $Y\in MR_{\mathfrak{a}}$.
			That is, we will show if $Y\in MR_{\mathfrak{a}}$ the strong solution of \eqref{ccal1}, then $Y(t)\in D(A(t))\;\text{and}\; A(t)(Y(t))=\mathcal{A}(t)(Y(t))$ for all $t\in [0,\tau]$.\\
			For all $v\in H_{\Gamma_{0}}^{1}(\mathcal{U})$, the strong solution of \eqref{ccal1}, satisfies
			\begin{eqnarray*}
				\int_{\mathcal{U}}Y^{'}(t)v\mathrm{~d}x+\int_{\mathcal{U}}\mathcal{A}(t)Y(t)v\mathrm{~d}x=\int_{\mathcal{U}}F(t)v\mathrm{~d}x.
			\end{eqnarray*}
			Then 
			\begin{eqnarray*}
				\int_{\mathcal{U}}Y^{'}(t)v\mathrm{~d}x+\mathfrak{a}(t,Y(t),v)=\int_{\mathcal{U}}F(t)v\mathrm{~d}x.
			\end{eqnarray*}
			So 
			\begin{eqnarray}
				&&\int_{\mathcal{U}}Y^{'}(t)v\mathrm{~d}x+\varepsilon\int_{\mathcal{U}}\nabla Y(t)\cdot\nabla v \mathrm{~d}x+\int_{\Gamma}\hat{\mathbf{B}}(x,t)\cdot \n(x)Y(t)v\mathrm{~d}\sigma \\ 
				&&\quad -\int_{\mathcal{U}}\left(\hat{\mathbf{B}}(x,t)\cdot\nabla Y(t)\right)v\mathrm{~d}x -\int_{\mathcal{U}}\left(\nabla\cdot \hat{\mathbf{B}}(x,t) \right)Y(t)v\mathrm{~d}x=\int_{\mathcal{U}}F(t)v\mathrm{~d}x. \nonumber
				\label{cc4}
			\end{eqnarray}
			In particular, for all $v\in \mathcal{D}(\mathcal{U})$
			\begin{eqnarray*}
				\int_{\mathcal{U}}Y^{'}(t)v\mathrm{~d}x &+&\varepsilon\int_{\mathcal{U}}\nabla Y(t)\cdot\nabla v \mathrm{~d}x-\int_{\mathcal{U}}\left(\hat{\mathbf{B}}(x,t)\cdot\nabla Y(t)\right)v\mathrm{~d}x \nonumber\\
				&-&\int_{\mathcal{U}}\left(\nabla\cdot \hat{\mathbf{B}}(x,t) \right)Y(t)v\mathrm{~d}x=\int_{\mathcal{U}}F(t)v\mathrm{~d}x.
			\end{eqnarray*}
			Then
			\begin{eqnarray*}
				\int_{\mathcal{U}}  \left[Y^{'}(t)-\hat{\mathbf{B}}(x,t)\cdot\nabla Y(t)-\nabla\cdot \hat{\mathbf{B}}(x,t) Y(t)-F(t)\right]v\mathrm{~d}x
				= -\varepsilon\int_{\mathcal{U}}\nabla Y(t)\cdot\nabla v \mathrm{~d}x,
			\end{eqnarray*}
			for all $v\in D(\mathcal{U})$, therfore $\Delta Y(t)\in L^{2}(\mathcal{U})$, and 
			\begin{eqnarray*}
				\varepsilon\Delta Y(t)=Y^{'}(t)-\hat{\mathbf{B}}(x,t)\cdot\nabla Y(t)-\nabla\cdot \hat{\mathbf{B}}(x,t)Y(t)-F(t).
			\end{eqnarray*}
			Thus 
			\begin{eqnarray}\label{cc5}
				Y^{'}(t)=\varepsilon\Delta Y(t)+\hat{\mathbf{B}}(x,t)\cdot\nabla Y(t)+\nabla\cdot \hat{\mathbf{B}}(x,t) Y(t)+F(t).
			\end{eqnarray}
			By replacing \eqref{cc5} in \eqref{cc4}, we deduce that
			\begin{eqnarray*}
				\forall v\in H_{\Gamma_{0}}^{1}(\mathcal{U}),\; \varepsilon\int_{\mathcal{U}}\Delta Y(t)v\mathrm{~d}x+\varepsilon\int_{\mathcal{U}}\nabla Y(t)\cdot\nabla v \mathrm{~d}x+\int_{\Gamma}\hat{\mathbf{B}}(x,t)\cdot\n(x)Y(t)\mathrm{~d}\sigma =0.
			\end{eqnarray*}
			Then $\partial_{\n}Y(t)\in L^{2}(\Gamma)$ and $\varepsilon\partial_{\n}Y(t)=-\hat{\mathbf{B}}(x,t)\cdot\n(x)Y(t).$ Thus $Y(t)\in D(A(t))$. By a simple integration by parts, we have
			\begin{eqnarray*}
				\left( A(t)(Y(t)),v\right)=\mathfrak{a}(t,Y(t),v),
			\end{eqnarray*}
			for all  $v\in H_{\Gamma_{0}}^{1}(\mathcal{U})$, then $A(t)(Y(t))=\mathcal{A}(t)(Y(t))$. Then the Cauchy problem \eqref{c1}, and hence system \eqref{s3} has a unique strong 
			solution $u\in MR_{\mathfrak{a}}(t_{1},t_{2})$.
			From the Proposition \ref{embedded}, we obtain $u\in\mathcal{C}([t_{1},t_{2}];H_{\Gamma_{0}}^{1}(\mathcal{U}))$.
		\end{proof}
		\begin{remark}
			Propositions \ref{pw} and \ref{ps} are valid if $\Omega_{0}=\varnothing$.
		\end{remark}
		\section{Agmon's inequalities}
		\label{Section 4} In this section, we wil present somme technical results and we will prove somme new Agmon inequalities which will be the key to establish very interesting dissipativity estimates.
		\par
		Let us start with the following notation, which will be useful in what follows 
		\begin{notation}
			Let $0\leqslant t_{1}\leqslant t_{2}\leqslant T$ and $x_{0}\in\mathbb{R}^{d}$. For $r>0$, we note $\mathcal{D}_{r}(t_{1},t_{2},x_{0})$ the union of
				trajectories starting at $t_{2}$ in the ball $\overline{B}(x_{0},r)$:
			$$\mathcal{D}_{r}(x_{0},t_{1},t_{2})=\left\lbrace (\varPhi(t,t_{2},y),t)\;:\; y\in\overline{B}(x_{0},r)\;\;\mbox{and}\;\; t\in [t_{1},t_{2}]\right\rbrace,$$ 
			where $t\longmapsto \varPhi(t,\cdot,\cdot)$ the trajectories of ordinary equation \eqref{OD}.
		\end{notation}
		
		The following Lemma asserts the existence of a Lipchitzian function that verifies certain conditions associated with the trajectories of vector $\B$, the construction of this function is based on the change of coordinates by the trajectories of the vector field $\mathbf{B}$ and the use of radial functions. For more details, see \cite[Section 2.2]{guerrero2007singular}.
		
		\begin{lemma} \label{lemma}
			Let $\mathbf{B}\in L^{\infty}(0,T;W^{1,\infty}(\mathbb{R}^{d})^{d})$, then for all $0\leqslant t_{1}\leqslant t_{2}\leqslant T$, $x_{0}\in\mathbb{R}^{d}$ and $r>0$
			there exists a nonnegative Lipschitz function $\theta$ on $\mathbb{R}^{d}\times [t_{1},t_{2}]$ such that 
			\begin{eqnarray}
				\partial_{t}\theta-|\nabla \theta|^2 +\mathbf{B} \cdot \nabla \theta \geq 0 \quad \text { a.e in } \mathbb{R}^{d} \times[t_{1},t_{2}],\\
				\theta(x,t)=0\;\;\forall (x,t)\in \mathcal{D}_{r}(x_{0},t_{1},t_{2}),\\
				\theta(x,t)\geq c_{0}r^{2}\;\;\forall (x,t)\notin \mathcal{D}_{2r}(x_{0},t_{1},t_{2}),
			\end{eqnarray}
			where $c_{0}>0$ depends only on $t_{2}-t_{1}$ and $\displaystyle\int_{t_{1}}^{t_{2}}\|\nabla \mathbf{B}(\cdot, s)\|_{\infty}\mathrm{~d}s$.\\
		\end{lemma}
		Now we are ready to present and prove some Agmon inequalities. 
		\begin{proposition} \label{propagmon}
			Let $\Omega$ be a domain with Lipschitz boundary $\Gamma$, $\mathbf{B}\in L^{\infty}(0,T;W^{1,\infty}(\Omega)^{d})$
			and let $\theta$ be a Lipschitz function on $\overline{\Omega} \times[0, T]$ such that
			\begin{equation}
				\partial_{t}\theta -|\nabla \theta|^2 +\mathbf{B} \cdot \nabla \theta \geq 0,\;\;  \text { a.e in } \overline{\Omega} \times[0, T].
			\end{equation}
			Then, we have the following estimates:
			\begin{enumerate}[label=(\arabic*)]
				\item There exists a constant $C>0$ independent of $\varepsilon$ such that for all $\varepsilon\in (0,1]$ and
				any solution $\varphi$ of system  $\mathcal{S}(\Omega_{0},t_{1},t_{2},0,g,\varepsilon,\mathbf{B})$ with data $g \in L^{2}(\mathcal{U})$, the following Agmon-type inequality holds true for all $t \in[t_{1},t_{2}]$,
				\begin{eqnarray}
					\exp\left(-\frac{C}{\varepsilon}(t_{2}-t)\right)\int_{\mathcal{U}}|\psi(x, t)|^2 \mathrm{~d} x &+&\varepsilon \int_t^{t_{2}} \int_{\mathcal{U}} \exp\left(-\frac{C}{\varepsilon}(t_{2}-s)\right)|\nabla \psi(x, s)|^2 \mathrm{~d} x \mathrm{~d} s \nonumber\\
					&& \leq \int_{\mathcal{U}}|\psi(x, t_{2})|^2 \mathrm{~d}x, \label{A1}
				\end{eqnarray}
				where $ \psi=\exp\left(\frac{\theta}{\varepsilon}\right)\varphi$ and $\mathcal{U}=\Omega\setminus\overline{\Omega_{0}}$.\\
				\item If moreover, $\mathbf{B}(x,t)\cdot\mathbf{n}(x)\geq 0$ on $\Gamma_{T}$, then 
				\begin{eqnarray} \label{A2}
					\exp\left(-C_{\mathbf{B}}(t_{2}-t)\right)\int_{\mathcal{U}}|\psi(x, t)|^2 \mathrm{~d} x&+&2\varepsilon \int_t^{t_{2}} \int_{\mathcal{U}} \exp\left(-C_{\mathbf{B}}(t_{2}-s)\right)|\nabla \psi(x, s)|^2 \mathrm{~d} x \mathrm{~d} s \nonumber\\
					&& \leq \int_{\mathcal{U}}|\psi(x, t_{2})|^2 \mathrm{~d}x,
				\end{eqnarray}
				where $C_{\mathbf{B}}:=\|\nabla\cdot\mathbf{B}\|_{L^{\infty}(\Omega_{T})}$.
			\end{enumerate}
		\end{proposition}
		\begin{proof}
			Let $g\in L^{2}(\mathcal{U})$,
			$t\in [t_{1},t_{2}]$ and $E(t)=\dfrac{1}{2}\displaystyle\int_{\mathcal{U}}|\psi(x,t)|^{2}dx$. By integration by parts, we have 
			\begin{eqnarray}
				E^{'}(t)&=& \frac{1}{\varepsilon}\int_{\mathcal{U}}\left(\partial_{t}\theta-\lvert\nabla\xi\lvert^{2}+\mathbf{B}(x,t)\cdot\nabla\theta\right)|\psi|^{2}\mathrm{~d}x+\varepsilon\int_{\mathcal{U}}\lvert \nabla \psi\lvert^{2}\mathrm{~d}x \nonumber\\
				&&-\frac{1}{2}\int_{\mathcal{U}}\nabla\cdot \mathbf{B}(x,t)|\psi|^{2}\mathrm{~d}x  +\frac{1}{2}\int_{\Gamma}\mathbf{B}(x,t)\cdot\n(x)|\psi|^{2} \mathrm{~d}\sigma. \label{en} 
			\end{eqnarray}
			\begin{enumerate}[label=(\arabic*)]
				\item By the hypothesis verified by the function $\theta$, we have  
				\begin{eqnarray}
					E^{'}(t)&\geq& \varepsilon\int_{\mathcal{U}}\lvert \nabla \psi\lvert^{2}\mathrm{~d}x-\|\nabla\cdot\mathbf{B}\|_{L^{\infty}(\Omega_{T})}E(t)-\frac{\|\mathbf{B}\|_{L^{\infty}(\Gamma_{T})}}{2}\int_{\Gamma}|\psi|^{2} \mathrm{~d}\sigma. \nonumber
				\end{eqnarray}
				By the trace estimate \eqref{Trace estimate}
				and Young’s inequality, we obtain:
				\begin{eqnarray}
					E^{'}(t)&\geq& \frac{\varepsilon}{2}\int_{\mathcal{U}}\lvert \nabla \psi\lvert^{2}\mathrm{~d}x-\left(\|\nabla\cdot \mathbf{B}\|_{L^{\infty}(\Omega_{T})}+\varepsilon+\frac{C^{2}\|\mathbf{B}\|^{2}_{L^{\infty}(\Gamma_{T})}}{4\varepsilon}\right)E(t). \nonumber
				\end{eqnarray}
				Thus, there exists a constant $C>0$ independent of $\varepsilon$ such that
				\begin{eqnarray}
					E^{'}(t)&\geq& \frac{\varepsilon}{2}\int_{\mathcal{U}}\lvert \nabla \psi\lvert^{2}\mathrm{~d}x-\frac{C}{\varepsilon}E(t). \nonumber
				\end{eqnarray}	
				By applying Grönwall's lemma, we deduce the inequality \eqref{A1}.
				\item If $\mathbf{B}(x,t)\cdot\mathbf{n}(x)\geq 0$ on $\Gamma_{T}$, then \eqref{en} gives 
				\begin{eqnarray*}
					E^{'}(t)\geq \varepsilon\int_{\mathcal{U}}\lvert \nabla \psi\lvert^{2}\mathrm{~d}x-\frac{1}{2}\int_{\mathcal{U}}\nabla\cdot \mathbf{B}(x,t)|\psi|^{2}\mathrm{~d}x.
				\end{eqnarray*}
				Hence
				\begin{eqnarray*}
					E^{'}(t)\geq \varepsilon\int_{\mathcal{U}}\lvert \nabla \psi\lvert^{2}\mathrm{~d}x-C_{\mathbf{B}}E(t),
				\end{eqnarray*}
				where $C_{\mathbf{B}}=\|\nabla\cdot\mathbf{B}\|_{L^{\infty}(\Omega_{T})}$. The 
				Grönwall's lemma gives directly the inequality \eqref{A2}.
			\end{enumerate}
		\end{proof}
		Considering $\theta=0$ in the previous lemma, we obtain the following corollary:
		\begin{corollary}
			Assume that $\mathbf{B}\in L^{\infty}(0,T;W^{1,\infty}(\Omega)^{d})$ such that $\mathbf{B}(x,t)\cdot\mathbf{n}(x)\geq 0$ on $\Gamma_{T}$.
			Then, any solution $\varphi$ of system  $\mathcal{S}(\Omega_{0},t_{1},t_{2},0,g,\varepsilon,\mathbf{B})$ with data $g \in L^{2}(\mathcal{U})$  satisfies
			\begin{eqnarray} \label{A3}
				\int_{\mathcal{U}}|\varphi(x, t)|^2 \mathrm{~d} x\leq \exp\left(C_{\mathbf{B}}(t_{2}-t)\right)\int_{\mathcal{U}}|\varphi(x, t_{2})|^2 \mathrm{~d}x,
			\end{eqnarray}
			where $C_{\mathbf{B}}:=\|\nabla\cdot\mathbf{B}\|_{L^{\infty}(\Omega_{T})}$.
		\end{corollary}
		\section{Dissipation results by Agmon inequality} \label{Section 5}
		In this section, we will assume that
		\begin{enumerate}[label=\textcolor{red}{(\Alph*)}]
			\item $T>0$, $\varepsilon\in (0,1]$ and $\mathbf{B}\in L^{\infty}(0,T;W^{1,\infty}(\Omega)^{d})$, \label{A}
			\item $\mathbf{B}(x,t)\cdot\mathbf{n}(x)\geq 0$ on $\Gamma_{T}$, \label{B}
			\item $\exists T_{0}\in (0,T)$ and $r_{0}>0$ such that $(T,T_{0},r_{0},\mathbf{B},\Omega)$  satisfies the condition \eqref{Flushing Condition} for the control region $\omega$ \label{C}
		\end{enumerate}
		and we will prove some very important dissipation results.
		\par   
		Applying Proposition \ref{P1}, to Hypothesis \ref{C}, there exists $\omega_{0}\subset\subset\omega$ a regular open such that $(T,T_{0},\frac{r_{0}}{2},\mathbf{B},\Omega)$ satissfies the condition \eqref{Flushing Condition} for $\omega_{0}$. Hence
		\begin{eqnarray}
			\forall x_{0}\in\overline{\Omega},\; \forall t_{0}\in [T_{0},T],\; \exists t\in (t_{0}-T_{0},t_{0}),\; \forall x\in \overline{B}\left(x_{0},\frac{r_{0}}{2}\right),\; \varPhi(t,t_{0},x)\in \omega_{0}. \label{o4}
		\end{eqnarray} 
		Let
		\begin{equation*}
			\mathcal{U}:=\Omega\setminus\overline{\omega_{0}}.
		\end{equation*} 
		The assertion \eqref{o4}, implies that 
		\begin{eqnarray}
			\forall x_{0}\in\overline{\Omega},\; \forall t_{0}\in [T_{0},T],\; \exists t\in (t_{0}-T_{0},t_{0}),\; \forall x\in \overline{B}\left(x_{0},\frac{r_{0}}{2}\right),\; \varPhi(t,t_{0},x)\notin \overline{\mathcal{U}}. \label{o5}
		\end{eqnarray}
		\subsection{Dissipation result outside the region $\omega_{0}$}
		\begin{proposition} \label{P2}
			Under Hypotheses \ref{A}, \ref{B} and \ref{C}, there are constants $C_{0},C>0$  independent of $\varepsilon$ such that for any $t_{0}\in [T_{0},T]$ and all solution $\varphi$ of  $\mathcal{S}(\omega_{0},t_{0}-T_{0},t_{0},0,g,\varepsilon,\mathbf{B})$ with data $g\in L^{2}(\mathcal{U})$ verify the following dissipation estimates:
			\begin{equation} \label{disspation 1}
				\|\varphi(\cdot, t_{0}-T_{0})\|^{2}_{L^{2}(\mathcal{U})}\leqslant C\exp\left(\frac{-C_{0}}{\varepsilon}\right)\|g\|^{2}_{L^{2}(\mathcal{U})}.
			\end{equation}
		\end{proposition}
		\begin{proof}
			From \eqref{o5}, we have 
			\begin{eqnarray}
				\forall (x_{0},t_{0})\in\overline{\mathcal{U}}\times [T_{0},T],\; \exists t:=t(x_{0},t_{0})\in (t_{0}-T_{0},t_{0}),\; \forall x\in \overline{B}\left(x_{0},\frac{r_{0}}{2}\right),\; \varPhi(t,t_{0},x)\notin \mathcal{\overline{\mathcal{U}}}.\nonumber\\ \label{d1}
			\end{eqnarray}
			Since $\overline{\mathcal{U}}$ is compact, then it admits a finite partition by the balls $B_{j}\left(x_{j},\frac{r_{0}}{2}\right)$,\\ $j=1,\cdots,J$ and a
			partition of unity $\chi_{j}$ associated to this finite covering.\\
			For all $j=1,\cdots,J$, we consider $\theta_{j}$ the function that satisfies the Lemma \ref{lemma} with the choice $x_{0}=x_{j}$, $t_{1}=t_{0}-T_{0}$ and $t_{2}=t_{0}$ and $r=\frac{r_{0}}{2}$. Let $\varphi$ the weak solution of $\mathcal{S}(\omega_{0},t_{0}-T_{0},t_{0},0,g,\varepsilon,\mathbf{B})$ and $\varphi_{j}$ the weak solution of $\mathcal{S}(\omega_{0},t_{0}-T_{0},t_{0},0,\chi_{j}(\;.\;)g,\varepsilon,\mathbf{B})$.\\
			By Agmon inequality \eqref{A2}, for all $t\in [t_{0}-T_{0},t_{0}]$, we obtain
			\begin{equation}
				\| \psi_{j}(\cdot,t)\|_{L^{2}(\mathcal{U})}\leqslant \exp\left(\frac{C_{\mathbf{B}}}{2}(t_{0}-t)\right)	\| \psi_{j}(\cdot,t_{0})\|_{L^{2}(\mathcal{U})},
			\end{equation}
			where $\psi_{j}(\cdot, t)=\exp\left(\frac{\theta_{j}(\cdot, t)}{\varepsilon}\right)\varphi_{j}(\cdot,t)$.\\
			The properties of $\theta_{j}$ in the Lemma \ref{lemma} and \eqref{d1} give
			$$	\begin{cases}
				\theta_{j}(x,t_{0})=0,\quad &\text{if}\;x\in B_{j}\left(x_{0},\frac{r_{0}}{2}\right),\\
				\theta_{j}(x,t(x_{j},t_{0}))\geq \frac{c_{0}}{4}r_{0}^{2},\quad &\text{if}\; x\in\overline{\mathcal{U}}.
			\end{cases}$$
			Hence
			\begin{eqnarray}
				\| \varphi_{j}(\cdot,t(x_{j},t_{0}))\|_{L^{2}(\mathcal{U})}&\leqslant & \exp\left(\frac{-c_{0}r_{0}^{2}}{4\varepsilon}\right)\exp\left(\frac{C_{\mathbf{B}}}{2}(t_{0}-t(x_{j},t_{0}))\right)	\| \varphi_{j}(\cdot,t_{0})\|_{L^{2}(\mathcal{U})}\nonumber\\
				&\leqslant & \exp\left(\frac{-c_{0}r_{0}^{2}}{4\varepsilon}\right)\exp\left(\frac{C_{\mathbf{B}}}{2}T_{0}\right)	\| \varphi_{j}(\cdot,t_{0})\|_{L^{2}(\mathcal{U})}. \label{d2}
			\end{eqnarray}
			Using Agmon inequality \eqref{A3}, we have 
			\begin{eqnarray}
				\| \varphi_{j}(\cdot,t_{0}-T_{0})\|_{L^{2}(\mathcal{U})}
				&\leqslant & \exp\left(\frac{C_{\mathbf{B}}}{2}\left(t(x_{j},t_{0})-(t_{0}-T_{0})\right)\right)\| \varphi_{j}(\cdot,t(x_{j},t_{0}))\|_{L^{2}(\mathcal{U})} \nonumber\\
				&\leqslant & \exp\left(\frac{C_{\mathbf{B}}}{2}T_{0}\right)\| \varphi_{j}(\cdot,t(x_{j},t_{0}))\|_{L^{2}(\mathcal{U})}. \label{d3}
			\end{eqnarray}
			From \eqref{d2} and \eqref{d3}, we deduce that
			\begin{eqnarray} \label{d4}
				\| \varphi_{j}(\cdot,t_{0}-T_{0})\|_{L^{2}(\mathcal{U})}
				\leqslant \exp\left(\frac{-c_{0}r_{0}^{2}}{4\varepsilon}\right)\exp\left(C_{\mathbf{B}}T_{0})\right)	\| \varphi_{j}(\cdot,t_{0})\|_{L^{2}(\mathcal{U})}.
			\end{eqnarray}
			The fact that the systems considered are linear and  $\varphi(\cdot,t_{0})=\displaystyle\sum_{j=1}^{J}\chi_{j}(\cdot)g$, we get 
			\begin{equation*}
				\varphi(\cdot,t)=\sum_{j=1}^{J}\varphi_{j}(\cdot,t),\quad\text{for all}\;t_{0}-T_{0}\leqslant t\leqslant t_{0}.
			\end{equation*}
			Using this decomposition of $\varphi$ and \eqref{d4}, we then obtain
			\begin{eqnarray}
				\| \varphi(\cdot,t_{0}-T_{0})\|_{L^{2}(\mathcal{U})}
				\leqslant J\exp\left(\frac{-c_{0}r_{0}^{2}}{4\varepsilon}\right)\exp\left(C_{\mathbf{B}}T_{0}\right)	\| g\|_{L^{2}(\mathcal{U})}.
			\end{eqnarray}
			Hence the estimate \eqref{disspation 1} with the choice $C=J^{2}\exp\left(2C_{\mathbf{B}}T_{0}\right)$ and $C_{0}=\frac{c_0}{2}r_{0}^{2}$.
		\end{proof}
		\subsection{Global dissipation result}
		\begin{proposition} \label{P3}
			Under Hypotheses \ref{A}, \ref{B} and \ref{C}, there are constants $C_{0},C>0$  independent of $\varepsilon$ such that for any $\varphi$ of  $\mathcal{S}(\varnothing,0,T,0,g,\varepsilon,\mathbf{B})$ with data $g\in L^{2}(\Omega)$ verify the following dissipation estimates:
			\begin{enumerate}[label=\textcolor{red}{(\arabic*)}]
				\item For all $t_{0}\in [T_{0},T]$, we have
				\begin{equation} \label{disspation 2}
					\|\varphi(\cdot, t_{0}-T_{0})\|^{2}_{L^{2}(\Omega)}\leqslant C\left(\exp\left(\frac{-C_{0}}{\varepsilon}\right)\|\varphi(\cdot, t_{0})\|^{2}_{L^{2}(\Omega)}+\|\varphi\|^{2}_{L^{2}(t_{0}-T_{0},t_{0};L^{2}(\omega))}\right).
				\end{equation}
				\item For any integer $m$ such that $1\leqslant m\leqslant\frac{T}{T_{0}}$, there exists $C>0$  independent of $\varepsilon$ such that for all $t\in [mT_{0},T]$, we have
				\begin{equation} \label{disspation 3}
					\|\varphi(\cdot, 0)\|^{2}_{L^{2}(\Omega)}\leqslant C\left(\exp\left(\frac{-mC_{0}}{\varepsilon}\right)\|\varphi(\cdot, t)\|^{2}_{L^{2}(\Omega)}+\|\varphi\|^{2}_{L^{2}(0,T;L^{2}(\omega))}\right).
				\end{equation}
			\end{enumerate}
		\end{proposition}
		
		\begin{proof}
			Throughout this proof $C\geq 1$ will be an independent contant of $\varepsilon$ which will be changed from one line to another. $C_{0}$ and $T_{0}$ the constants in Proposition \ref{P2}.\\ 
			(1) The proof is based on the classic cut-off technique. Let us now consider 
			$\vartheta\in \mathcal{C}^{\infty}(\mathbb{R}^{d})$ that check $\vartheta=1$ in a neighborhood of $\omega_{0}$ and $\text{supp}(\vartheta)\subset\omega$. Define
			\begin{eqnarray}
				\varphi_{1}(x,t)=\vartheta(x)\varphi(x,t)\quad\text{and}\quad \varphi_{2}(x,t)=(1-\vartheta(x))\varphi(x,t),\quad\text{on}\;\;\Omega_{T}.
			\end{eqnarray}
			\textcolor{blue}{\textbf{Estimation of $\varphi_{1}$.}} 
			We will estimate $\varphi_{1}$ using  Proposition \ref{P0}. Firstly, one has
			\begin{eqnarray*}
				\partial_{t}\varphi_{1}+\varepsilon\Delta\varphi_{1}+\nabla\cdot(\varphi_{1}\mathbf{B}(x,t))= f(x,t),\quad \text{on}\;\;\Omega_{T},
			\end{eqnarray*}
			where
			\begin{eqnarray*}
				f(x,t):= \varepsilon\displaystyle\sum_{i=1}^{d}\partial_{x_{i}}\left(2\varphi\partial_{x_{i}}\vartheta\right)-(\varepsilon\Delta\vartheta-\mathbf{B}\cdot\nabla\vartheta)\varphi.
			\end{eqnarray*}
			To apply Proposition \ref{P0}, we truncate $\varphi_{1}$ by $\psi\in\mathcal{C}^{\infty}(\mathbb{R})$, $\psi(t_{0}-T_{0})=1$ and $\psi(t_{0})=0$. Taking  $\varphi_{3}(x,t):=\psi(t)\varphi_{1}(x,t)$. Then 
			\begin{eqnarray*}
				\partial_{t}\varphi_{3}+\varepsilon\Delta\varphi_{3}+\nabla\cdot(\varphi_{3}\mathbf{B}(x,t))&=& \psi^{'}(t)\varphi_{1}(x,t)+\psi(t)\left[\partial_{t}\varphi_{1}+\varepsilon\Delta\varphi_{1}+\nabla\cdot(\varphi_{1}\mathbf{B}(x,t))\right]\\
				&=& \psi^{'}(t)\varphi_{1}(x,t)+\psi(t)f(x,t)\\
				&=& \psi^{'}(t)\vartheta(x)\varphi(x,t)+\psi(t)f(x,t):=h(x,t)
			\end{eqnarray*}
			Hence $\varphi_{3}$ is the solution of $\mathcal{S}(\varnothing,t_{0}-T_{0},t_{0},h,0,\varepsilon,\mathbf{B}).$ Let us apply Proposition \ref{P0} with $f_{0}=[\psi^{'}(t)\vartheta(x)-\psi(t)(\varepsilon\Delta\vartheta-\mathbf{B}\cdot\nabla\vartheta)]\varphi$ and  $f_{i}=2\psi(t)\partial_{x_{i}}\vartheta\varphi,\;\; 1\leqslant i\leqslant d$, we obtain 
			\begin{eqnarray}
				\|\varphi_{3}(\cdot,t_{0}-T_{0})\|^{2}_{L^{2}(\Omega)}\leqslant C\sum_{i=0}^{d}\|f_{i}\|^{2}_{L^{2}(t_{0}-T_{0},t_{0};L^{2}(\Omega))}.
			\end{eqnarray}
			Since $\vartheta$ has support in $\omega$ and $\psi(t_{0}-T_{0})=1$, then 
			\begin{eqnarray}
				\|\varphi_{1}(\cdot,t_{0}-T_{0})\|^{2}_{L^{2}(\Omega)}\leqslant C\|\varphi\|^{2}_{L^{2}(t_{0}-T_{0},t_{0};L^{2}(\omega))}. \label{dd10}
			\end{eqnarray}
			\textcolor{blue}{\textbf{Estimation of $\varphi_{2}$.}} Now, we will estimate $\varphi_{2}$ by decomposing it into two solutions using Propositions \ref{P0} and \ref{P2}. Since $\varphi_{2}=\varphi-\varphi_{1}$, then 
			\begin{eqnarray*}
				\partial_{t}\varphi_{2}+\varepsilon\Delta\varphi_{2}+\nabla\cdot(\varphi_{2}\mathbf{B}(x,t))= -f(x,t),\quad \text{on}\;\;\Omega_{T}.
			\end{eqnarray*}
			Therefore, we decompose $\varphi_{2}$ on $\mathcal{U}\times (t_{0}-T_{0},t_{0})$, 
			 as follows 
			$$\begin{cases}
				\varphi_{2}=\varphi_{4}+\varphi_{5},\\
				\varphi_{4}\;\text{is the solution of}\; \mathcal{S}(\omega_{0},t_{0}-T_{0},t_{0},0,\varphi_{2}(\cdot,t_{0}),\varepsilon,\mathbf{B}),\\
				\varphi_{5}\;\text{is the solution of}\; \mathcal{S}(\omega_{0},t_{0}-T_{0},t_{0},-f,0,\varepsilon,\mathbf{B}).
			\end{cases}$$
			From Proposition \ref{P2} and $\varphi_{2}(x,t_{0})=(1-\vartheta(x))\varphi(x,t_{0})$, we obtain 
			\begin{equation} 
				\|\varphi_{4}(\cdot, t_{0}-T_{0})\|^{2}_{L^{2}(\mathcal{U})}\leqslant C\exp\left(\frac{-C_{0}}{\varepsilon}\right)\|\varphi(\cdot, t_{0})\|^{2}_{L^{2}(\mathcal{U})}. \label{d7}
			\end{equation}
			Concerning $\varphi_{5}$, by application of Proposition \ref{P0} with $f_{0}=(\varepsilon\Delta\vartheta-\mathbf{B}\cdot\nabla\vartheta)\varphi$ and  $f_{i}=-2\partial_{x_{i}}\vartheta\varphi,\;\; 1\leqslant i\leqslant d$, we get
			\begin{eqnarray*}
				\|\varphi_{5}(\cdot,t_{0}-T_{0})\|^{2}_{L^{2}(\mathcal{U})}\leqslant C\sum_{i=0}^{d}\|f_{i}\|^{2}_{L^{2}(t_{0}-T_{0},t_{0};L^{2}(\mathcal{U}))}.
			\end{eqnarray*}
			Since $\vartheta$ has support in $\omega$, then 
			\begin{eqnarray}
				\|\varphi_{5}(\;\cdot,t_{0}-T_{0})\|^{2}_{L^{2}(\mathcal{U})}&\leqslant & C\|\varphi\|^{2}_{L^{2}(t_{0}-T_{0},t_{0};L^{2}(\omega\setminus\overline{\omega_{0}}))} \nonumber\\
				&\leqslant & C\|\varphi\|^{2}_{L^{2}(t_{0}-T_{0},t_{0};L^{2}(\omega))}. \label{d8}
			\end{eqnarray}
			The function $\vartheta=1$ in a neighborhood of $\omega_{0}$ implies that $\varphi_{2}$ has a support in $\mathcal{U}$, thus form \eqref{d7} and \eqref{d8}, we obtain
			\begin{eqnarray}
				\|\varphi_{2}(\cdot,t_{0}-T_{0})\|^{2}_{L^{2}(\Omega)}&= & \|\varphi_{2}(\cdot,t_{0}-T_{0})\|^{2}_{L^{2}(\mathcal{U})} \nonumber\\
				&\leqslant & 2\left( \|\varphi_{4}(\cdot,t_{0}-T_{0})\|^{2}_{L^{2}(\mathcal{U})}+\|\varphi_{5}(\cdot,t_{0}-T_{0})\|^{2}_{L^{2}(\mathcal{U})}\right) \nonumber\\
				&\leqslant & C\left(\exp\left(\frac{-C_{0}}{\varepsilon}\right)\|\varphi(\cdot, t_{0})\|^{2}_{L^{2}(\mathcal{U})}+\|\varphi\|^{2}_{L^{2}(t_{0}-T_{0},t_{0};L^{2}(\omega))}\right).\nonumber\\ \label{d9}
			\end{eqnarray}
			Finally, using \eqref{dd10} and \eqref{d9}, we obtain 
			\begin{eqnarray}
				\|\varphi(\cdot,t_{0}-T_{0})\|^{2}_{L^{2}(\Omega)}&\leqslant & C\left(\exp\left(\frac{-C_{0}}{\varepsilon}\right)\|\varphi(\cdot, t_{0})\|^{2}_{L^{2}(\mathcal{U})}+\|\varphi\|^{2}_{L^{2}(t_{0}-T_{0},t_{0};L^{2}(\omega))}\right) \nonumber\\
				&\leqslant & C\left(\exp\left(\frac{-C_{0}}{\varepsilon}\right)\|\varphi(\cdot, t_{0})\|^{2}_{L^{2}(\Omega)}+\|\varphi\|^{2}_{L^{2}(t_{0}-T_{0},t_{0};L^{2}(\omega))}\right).\nonumber\\ \label{d11}
			\end{eqnarray}
			(2)  Let $m$ be an integer such that $1\leqslant m\leqslant \frac{T}{T_{0}}$. From the first dissipation estimate \eqref{disspation 2}, we get 
			\begin{equation*} 
				\|\varphi(\cdot, (k-1)T_{0})\|^{2}_{L^{2}(\mathcal{U})}\leqslant C\left(\exp\left(\frac{-C_{0}}{\varepsilon}\right)\|\varphi(\cdot, kT_{0})\|^{2}_{L^{2}(\Omega)}+\int_{(k-1)T_{0}}^{kT_{0}}\int_{\omega}|\varphi|^{2}\mathrm{~d}x\mathrm{~d}t\right),
			\end{equation*}
			for all $k=1,2,\cdots, m$. This last estimate gives
			\begin{eqnarray} \label{d5}
				&&C^{k-1}\exp\left(\frac{-(k-1)C_{0}}{\varepsilon}\right)\|\varphi(\cdot, (k-1)T_{0})\|^{2}_{L^{2}(\Omega)}-C^{k}\exp\left(\frac{-kC_{0}}{\varepsilon}\right)\|\varphi(\cdot, kT_{0})\|^{2}_{L^{2}(\Omega)} \nonumber\\
				&&\quad\quad\leqslant C^{m}\int_{(k-1)T_{0}}^{kT_{0}}\int_{\omega}|\varphi|^{2}\mathrm{~d}x\mathrm{~d}t.
			\end{eqnarray}
			Summing \eqref{d5} from $1$ to $m,$ we obtain 
			\begin{eqnarray} \label{d10}
				\|\varphi(\cdot, 0)\|^{2}_{L^{2}(\Omega)}&\leqslant & C^{m}\left(\exp\left(\frac{-mC_{0}}{\varepsilon}\right)\|\varphi(\cdot, mT_{0})\|^{2}_{L^{2}(\Omega)}+\int_{0}^{mT_{0}}\int_{\omega}|\varphi|^{2}\mathrm{~d}x\mathrm{~d}t\right)\nonumber\\
				&\leqslant & C^{m}\left(\exp\left(\frac{-mC_{0}}{\varepsilon}\right)\|\varphi(\cdot, mT_{0})\|^{2}_{L^{2}(\Omega)}+\int_{0}^{T}\int_{\omega}|\varphi|^{2}\mathrm{~d}x\mathrm{~d}t\right).\nonumber\\
			\end{eqnarray}
			Since $\mathbf{B}(x,t)\cdot\mathbf{n}(x)\geq 0$ on $\Gamma_{T}$. By Agmon inequality \eqref{A3}, we have 
			\begin{eqnarray}
				\| \varphi(\cdot,mT_{0})\|^{2}_{L^{2}(\Omega)}
				&\leqslant & \exp\left(C_{\B}(t-mT_{0})\right)\| \varphi(\cdot,t)\|^{2}_{L^{2}(\Omega)} \nonumber \\
				&\leqslant &  \exp\left(C_{\B}T_{0}\right)\| \varphi(\cdot,t)\|^{2}_{L^{2}(\Omega)}. \label{d6}
			\end{eqnarray}
			From \eqref{d10} and \eqref{d6}, we obtain the dissipation estimate \eqref{disspation 3}.
		\end{proof}
		\section{Carleman's estimation for the solutions of \eqref{s2} in the case of tangential transport} \label{Section 6}
		In order to establish an observability inequality with an observability constant $\exp\left(\frac{C}{\varepsilon}\right)$, we have to show a Carleman estimate for the solutions of the adjoint system \eqref{s2} while satisfying the constraint $s\geq \frac{C}{\varepsilon}$ and $\lambda\geq C$ (see parameters $s$ and $\lambda$ below), where $C>0$ is a constant independent of $\varepsilon$. However, the presence of a transport term and the constraint pose challenges in this endeavor. To address this, we tackle the problem within the framework of tangential transport. The techniques employed draw inspiration from previous works referenced as \cite{ barcena2021cost, fernandez2006exact, guerrero2007singular} and 
		we have to re-check certain steps from \cite{ettahri:hal-04131920}.
		\par
		Let us consider tangential transport $\B(x,t)$, i.e., $\B(x,t)\cdot\n (x)=0$ on $\Gamma_{T}$. In this case \eqref{s2} becomes 
		\begin{equation} 
			\left\{
			\begin{aligned}
				\partial_t \varphi+\varepsilon\Delta \varphi+\nabla\cdot\left(\varphi \mathbf{B}\right) &=0 & & \text { in } \Omega_{T}, \\
				\partial_{\n}\varphi &=0 & & \text { on } \Gamma_{T}, \\
				\varphi(\cdot,T)&=\varphi_{T} & & \text { in } \Omega.
				\label{2}
			\end{aligned}
			\right.
		\end{equation}
		We introduce the following positive weight functions $\alpha_{\pm}$ and $\xi_{\pm}$ which depend only on $\Omega$ and $\omega$:
		\begin{eqnarray*}
			\alpha_{\pm}(x, t):=\frac{\exp(6\lambda)-\exp(4\lambda\pm\lambda\eta(x))}{t(T-t)} \quad\mbox{and}\quad \xi_{\pm}(x, t):=\frac{\exp(4\lambda\pm\lambda\eta(x))}{t(T-t)},
		\end{eqnarray*}
		where $\lambda \geq 1$ and $\eta=\eta(x)$ is a function in $\mathcal{C}^2(\overline{\Omega})$ satisfying
		\begin{equation}
			\label{eta}
			\eta>0 \text { in } \Omega, \quad \eta=0 \text { on } \Gamma,\quad  \inf_{\Omega \backslash \omega^{\prime}}\left|\nabla \eta(x)\right|=\delta>0,\quad \|\eta\|_{\infty}=1,
		\end{equation}
		where $\omega^{\prime} \subset \subset \omega$ is a nonempty open set. If $\Omega$ is a domain with $\mathcal{C}^{2}$ smoothness, the paper \cite{fursikov1996controllability} provides a proof of the existence of $\eta$ that satisfies \eqref{eta}.
		\par
		The Carleman estimate we will use is as follows.
		\begin{proposition} \label{P4}
			Let $T>0$, $\varepsilon\in (0,1)$, $\Omega$ is a $\mathcal{C}^{2}$ domain, $\omega \subset \subset \Omega$ is a nonempty open set and assume that $\B\in W^{1,\infty}(\Omega_{T})^{d}$ such that $\B(x,t)\cdot\n(x)=0$ on $\Gamma_{T}$. Then there are constants $C>0$ and $\lambda_{1}, s_{1}\geq 1$ depend only on $\omega$ and $\Omega$ such that
			\begin{align}
				&s^{2}\lambda^{2}\int_{\Omega_{T}}\exp(-2s\alpha_{+})\xi_{+}^{3}|\varphi|^{2}\mathrm{~d}x\mathrm{~d}t +\int_{\Omega_{T}}|\nabla\varphi|^{2}\xi_{+}\mathrm{~d}x\mathrm{~d}t \nonumber\\
				&\quad\quad \leqslant Cs^{2}\lambda^{2}\int_{\omega\times (0,T)}\exp(-2s\alpha_{+})\xi_{+}^{3}|\varphi|^{2}\mathrm{~d}x\mathrm{~d}t,
				\label{Carleman}
			\end{align}
			for any $\varphi$ solution of \eqref{s2} with data $\varphi_{T}\in L^{2}(\Omega)$, $\lambda\geq \lambda_{1}$, $s\geq \frac{s_{1}}{\varepsilon}(T^{2}+T)\mathcal{B}_{T}$ and
			\begin{eqnarray}
				\mathcal{B}_{T}&=:& 1+\|\B\|_{L^{\infty}(\Omega_{T})}+\|\nabla\cdot \B\|_{L^{\infty}(\Omega_{T})}+\|\nabla \B\|_{L^{\infty}(\Omega_{T})}+\|\B\|_{L^{\infty}(\Omega_{T})}^{1/2} \nonumber\\
				&& + \|\nabla\cdot\B\|_{L^{\infty}(\Omega_{T})}^{1/2}+\|\nabla\cdot\B\|_{L^{\infty}(\Omega_{T})}^{2/3}+\|\partial_{t}\B\|_{L^{\infty}(\Omega_{T})}^{1/2}.\label{constante s}
			\end{eqnarray}
		\end{proposition}
		\begin{proof}
			 Throughout the proof $C$, $c$, $s_{1}$ and $\lambda_{1}$ will denote generic constants
			which are independent of $\varepsilon$, $s$, $\lambda$ and $B$. These constants may vary even from line to line. $\|\cdot\|_{\infty}$ designates the norm $\|\cdot\|_{L^{\infty}(\Omega_{T})}$.\\
			\par 
			To derive the global estimate, we will give the proof in several steps:
			Initially, a change of variables is implemented to acquire functions that display decay characteristics at both the initial time $t=0$ and the final time $t=T$.
			Subsequently, we assess and approximate the scalar product that arises naturally during the change of variables.
			Afterwards, we draw preliminary conclusions by examining the boundary terms on the left-hand side of the inequality.
			We then estimate the local gradient term.
			Additionally, we simplify the boundary terms and revert the change of variables to obtain the desired estimate.\\\\
			\textcolor{blue}{\textbf{Step 1: An auxiliary problem.}} For all strong solution $\varphi$ of the adjoint system \eqref{2} and for
			all $s\geq 1$, we define 
			\begin{eqnarray}
				\psi_{\pm}:=\exp(-s\alpha_{\pm})\varphi. \label{c-4}
			\end{eqnarray}
			\par
			
			Using the fact that the tangential derivative, noted $\nabla_{\Gamma}$ depends only on values on $\Gamma$ and $\eta=0$ on $\Gamma$, we have the following equations.
			\begin{eqnarray}
				&&\alpha_{+}=\alpha_{-}:=\alpha,\quad \xi_{+}=\xi_{-}:=\xi\;\;\text{on}\;\Gamma_{T}, \label{c4}\\
				&& \psi_{+}=\psi_{-}:=\psi,\quad \nabla_{\Gamma}\psi_{+}=\nabla_{\Gamma}\psi_{-}=\nabla_{\Gamma}\psi\;\;\text{on}\;\Gamma_{T}.\label{c5}
			\end{eqnarray}
			\par
			The problem solved by $\psi_{\pm}$ is given by
			\begin{eqnarray*}
				\partial_{t}\psi_{\pm}+\varepsilon\Delta\psi_{\pm}+\nabla\cdot(\psi_{\pm} \mathbf{B}) &=& -s\partial_{t}\alpha_{\pm}\psi_{\pm}
				-\varepsilon s^{2}\lambda^{2}|\nabla\eta|^{2}
				\xi^{2}_{\pm}\psi_{\pm} \pm 2\varepsilon s\lambda \xi_{\pm}(\nabla\eta\cdot\nabla\psi_{\pm})\\
				&& + \varepsilon s\lambda^{2}|\nabla\eta|^{2}\xi_{\pm}\psi_{\pm}
				\pm  \varepsilon s \lambda\Delta\eta\xi_{\pm}\psi_{\pm} \pm s\lambda (\nabla\eta\cdot\mathbf{B})\xi_{\pm}\psi_{\pm}.
			\end{eqnarray*}
			We rewrite this equation as 
			\begin{eqnarray}
				L_{1}\psi_{\pm}+L_{2}\psi_{\pm}=L_{3}\psi_{\pm}, \label{c1}
			\end{eqnarray}
			where 
			\begin{eqnarray}
				L_{1}\psi_{\pm}&=& - 2\varepsilon s\lambda^{2}|\nabla\eta|^{2}\xi_{\pm}\psi_{\pm}\mp 2\varepsilon s\lambda(\nabla\eta\cdot\nabla\psi_{\pm})\xi_{\pm}+\partial_{t}\psi_{\pm}+\mathbf{B}\cdot\nabla\psi_{\pm}, \label{c-2}\\
				L_{2}\psi_{\pm} &=& \varepsilon s^{2}\lambda^{2}|\nabla\eta|^{2}\xi_{\pm}^{2}\psi_{\pm}+\varepsilon\Delta\psi_{\pm}+s\partial_{t}\alpha_{\pm}\psi_{\pm} \mp s\lambda (\nabla\eta\cdot\mathbf{B})\xi_{\pm}\psi_{\pm}, \label{c-1}\\
				L_{3}\psi_{\pm}&=& \pm \varepsilon s\lambda\Delta\eta\xi_{\pm}\psi_{\pm}
				-\varepsilon s\lambda^{2}|\nabla\eta|^{2}\xi_{\pm}\psi_{\pm}-(\nabla\cdot\mathbf{B})\psi_{\pm}. \label{c0}
			\end{eqnarray}
			Applying $\|\cdot\|_{L^{2}(\Omega_{T})}$ to equation \eqref{c1}, we obtain 
			\begin{eqnarray}
				\|L_{1}\psi_{\pm}\|^{2}_{L^{2}(\Omega_{T})}+2	(L_{1}\psi_{\pm},L_{2}\psi_{\pm})_{L^{2}(\Omega_{T})}+\|L_{2}\psi_{\pm}\|^{2}_{L^{2}(\Omega_{T})}=\|L_{3}\psi_{\pm}\|^{2}_{L^{2}(\Omega_{T})}. \label{c2}
			\end{eqnarray}
			The main novelty will be on how to deal with the term $B\cdot\nabla\psi_\pm$.\\\\
			\textcolor{blue}{\textbf{Step 2. Estimating the mixed terms in \eqref{c2} from below.}}
			The main idea is to expand the term $(L_{1}\psi_{\pm},L_{2}\psi_{\pm})_{L^{2}(\Omega_{T})}$ and use the particular structure of $\alpha_{\pm}$ and the fact
			that $s$ is large enough in order to obtain large positive terms in this scalar product. Denoting by $\left(L_{i}\psi_{\pm}\right)_{j}$ the $j$-th term in the above expression of $L_{i}\psi_{\pm}$ and for all $1\leqslant i,j \leqslant 4$, we note
			\begin{eqnarray}
				I_{\pm}^{ij}&:=&((L_{1}\psi_{\pm})_{i},(L_{2}\psi_{\pm})_{j})_{L^{2}(\Omega_{T})}\;\mbox{and}\;
				I_{\pm}^{j}:=(L_{1}\psi_{\pm},(L_{2}\psi_{\pm})_{j})_{L^{2}(\Omega_{T})}.
			\end{eqnarray}
			We obtain
			\begin{eqnarray}
				I_{\pm}^{j}=\sum_{i=1}^{4}I_{\pm}^{ij}\quad \mbox{and}\quad 
				(L_{1}\psi_{\pm},L_{2}\psi_{\pm})_{L^{2}(\Omega_{T})}=\sum_{j=1}^{4}I_{\pm}^{j}. \label{c3}
			\end{eqnarray}
			Let us recall some of the vector calculus formulas we will be using in this step.
			\begin{eqnarray}
				&&\nabla(\nabla\eta\cdot\nabla\psi_{\pm})\cdot\nabla\psi_{\pm}=\nabla^{2}\eta (\nabla\psi_{\pm},\nabla\psi_{\pm})+\frac{1}{2}\nabla\eta\cdot\nabla(|\nabla\psi_{\pm}|^{2}), \label{cv1}\\
				&&\nabla(\mathbf{B}\cdot\nabla\psi_{\pm})\cdot\nabla\psi_{\pm}=\frac{1}{2}\nabla(|\nabla\psi_{\pm}|^{2})\cdot\mathbf{B}+(\nabla \mathbf{B}\cdot\nabla\psi_{\pm})\nabla\psi_{\pm}, \label{cv2}\\
				&&\nabla\cdot((\nabla\eta\cdot\B)\nabla\eta)=\frac{1}{2}\nabla(|\nabla\eta|^{2})\cdot\B+(\nabla \B\cdot\nabla\eta)\nabla\eta +(\nabla\eta\cdot\B)\Delta\eta, \label{cv3}\\
				&&\nabla\cdot((\nabla\eta\cdot\B)\B)=\nabla^{2}\eta(\B,\B)+(\nabla \B\cdot\B)\cdot\nabla\eta+(\nabla\eta\cdot\B)\nabla\cdot\B \label{cv4},
			\end{eqnarray}
		where $\nabla^{2}\eta$ denotes the Hessian matrix of $\eta$ (it is considered as a symmetrical bilinear form) and $\nabla \mathbf{B}$ is the Jacobian matrix of $B$.\\
			Let us compute each term $I_{\pm}^{j}$.\\
			\textcolor{blue}{\textbf{Step $2a$. Estimate from below of $I_{\pm}^{1}$.}}\\
			Firstly, one has
			\begin{eqnarray}
				I_{\pm}^{11}=- 2\varepsilon^{2}s^{3}\lambda^{4}\int_{\Omega_{T}}|\nabla\eta|^{4}\xi_{\pm}^{3}|\psi_{\pm}|^{2}\mathrm{~d}x\mathrm{~d}t.
			\end{eqnarray}
			A simple calculation and an integrating by parts yield 
			\begin{eqnarray*}
				I_{\pm}^{21}&=&\mp 2\varepsilon^{2}s^{3}\lambda^{3}\int_{\Omega_{T}}(\nabla\eta\cdot\nabla\psi_{\pm})|\nabla\eta|^{2}\xi^{3}_{\pm}\psi_{\pm}\mathrm{~d}x\mathrm{~d}t\\
				&=& 3\varepsilon^{2}s^{3}\lambda^{4}\int_{\Omega_{T}} |\nabla\eta|^{4}\xi_{\pm}^{3}|\psi_{\pm}|^{2}\mathrm{~d}x\mathrm{~d}t \pm \varepsilon^{2}s^{3}\lambda^{3}\int_{\Omega_{T}}\nabla\cdot( |\nabla\eta|^{2}\nabla\eta)\xi_{\pm}^{3}|\psi_{\pm}|^{2}\mathrm{~d}x\mathrm{~d}t\\
				&& \mp\varepsilon^{2}s^{3}\lambda^{3}\int_{\Gamma_{T}} |\nabla\eta|^{2}\partial_{\n}\eta\xi^{3}|\psi|^{2}\mathrm{~d}\sigma\mathrm{~d}t\\
				&:=& (I_{\pm}^{21})_{1}+(I_{\pm}^{21})_{2}+(I_{\pm}^{21})_{3}.
			\end{eqnarray*}
			Clearly  $I_{\pm}^{11}+(I_{\pm}^{21})_{1}=\varepsilon^{2}s^{3}\lambda^{4}\displaystyle\int_{\Omega_{T}} |\nabla\eta|^{4}\xi_{\pm}^{3}|\psi_{\pm}|^{2}\mathrm{~d}x\mathrm{~d}t\geq 0$. From Properties \eqref{eta} of $\eta$, we obtain
			\begin{eqnarray}
				I_{\pm}^{11}+(I_{\pm}^{21})_{1} 
				&\geq&  c\, \varepsilon^{2}s^{3}\lambda^{4}\int_{\Omega_{T}}\xi_{\pm}^{3}|\psi_{\pm}|^{2}\mathrm{~d}x\mathrm{~d}t-C\varepsilon^{2}s^{3}\lambda^{4}\int_{\omega^{'}\times (0,T)}\xi_{\pm}^{3}|\psi_{\pm}|^{2}\mathrm{~d}x\mathrm{~d}t \nonumber\\
				&:=& A_{\pm}-B_{\pm},
			\end{eqnarray}  
			for all $c\in (0,\delta^{4})$.
			The term $(I_{\pm}^{21})_{2}$ can be absorbed by $\tilde{A}_{\pm}$ if $\lambda\geq \lambda_{1}$ for large $\lambda_{1}$. Then 
			\begin{eqnarray}
				I_{\pm}^{11}+I_{\pm}^{21}&\geq & A_{\pm}-B_{\pm}+(I_{\pm}^{21})_{3}. \label{c7}
			\end{eqnarray}
			By integration by parts in time and $\psi_{\pm}(\cdot,0)=0=\psi_{\pm}(\cdot,T)$, we have 
			\begin{eqnarray*}
				I_{\pm}^{31}&=& \varepsilon s^{2}\lambda^{2}\int_{\Omega_{T}}|\nabla\eta|^{2}\xi_{\pm}^{2}\psi_{\pm}\partial_{t}\psi_{\pm}\mathrm{~d}x\mathrm{~d}t\\
				&=& -\varepsilon s^{2}\lambda^{2}\int_{\Omega_{T}}|\nabla\eta|^{2}\xi_{\pm}\partial_{t}\xi_{\pm}|\psi_{\pm}|^{2}\mathrm{~d}x\mathrm{~d}t.
			\end{eqnarray*}
			Since $|\partial_{t}\xi_{\pm}|\leqslant T\xi_{\pm}^{2}$, then 
			\begin{eqnarray*}
				|I_{\pm}^{31}| &\leqslant & C\varepsilon s^{2}\lambda^{2}T\int_{\Omega_{T}}\xi_{\pm}^{3}|\psi_{\pm}|^{2}\mathrm{~d}x\mathrm{~d}t.
			\end{eqnarray*}
			This last integral can be absorbed by $A_{\pm}$ if we take $\lambda\geq 1$ and $s\geq \frac{s_{1}}{\varepsilon}T$.\\
			Using $\B(x,t)\cdot\n(x)=0$ on $\Gamma_{T}$ and integration by parts, one has 
			\begin{eqnarray*}
				I_{\pm}^{41}
				&=& \frac{\varepsilon s^{2}\lambda^{2}}{2}\int_{\Omega_{T}}|\nabla\eta|^{2}\xi_{\pm}^{2}(\mathbf{B}\cdot\nabla|\psi_{\pm}|^{2})\mathrm{~d}x\mathrm{~d}t\\
			    &=& -\frac{\varepsilon s^{2}\lambda^{2}}{2}\int_{\Omega_{T}}(\nabla(|\nabla\eta|^{2})\cdot\mathbf{B})\xi_{\pm}^{2}|\psi_{\pm}|^{2}\mathrm{~d}x\mathrm{~d}t \mp \varepsilon s^{2}\lambda^{3}\int_{\Omega_{T}}|\nabla\eta|^{2}\xi_{\pm}^{2}(\nabla\eta\cdot\mathbf{B})|\psi_{\pm}|^{2}\mathrm{~d}x\mathrm{~d}t\\
				&& -\frac{\varepsilon s^{2}\lambda^{2}}{2}\int_{\Omega_{T}}|\nabla\eta|^{2}\xi_{\pm}^{2}\nabla\cdot\mathbf{B}|\psi_{\pm}|^{2}\mathrm{~d}x\mathrm{~d}t\\
				&\leqslant & C\varepsilon s^{2}\lambda^{2}(\|\mathbf{B}\|_{\infty}+\lambda\|\mathbf{B}\|_{\infty}+\|\nabla\cdot\mathbf{B}\|_{\infty})\int_{\Omega_{T}}\xi_{\pm}^{2}|\psi_{\pm}|^{2}\mathrm{~d}x\mathrm{~d}t.
			\end{eqnarray*}
			Using $\xi_{\pm}\geq\frac{4}{T^{2}}$, we get 
			\begin{eqnarray*}
				I_{\pm}^{41}
				&\leqslant & C\varepsilon s^{2}\lambda^{2}T^{2}(\|\mathbf{B}\|_{\infty}+\lambda\|\mathbf{B}\|_{\infty}+\|\nabla\cdot\mathbf{B}\|_{\infty})\int_{\Omega_{T}}\xi_{\pm}^{3}|\psi_{\pm}|^{2}\mathrm{~d}x\mathrm{~d}t,
			\end{eqnarray*}
			which is absorbed by $A_{\pm}$ by taking $\lambda\geq 1$ and $s\geq \frac{s_{1}}{\varepsilon}T^{2}(\|\mathbf{B}\|_{\infty}+\|\nabla\cdot\mathbf{B}\|_{\infty})$.
			Consequently, we obtain 
			\begin{eqnarray}
				I_{\pm}^{1}=\sum_{i=1}^{4}I_{\pm}^{i1}&\geq &  c\,\varepsilon^{2}s^{3}\lambda^{4}\int_{\Omega_{T}}\xi_{\pm}^{3}|\psi_{\pm}|^{2}\mathrm{~d}x\mathrm{~d}t-C\varepsilon^{2}s^{3}\lambda^{4}\int_{\omega^{'}\times (0,T)}\xi_{\pm}^{3}|\psi_{\pm}|^{2}\mathrm{~d}x\mathrm{~d}t \nonumber\\
				&& \quad\quad\mp\varepsilon^{2}s^{3}\lambda^{3}\int_{\Gamma_{T}} |\nabla\eta|^{2}\partial_{\n}\eta\xi^{3}|\psi|^{2}\mathrm{~d}\sigma\mathrm{~d}t, \label{c8}
			\end{eqnarray}
			for any $\lambda\geq \lambda_{1}$ and any $s\geq \frac{s_{1}}{\varepsilon}\left(T+T^{2}\left(\|\mathbf{B}\|_{\infty}+\|\nabla\cdot\mathbf{B}\|_{\infty}\right)\right)$.\\\\
			\textcolor{blue}{\textbf{Step $2b$. Estimate from below of $I_{\pm}^{2}$.}}
			By integration by parts, we get
			\begin{eqnarray*}
				I_{\pm}^{12}&=& - 2\varepsilon^{2}s\lambda^{2}\int_{\Omega_{T}}|\nabla\eta|^{2}\xi_{\pm}\psi_{\pm}\Delta\psi_{\pm}\mathrm{~d}x\mathrm{~d}t\\
				&=&  - 2\varepsilon^{2}s\lambda^{2}\int_{\Gamma_{T}}|\nabla\eta|^{2}\xi\psi\partial_{\n}\psi_{\pm}\mathrm{~d}\sigma\mathrm{~d}t + 2\varepsilon^{2}s\lambda^{2}\int_{\Omega_{T}}(\nabla(|\nabla\eta|^{2})\cdot\nabla\psi_{\pm})  \xi_{\pm}\psi_{\pm}\mathrm{~d}x\mathrm{~d}t\\
				& & \pm  2\varepsilon^{2}s\lambda^{3}\int_{\Omega_{T}}|\nabla\eta|^{2}(\nabla\eta\cdot\nabla\psi_{\pm})\xi_{\pm}\psi_{\pm}\mathrm{~d}x\mathrm{~d}t + 2\varepsilon^{2}s\lambda^{2}\int_{\Omega_{T}}|\nabla\eta|^{2}|\nabla\psi_{\pm}|^{2}\xi_{\pm}\mathrm{~d}x\mathrm{~d}t\\
				&=:& (	I_{\pm}^{12})_{1}+(	I_{\pm}^{12})_{2}+(	I_{\pm}^{12})_{3}+(	I_{\pm}^{12})_{4}.
			\end{eqnarray*}
			We will leave $(	I_{\pm}^{12})_{1}$ and $(	I_{\pm}^{12})_{4}$ on the left-hand side. Concerning $(I_{\pm}^{12})_{2}$ and $(	I_{\pm}^{12})_{3}$, we can apply Holder's inequality
			\begin{eqnarray*}
				|(I_{\pm}^{12})_{2} | 
				&\leqslant & C\varepsilon^{2}s\lambda^{4}\int_{\Omega_{T}} \xi_{\pm}|\psi_{\pm}|^{2}\mathrm{~d}x\mathrm{~d}t+C\varepsilon^{2}s\int_{\Omega_{T}}|\nabla\psi_{\pm}|^{2} \xi_{\pm}\mathrm{~d}x\mathrm{~d}t.
			\end{eqnarray*}
			Since $\xi_{\pm}\geq \frac{4}{T^{2}}$, then for all $s\geq s_{1}T^{2}$ 
			\begin{eqnarray*}
				|(I_{\pm}^{12})_{2} | &\leqslant& C\varepsilon^{2}s^{2}\lambda^{4}\int_{\Omega_{T}} \xi_{\pm}^{2}|\psi_{\pm}|^{2}\mathrm{~d}x\mathrm{~d}t+C\varepsilon^{2}s\int_{\Omega_{T}}|\nabla\psi_{\pm}|^{2} \xi_{\pm}\mathrm{~d}x\mathrm{~d}t.
			\end{eqnarray*}
			And 
				
			\begin{eqnarray*}
				|(I_{\pm}^{12})_{3} | 
				&\leqslant & C\varepsilon^{2}s^{2}\lambda^{4}\int_{\Omega_{T}} \xi_{\pm}^{2}|\psi_{\pm}|^{2}\mathrm{~d}x\mathrm{~d}t+C\varepsilon^{2}\lambda^{2}\int_{\Omega_{T}}|\nabla\psi_{\pm}|^{2}\mathrm{~d}x\mathrm{~d}t.
			\end{eqnarray*}
			We conclude that 
			\begin{eqnarray}
				I_{\pm}^{12}&\geq& - 2\varepsilon^{2}s\lambda^{2}\int_{\Gamma_{T}}|\nabla\eta|^{2}\xi\psi\partial_{\n}\psi_{\pm}\mathrm{~d}\sigma\mathrm{~d}t-C\varepsilon^{2}s^{2}\lambda^{4}\int_{\Omega_{T}} \xi_{\pm}^{2}|\psi_{\pm}|^{2}\mathrm{~d}x\mathrm{~d}t \nonumber\\
				&& - C\varepsilon^{2}\int_{\Omega_{T}}(s\xi_{\pm}+\lambda^{2})|\nabla\psi_{\pm}|^{2}\mathrm{~d}x\mathrm{~d}t + 2\varepsilon^{2}s\lambda^{2}\int_{\Omega_{T}}|\nabla\eta|^{2}|\nabla\psi_{\pm}|^{2}\xi_{\pm}\mathrm{~d}x\mathrm{~d}t,\nonumber \\
				\label{c9}
			\end{eqnarray}
			for all $s\geq s_{1}T^{2}$. On the other hand, by integrations by parts, we obtain
			\begin{eqnarray*}
				I_{\pm}^{22}&=& \mp 2\varepsilon^{2}s\lambda\int_{\Omega_{T}}(\nabla\eta\cdot\nabla\psi_{\pm})\xi_{\pm}\Delta\psi_{\pm}\mathrm{~d}x\mathrm{~d}t\\
				&=& \mp 2\varepsilon^{2}s\lambda\int_{\Gamma_{T}}(\nabla\eta\cdot\nabla\psi_{\pm})\xi\partial_{\n}\psi_{\pm}\mathrm{~d}\sigma\mathrm{~d}t \pm 2\varepsilon^{2}s\lambda\int_{\Omega_{T}}(\nabla(\nabla\eta\cdot\nabla\psi_{\pm})\cdot\nabla\psi_{\pm} ) \xi_{\pm}\mathrm{~d}x\mathrm{~d}t\\
				&&\pm   2\varepsilon^{2}s\lambda\int_{\Omega_{T}}(\nabla\eta\cdot\nabla\psi_{\pm})(\nabla\xi_{\pm}\cdot\nabla\psi_{\pm})\mathrm{~d}x\mathrm{~d}t.
			\end{eqnarray*}
			By utilizing formula \eqref{cv1}, as well as employing another integration by parts, one has 
			\begin{eqnarray*}
				I_{\pm}^{22} 
				&=&  \mp 2\varepsilon^{2}s\lambda\int_{\Gamma_{T}}(\nabla\eta\cdot\nabla\psi_{\pm})\xi\partial_{\n}\psi_{\pm}\mathrm{~d}\sigma\mathrm{~d}t \pm 2\varepsilon^{2}s\lambda\int_{\Omega_{T}}\nabla^{2}\eta (\nabla\psi_{\pm},\nabla\psi_{\pm})\xi_{\pm}\mathrm{~d}x\mathrm{~d}t\\
				& &  \pm \varepsilon^{2}s\lambda\int_{\Gamma_{T}}\xi\partial_{\n}\eta|\nabla\psi_{\pm}|^{2}\mathrm{~d}\sigma\mathrm{~d}t - \varepsilon^{2}s\lambda^{2}\int_{\Omega_{T}}|\nabla\eta|^{2}|\nabla\psi_{\pm}|^{2}\xi_{\pm}\mathrm{~d}x\mathrm{~d}t\\
				& & \mp\varepsilon^{2}s\lambda\int_{\Omega_{T}}\Delta\eta\xi_{\pm}|\nabla\psi_{\pm}|^{2}\mathrm{~d}x\mathrm{~d}t+2\varepsilon^{2}s\lambda^{2}\int_{\Omega_{T}}|\nabla\eta\cdot\nabla\psi_{\pm}|^{2}\xi_{\pm}\mathrm{~d}x\mathrm{~d}t\\
				&:=& (I_{\pm}^{22})_{1}+(I_{\pm}^{22})_{2}+(I_{\pm}^{22})_{3}+(I_{\pm}^{22})_{4}+(I_{\pm}^{22})_{5}+(I_{\pm}^{22})_{6}.
			\end{eqnarray*}
			Using $|\nabla\psi_{\pm}|^{2}=|\nabla_{\Gamma}\psi|^{2}+(\partial_{\nu}\psi_{\pm})^{2}$ and $\eta=0$ on $\Gamma$, we obtain 
			$$(I_{\pm}^{22})_{1}+(I_{\pm}^{22})_{3}=\mp \varepsilon^{2}s\lambda\displaystyle\int_{\Gamma_{T}}\partial_{\n}\eta\xi|\partial_{\n}\psi_{\pm}|^{2}\mathrm{~d}\sigma\mathrm{~d}t \pm \varepsilon^{2}s\lambda\displaystyle\int_{\Gamma_{T}}\partial_{\n}\eta\xi|\nabla_{\Gamma}\psi|^{2}\mathrm{~d}\sigma\mathrm{~d}t.$$
			Clearly, we have , $(I_{\pm}^{22})_{6}\geq 0$ and
			$|(I_{\pm}^{22})_{2}+(I_{\pm}^{22})_{5}|\leqslant C\varepsilon^{2}s\lambda\displaystyle\int_{\Omega_{T}}\xi_{\pm}|\nabla\psi_{\pm}|^{2}\mathrm{~d}x\mathrm{~d}t$. Thus
			\begin{eqnarray}
				I_{\pm}^{22} &\geq& \mp \varepsilon^{2}s\lambda\displaystyle\int_{\Gamma_{T}}\partial_{\n}\eta\xi|\partial_{\n}\psi_{\pm}|^{2}\mathrm{~d}\sigma\mathrm{~d}t \pm \varepsilon^{2}s\lambda\displaystyle\int_{\Gamma_{T}}\partial_{\n}\eta\xi|\nabla_{\Gamma}\psi|^{2}\mathrm{~d}\sigma\mathrm{~d}t \label{c10}\\
				&&-\varepsilon^{2}s\lambda^{2}\int_{\Omega_{T}}|\nabla\eta|^{2}|\nabla\psi_{\pm}|^{2}\xi_{\pm}\mathrm{~d}x\mathrm{~d}t-C\varepsilon^{2}s\lambda\displaystyle\int_{\Omega_{T}}\xi_{\pm}|\nabla\psi_{\pm}|^{2}\mathrm{~d}x\mathrm{~d}t. \nonumber
			\end{eqnarray}
			An integration by parts in space and time and $\nabla\psi_{\pm}$  vanishes at $t=0$ and at $t=T$ yield
			\begin{eqnarray} \label{c11}
				I_{\pm}^{32}=\varepsilon\int_{\Omega_{T}}\partial_{t}\psi_{\pm}\Delta\psi_{\pm}\mathrm{~d}\sigma\mathrm{~d}t= \varepsilon\int_{\Gamma_{T}}\partial_{t}\psi\partial_{\n}\psi_{\pm}\mathrm{~d}\sigma\mathrm{~d}t.
			\end{eqnarray} 
			Sine $\B$ is tangential transport, we also have
			\begin{eqnarray*} 
				I_{\pm}^{42}&=&\varepsilon\int_{\Omega_{T}}(\mathbf{B}\cdot\nabla\psi_{\pm})\Delta\psi_{\pm}\mathrm{~d}x\mathrm{~d}t\\
				&=& \varepsilon\int_{\Gamma_{T}}(\mathbf{B}\cdot\nabla\psi_{\pm})\partial_{\n}\psi_{\pm}\mathrm{~d}\sigma\mathrm{~d}t-\varepsilon\int_{\Omega_{T}}\nabla(\mathbf{B}\cdot\nabla\psi_{\pm})\cdot\nabla\psi_{\pm}\mathrm{~d}x\mathrm{~d}t.
			\end{eqnarray*}
			From the fact that
			 $\mathbf{B}(x,t)\cdot\n(x)=0$ on $\Gamma_{T}$ and formula \eqref{cv2}, we obtain
			 
			\begin{eqnarray} 
				I_{\pm}^{42} 
				&=& \varepsilon\int_{\Gamma_{T}}(\mathbf{B}\cdot\nabla_{\Gamma}\psi)\partial_{\mathbf{n}}\psi_{\pm}\mathrm{~d}\sigma\mathrm{~d}t-\varepsilon\int_{\Omega_{T}}\nabla(\mathbf{B}\cdot\nabla\psi_{\pm})\cdot\nabla\psi_{\pm}\mathrm{~d}x\mathrm{~d}t \nonumber\\
				&=& \varepsilon\int_{\Gamma_{T}}(\mathbf{B}\cdot\nabla_{\Gamma}\psi)\partial_{\mathbf{n}}\psi_{\pm}\mathrm{~d}\sigma\mathrm{~d}t+\frac{\varepsilon}{2}\int_{\Omega_{T}}\nabla\cdot\mathbf{B}|\nabla\psi_{\pm}|^{2}\mathrm{~d}x\mathrm{~d}t \nonumber\\
				&&-\varepsilon\int_{\Omega_{T}}(\nabla \mathbf{B}\cdot\nabla\psi_{\pm})\cdot\nabla\psi_{\pm}\mathrm{~d}x\mathrm{~d}t \nonumber\\
				&&\geq  \varepsilon\int_{\Gamma_{T}}(\mathbf{B}\cdot\nabla_{\Gamma}\psi)\partial_{\mathbf{n}}\psi_{\pm}\mathrm{~d}\sigma\mathrm{~d}t-C\varepsilon\|\nabla \mathbf{B}\|_{\infty}\int_{\Omega_{T}}|\nabla\psi_{\pm}|^{2}\mathrm{~d}x\mathrm{~d}t. \label{cc12}
			\end{eqnarray}
			Considering \eqref{c9}, \eqref{c10}, \eqref{c11} and \eqref{cc12}, we get
			\begin{eqnarray*}
				I_{\pm}^{2}&=&\sum_{i=1}^{4}I_{\pm}^{i2}
				\geq -2\varepsilon^{2}s\lambda^{2}\int_{\Gamma_{T}}|\nabla\eta|^{2}\xi\psi\partial_{\mathbf{n}}\psi_{\pm}\mathrm{~d}\sigma\mathrm{~d}t-C\varepsilon^{2}s^{2}\lambda^{4}\int_{\Omega_{T}} \xi_{\pm}^{2}|\psi_{\pm}|^{2}\mathrm{~d}x\mathrm{~d}t \nonumber\\
				&&- C\varepsilon^{2}\int_{\Omega_{T}}(s\xi_{\pm}+\lambda^{2})|\nabla\psi_{\pm}|^{2}\mathrm{~d}x\mathrm{~d}t+2\varepsilon^{2}s\lambda^{2}\int_{\Omega_{T}}|\nabla\eta|^{2}|\nabla\psi_{\pm}|^{2}\xi_{\pm}\mathrm{~d}x\mathrm{~d}t\\
				&&\mp \varepsilon^{2}s\lambda\displaystyle\int_{\Gamma_{T}}\partial_{\mathbf{n}}\eta\xi|\partial_{\mathbf{n}}\psi_{\pm}|^{2}\mathrm{~d}\sigma\mathrm{~d}t \pm \varepsilon^{2}s\lambda\displaystyle\int_{\Gamma_{T}}\partial_{\n}\eta\xi|\nabla_{\Gamma}\psi|^{2}\mathrm{~d}\sigma\mathrm{~d}t\nonumber\\
				&&-\varepsilon^{2}s\lambda^{2}\int_{\Omega_{T}}|\nabla\eta|^{2}|\nabla\psi_{\pm}|^{2}\xi_{\pm}\mathrm{~d}x\mathrm{~d}t-C\varepsilon^{2}s\lambda\displaystyle\int_{\Omega_{T}}\xi_{\pm}|\nabla\psi_{\pm}|^{2}\mathrm{~d}x\mathrm{~d}t\\
				&&+\varepsilon\int_{\Gamma_{T}}\partial_{t}\psi\partial_{\n}\psi_{\pm}\mathrm{~d}\sigma\mathrm{~d}t+\varepsilon\int_{\Gamma_{T}}(\mathbf{B}\cdot\nabla_{\Gamma}\psi)\partial_{\mathbf{n}}\psi_{\pm}\mathrm{~d}\sigma\mathrm{~d}t \nonumber\\
				&&-C\varepsilon\|\nabla \B\|_{\infty}\int_{\Omega_{T}}|\nabla\psi_{\pm}|^{2}\mathrm{~d}x\mathrm{~d}t,
			\end{eqnarray*}
			for all $s\geq s_{1}T^{2}$. Then for any $\lambda\geq 1$ and $s\geq s_{1}T^{2}$, we have \\
			\begin{eqnarray}
				I_{\pm}^{2} &\geq& -2\varepsilon^{2}s\lambda^{2}\int_{\Gamma_{T}}|\nabla\eta|^{2}\xi\psi\partial_{\n}\psi_{\pm}\mathrm{~d}\sigma\mathrm{~d}t -C\varepsilon^{2}s^{2}\lambda^{4}\int_{\Omega_{T}} \xi_{\pm}^{2}|\psi_{\pm}|^{2}\mathrm{~d}x\mathrm{~d}t \nonumber\\
				&&
				-C\varepsilon^{2}\int_{\Omega_{T}}(s\lambda\xi_{\pm}+\lambda^{2})|\nabla\psi_{\pm}|^{2}\mathrm{~d}x\mathrm{~d}t +\varepsilon^{2}s\lambda^{2}\int_{\Omega_{T}}|\nabla\eta|^{2}|\nabla\psi_{\pm}|^{2}\xi_{\pm}\mathrm{~d}x\mathrm{~d}t \nonumber\\
				&& \mp \varepsilon^{2}s\lambda\displaystyle\int_{\Gamma_{T}}\partial_{\n}\eta\xi|\partial_{\n}\psi_{\pm}|^{2}\mathrm{~d}\sigma\mathrm{~d}t \pm \varepsilon^{2}s\lambda\displaystyle\int_{\Gamma_{T}}\partial_{\n}\eta\xi|\nabla_{\Gamma}\psi|^{2}\mathrm{~d}\sigma\mathrm{~d}t \nonumber\\
				&&
				+\varepsilon\int_{\Gamma_{T}}\partial_{t}\psi\partial_{\n}\psi_{\pm}\mathrm{~d}\sigma\mathrm{~d}t + \varepsilon\int_{\Gamma_{T}}(\B\cdot\nabla_{\Gamma}\psi)\partial_{\n}\psi_{\pm}\mathrm{~d}\sigma\mathrm{~d}t \nonumber\\
				&&
				-C\varepsilon\|\nabla \B\|_{\infty}\int_{\Omega_{T}}|\nabla\psi_{\pm}|^{2}\mathrm{~d}x\mathrm{~d}t := \sum_{i=1}^{9}(I_{\pm}^{2})_{i}. \label{c36}
			\end{eqnarray}
			Thanks to the Properties \eqref{eta} of $\eta$, $(I_{\pm}^{2})_{4}$ can be reduced as follows
			\begin{eqnarray*}
				(I_{\pm}^{2})_{4}
				&\geq & c\,\varepsilon^{2}s\lambda^{2}\int_{\Omega_{T}}|\nabla\psi_{\pm}|^{2}\xi_{\pm}\mathrm{~d}x\mathrm{~d}t-C\varepsilon^{2}s\lambda^{2}\int_{\omega^{'}\times (0,T)}|\nabla\psi_{\pm}|^{2}\xi_{\pm}\mathrm{~d}x\mathrm{~d}t\\
				&:=& C_{\pm}-D_{\pm},
			\end{eqnarray*}
			for all $c\in (0,\delta^{2})$.\\
			Since  $\xi_{\pm}\geq \frac{4}{T^{2}}$, $(I_{\pm}^{2})_{3}$ is absorbed by $C_{\pm}$ if $\lambda\geq\lambda_{1}$ and $s\geq s_{1}T^{2}$ for $\lambda_{1}$ and $s_{1}$ are large enough. The same for $(I_{\pm}^{2})_{9}$ by taking $\lambda\geq 1$ and $s\geq  \frac{s_{1}}{\varepsilon}T^{2}\|\nabla \B\|_{\infty}$.
			Therefore, for $\lambda\geq\lambda_{1}$ and $s\geq \frac{s_{1}}{\varepsilon}T^{2}(1+\|\nabla \B\|_{\infty})$, we obtain 
			\begin{eqnarray}
				I_{\pm}^{2} &\geq& -2\varepsilon^{2}s\lambda^{2}\int_{\Gamma_{T}}|\nabla\eta|^{2}\xi\psi\partial_{\n}\psi_{\pm}\mathrm{~d}\sigma\mathrm{~d}t-C\varepsilon^{2}s^{2}\lambda^{4}\int_{\Omega_{T}} \xi_{\pm}^{2}|\psi_{\pm}|^{2}\mathrm{~d}x\mathrm{~d}t \nonumber\\
				&&+ c\,\varepsilon^{2}s\lambda^{2}\int_{\Omega_{T}}|\nabla\psi_{\pm}|^{2}\xi_{\pm}\mathrm{~d}x\mathrm{~d}t-C\varepsilon^{2}s\lambda^{2}\int_{\omega^{'}\times (0,T)}|\nabla\psi_{\pm}|^{2}\xi_{\pm}\mathrm{~d}x\mathrm{~d}t \nonumber\\
				&& \mp\varepsilon^{2}s\lambda\displaystyle\int_{\Gamma_{T}}\partial_{\n}\eta\xi|\partial_{\n}\psi_{\pm}|^{2}\mathrm{~d}\sigma\mathrm{~d}t \pm \varepsilon^{2}s\lambda\displaystyle\int_{\Gamma_{T}}\partial_{\n}\eta\xi|\nabla_{\Gamma}\psi|^{2}\mathrm{~d}\sigma\mathrm{~d}t \nonumber\\
				&&+\varepsilon\int_{\Gamma_{T}}\partial_{t}\psi\partial_{\n}\psi_{\pm}\mathrm{~d}\sigma\mathrm{~d}t+ \varepsilon\int_{\Gamma_{T}}(\B\cdot\nabla_{\Gamma}\psi)\partial_{\n}\psi_{\pm}\mathrm{~d}\sigma\mathrm{~d}t. \label{c12}
			\end{eqnarray}
			\textcolor{blue}{\textbf{Step 2c. Estimate from below of $I_{\pm}^{3}$.}}
			From $|\partial_{t}\alpha_{\pm}|\leqslant T\xi_{\pm}^{2}$, one has
			\begin{eqnarray}
				\left| I_{\pm}^{13} \right|&\leqslant& 2\varepsilon s^{2}\lambda^{2}\int_{\Omega_{T}} |\nabla\eta|^{2}|\partial_{t}\alpha_{\pm}|\xi_{\pm}|\psi_{\pm}|^{2}\mathrm{~d}x\mathrm{~d}t \nonumber\\
				&\leqslant & C\varepsilon s^{2}\lambda^{2}T\int_{\Omega_{T}} \xi_{\pm}^{3}|\psi_{\pm}|^{2}\mathrm{~d}x\mathrm{~d}t\nonumber\\
				&\leqslant & C\varepsilon^{2} s^{3}\lambda^{2}\int_{\Omega_{T}} \xi_{\pm}^{3}|\psi_{\pm}|^{2}\mathrm{~d}x\mathrm{~d}t,
				\label{c13}
			\end{eqnarray}
			for all $s\geq \frac{s_{1}}{\varepsilon}T$. We also have 
			\begin{eqnarray*}
				I_{\pm}^{23}
				&= & \mp \varepsilon s^{2}\lambda\int_{\Omega_{T}}  \partial_{t}\alpha_{\pm}\xi_{\pm}\nabla\eta\cdot\nabla|\psi_{\pm}|^{2}\mathrm{~d}x\mathrm{~d}t\\
				&=& \mp \varepsilon s^{2}\lambda\int_{\Gamma_{T}}  \partial_{t}\alpha\xi\partial_{\n}\eta|\psi|^{2}\mathrm{~d}\sigma\mathrm{~d}t \pm \varepsilon s^{2}\lambda\int_{\Omega_{T}} (\nabla( \partial_{t}\alpha_{\pm}).\nabla\eta) \xi_{\pm}|\psi_{\pm}|^{2}\mathrm{~d}x\mathrm{~d}t\\
				&&+ \varepsilon s^{2}\lambda^{2}\int_{\Omega_{T}} \partial_{t}\alpha_{\pm}|\nabla\eta|^{2}\xi_{\pm}|\psi_{\pm}|^{2}\mathrm{~d}x\mathrm{~d}t \pm \varepsilon s^{2}\lambda\int_{\Omega_{T}}  \partial_{t}\alpha_{\pm}\Delta\eta\xi_{\pm}|\psi_{\pm}|^{2}\mathrm{~d}x\mathrm{~d}t\\
				&:=& (I_{\pm}^{23})_{1}+(I_{\pm}^{23})_{2}+(I_{\pm}^{23})_{3}+(I_{\pm}^{23})_{4}.
			\end{eqnarray*}
			Using $|\partial_{t}\alpha_{\pm}|\leqslant  T\xi_{\pm}^{2}$ and $|\nabla(\partial_{t}\alpha_{\pm})|\leqslant \|\nabla\eta\|_{\infty} \lambda T\xi_{\pm}^{2}$, we obtain 
			\begin{eqnarray*}
				|(I_{\pm}^{23})_{2}+(I_{\pm}^{23})_{3}+(I_{\pm}^{23})_{4}|
				&\leqslant & C\varepsilon^{2} s^{3}\lambda^{2}\int_{\Omega_{T}}\xi_{\pm}^{3}|\psi_{\pm}|^{2}\mathrm{~d}x\mathrm{~d}t,\;\text{for all}\;\lambda\geq 1\;\text{and}\; s\geq \frac{s_{1}}{\varepsilon}T,
			\end{eqnarray*}
			it follows that
			\begin{eqnarray}
				I_{\pm}^{23} &\geq& \mp \varepsilon s^{2}\lambda\int_{\Gamma_{T}}  \partial_{t}\alpha\xi\partial_{\n}\eta|\psi|^{2}\mathrm{~d}x\mathrm{~d}t-C\varepsilon^{2} s^{3}\lambda^{2}\int_{\Omega_{T}}\xi_{\pm}^{3}|\psi_{\pm}|^{2}\mathrm{~d}x\mathrm{~d}t. \label{c14}
			\end{eqnarray}
			By intergation in time and $\psi_{\pm}(\cdot,0)=0=\psi_{\pm}(\cdot,T)$, we obtain
			\begin{eqnarray*}
				I_{\pm}^{33} 
				&=&\frac{s}{2}\int_{\Omega_{T}}\partial_{t}\alpha_{\pm}\partial_{t}|\psi_{\pm}|^{2}\mathrm{~d}x\mathrm{~d}t=-\frac{s}{2}\int_{\Omega_{T}}\partial^{2}_{t}\alpha_{\pm}|\psi_{\pm}|^{2}\mathrm{~d}x\mathrm{~d}t.
			\end{eqnarray*}
			Since $|\partial^{2}_{t}\alpha_{\pm}|\leqslant 2T^{2}\xi_{\pm}^{3}$, then 
			\begin{eqnarray}
				\left|I_{\pm}^{33}\right|&\leqslant& sT^{2}\int_{\Omega_{T}}\xi_{\pm}^{3}|\psi_{\pm}|^{2}\mathrm{~d}x\mathrm{~d}t\nonumber\\
				&\leqslant &  C\varepsilon^{2}s^{3}\lambda^{2}\int_{\Omega_{T}}\xi_{\pm}^{3}|\psi_{\pm}|^{2}\mathrm{~d}x\mathrm{~d}t,\;\text{for all}\;\lambda\geq 1\;\text{and}\; s\geq \frac{s_{1}}{\varepsilon}T. \label{c15}
			\end{eqnarray}
			Using $\B(x,t)\cdot\n(x)=0$ on $\Gamma_{T}$, the last term of this sub-step is processed by integration by parts as follow
			\begin{eqnarray}
				I_{\pm}^{43}
				&=& \frac{s}{2}\int_{\Omega_{T}}\partial_{t}\alpha_{\pm}\B\cdot\nabla|\psi_{\pm}|^{2}\mathrm{~d}x\mathrm{~d}t\nonumber\\
				&=&-\frac{s}{2}\int_{\Omega_{T}}(\nabla(\partial_{t}\alpha_{\pm})\cdot\B)|\psi_{\pm}|^{2}\mathrm{~d}x\mathrm{~d}t-\frac{s}{2}\int_{\Omega_{T}}\partial_{t}\alpha_{\pm}(\nabla\cdot\B)|\psi_{\pm}|^{2}\mathrm{~d}x\mathrm{~d}t. \nonumber
			\end{eqnarray}
			Since $|\partial_{t}\alpha_{\pm}|\leqslant T\xi_{\pm}^{2}$ and $|\nabla(\partial_{t}\alpha_{\pm})|\leqslant \|\nabla\eta\|_{\infty}\lambda T\xi_{\pm}^{2}$ and $\xi_{\pm}\geq\frac{4}{T^{2}}$, then 
			\begin{eqnarray}
				|I_{\pm}^{43} |&\leqslant& C s T(\lambda\|\B\|_{\infty}+\|\nabla\cdot\B\|_{\infty})\int_{\Omega_{T}}\xi_{\pm}^{2}|\psi_{\pm}|^{2}\mathrm{~d}x\mathrm{~d}t \nonumber\\
				&\leqslant& C s T^{3}(\lambda\|\B\|_{\infty}+\|\nabla\cdot\B\|_{\infty})\int_{\Omega_{T}}\xi_{\pm}^{3}|\psi_{\pm}|^{2}\mathrm{~d}x\mathrm{~d}t \nonumber \\
				&\leqslant & C\varepsilon^{2} s^{3}\lambda^{2}\int_{\Omega_{T}} \xi_{\pm}^{3}|\psi_{\pm}|^{2}\mathrm{~d}x\mathrm{~d}t, \label{cc16}
			\end{eqnarray}
			for all $\lambda\geq 1$ and $s\geq  \frac{s_{1}}{\varepsilon}(T+T^{2})(\|\B\|_{\infty}^{1/2}+\|\nabla\cdot\B\|_{\infty}^{1/2})$.\\
			From \eqref{c13}, \eqref{c14}, \eqref{c15} and \eqref{cc16}, we deduce for $\lambda\geq 1$ and $s\geq  \frac{s_{1}}{\varepsilon}(T+T^{2})(1+\|\B\|_{\infty}^{1/2}+\|\nabla\cdot\B\|_{\infty}^{1/2})$ that
			\begin{eqnarray}
				I_{\pm}^{3}&=&\sum_{i=1}^{4} I_{\pm}^{i3} \geq  -C\varepsilon^{2} s^{3}\lambda^{2}\int_{\Omega_{T}} \xi_{\pm}^{3}|\psi_{\pm}|^{2}\mathrm{~d}x\mathrm{~d}t \mp \varepsilon s^{2}\lambda\int_{\Gamma_{T}}  \partial_{t}\alpha\xi\partial_{\n}\eta|\psi|^{2}\mathrm{~d}\sigma\mathrm{~d}t.\nonumber\\
				\label{c16}
			\end{eqnarray}.\\
			\textcolor{blue}{\textbf{Steb 2d. Estimate from below of $I_{\pm}^{4}$.}}
			Let us now consider double products involving $(L_{2}\psi_{\pm})_{4}$. Firstly 
			\begin{eqnarray*}
				I_{\pm}^{14}=\pm 2\varepsilon s^{2}\lambda^{3}\int_{\Omega_{T}} |\nabla\eta|^{2}\xi_{\pm}^{2}|\psi_{\pm}|^{2}(\nabla\eta\cdot\B)\mathrm{~d}x\mathrm{~d}t.
			\end{eqnarray*}
			Using $\xi_{\pm}\geq \frac{4}{T^{2}}$, we obtain 
			\begin{eqnarray}
				|I_{\pm}^{14}|
				&\leqslant & C\varepsilon s^{2}\lambda^{3}T^{2}\|\B\|_{\infty}\int_{\Omega_{T}} \xi_{\pm}^{3}|\psi_{\pm}|^{2}\mathrm{~d}x\mathrm{~d}t \nonumber\\
				&\leqslant & C\varepsilon^{2} s^{3}\lambda^{3}\int_{\Omega_{T}} \xi_{\pm}^{3}|\psi_{\pm}|^{2}\mathrm{~d}x\mathrm{~d}t,\; \mbox{for all}\; s\geq \frac{s_{1}}{\varepsilon}T^{2}\|\B\|_{\infty}. \label{T1}
			\end{eqnarray}
			Using $\nabla\eta(x)\cdot\B(x,t)=0$ on $\Gamma_{T}$ and integration by parts, one has
			\begin{eqnarray*}
				I_{\pm}^{24}
				&=& \varepsilon s^{2}\lambda^{2}\int_{\Omega_{T}}(\nabla\eta\cdot\B)\xi_{\pm}^{2}\nabla\eta\cdot\nabla|\psi_{\pm}|^{2}\mathrm{~d}x\mathrm{~d}t\\
				&=& -\varepsilon s^{2}\lambda^{2}\int_{\Omega_{T}}\nabla\cdot((\nabla\eta\cdot\B)\nabla\eta)\xi_{\pm}^{2}|\psi_{\pm}|^{2}\mathrm{~d}x\mathrm{~d}t\\
				&& \mp 2\varepsilon s^{2}\lambda^{3}\int_{\Omega_{T}}(\nabla\eta\cdot\B)|\nabla\eta|^{2}\xi_{\pm}^{2}|\psi_{\pm}|^{2}\mathrm{~d}x\mathrm{~d}t.
			\end{eqnarray*}	
			
			From formula \eqref{cv3} and the fact that $\xi_{\pm}\geq\frac{4}{T^{2}}$, we get
			\begin{eqnarray}
				|I_{\pm}^{24}| 
				&\leqslant & C\varepsilon s^{2}\lambda^{3}(\|\B\|_{\infty}+\|\nabla \B\|_{\infty})\int_{\Omega_{T}}\xi_{\pm}^{2}|\psi_{\pm}|^{2}\mathrm{~d}x\mathrm{~d}t\nonumber\\
				&\leqslant & C\varepsilon s^{2}\lambda^{3}T^{2}(\|\B\|_{\infty}+\|\nabla \B\|_{\infty})\int_{\Omega_{T}}\xi_{\pm}^{3}|\psi_{\pm}|^{2}\mathrm{~d}x\mathrm{~d}t \nonumber\\
				&\leqslant & C\varepsilon^{2} s^{3}\lambda^{3}\int_{\Omega_{T}}\xi_{\pm}^{3}|\psi_{\pm}|^{2}\mathrm{~d}x\mathrm{~d}t, \label{T2}
			\end{eqnarray}	
			for all $s\geq \frac{s_{1}}{\varepsilon}T^{2}(\|\B\|_{\infty}+\|\nabla \B\|_{\infty})$.\\
			Next, integration by parts in time, $\psi_{\pm}(\cdot,0)=0=\psi_{\pm}(\cdot,T)$ and the regularity of $\B$, gives
			\begin{eqnarray*}
				I_{\pm}^{34}		
				&=&\mp \frac{s\lambda}{2}\int_{\Omega_{T}}(\nabla\eta\cdot\B) \xi_{\pm}\partial_{t}|\psi_{\pm}|^{2}\mathrm{~d}x\mathrm{~d}t \\
				&=& \pm \frac{s\lambda}{2}\int_{\Omega_{T}}
				((\nabla\eta\cdot\partial_{t}\B)\xi_{\pm}+(\nabla\eta\cdot\B)\partial_{t}\xi_{\pm}
				\partial_{t}((\nabla\eta\cdot\B)|\psi_{\pm}|^{2}\mathrm{~d}x\mathrm{~d}t.
			\end{eqnarray*}
			Using $|\partial_{t}\xi_{\pm}|\leqslant T\xi_{\pm}^{2}$ and $\xi_{\pm}\geq \frac{4}{T^{2}}$, we obtain 
			\begin{eqnarray}
				|I_{\pm}^{34}|		
				&\leqslant& Cs\lambda 
				\int_{\Omega_{T}}(\|\partial_{t}\B\|_{\infty}\xi_{\pm}+T\|\B\|_{\infty}\xi_{\pm}^{2})|\psi_{\pm}|^{2}\mathrm{~d}x\mathrm{~d}t \nonumber\\
				&\leqslant & Cs\lambda (\|\partial_{t}\B\|_{\infty}T^{4}+\|\B\|_{\infty}T^{3})
				\int_{\Omega_{T}}\xi_{\pm}^{3}|\psi_{\pm}|^{2}\mathrm{~d}x\mathrm{~d}t \nonumber\\
				&\leqslant & C\varepsilon^{2}s^{3}\lambda
				\int_{\Omega_{T}}\xi_{\pm}^{3}|\psi_{\pm}|^{2}\mathrm{~d}x\mathrm{~d}t, \label{T3}
			\end{eqnarray}
			for all $s\geq\frac{s_{1}}{\varepsilon}(T^{2}\|\partial_{t}\B\|_{\infty}^{1/2}+(T^{2}+T)\|\B\|_{\infty}^{1/2})$.\\
			Using the fact that $\B(x,t)\cdot\n(x)=0$ on $\Gamma_{T}$ and integration by parts, we also have
			\begin{eqnarray*}
				I_{\pm}^{44} 
				&=& \mp \frac{s\lambda}{2}
				\int_{\Omega_{T}}(\nabla\eta\cdot\B)\xi_{\pm}\B\cdot\nabla|\psi_{\pm}|^{2}\mathrm{~d}x\mathrm{~d}t\\
				&=& \pm \frac{s\lambda}{2}
				\int_{\Omega_{T}}\nabla\cdot((\nabla\eta\cdot\B)\B)\xi_{\pm}|\psi_{\pm}|^{2}\mathrm{~d}x\mathrm{~d}t + \frac{s\lambda^{2}}{2}
				\int_{\Omega_{T}}|\nabla\eta\cdot\B|^{2}\xi_{\pm}|\psi_{\pm}|^{2}\mathrm{~d}x\mathrm{~d}t.
			\end{eqnarray*}
			From $\xi_{\pm}\geq\frac{4}{T^{2}}$ and formula \eqref{cv4}, one has
			\begin{eqnarray}
				|I_{\pm}^{44}| &\leqslant& Cs\lambda (\|\B\|^{2}_{\infty}+\|\nabla \B\|^{2}_{\infty}+\|\nabla\cdot\B\|^{2}_{\infty}+\lambda\|\B\|^{2}_{\infty})
				\int_{\Omega_{T}}\xi_{\pm}|\psi_{\pm}|^{2}\mathrm{~d}x\mathrm{~d}t \nonumber\\
				&\leqslant & Cs\lambda^{2} (\|\B\|^{2}_{\infty}+\|\nabla \B\|^{2}_{\infty}+\|\nabla\cdot\B\|^{2}_{\infty})
				\int_{\Omega_{T}}\xi_{\pm}|\psi_{\pm}|^{2}\mathrm{~d}x\mathrm{~d}t \nonumber\\
				&\leqslant & Cs\lambda^{2} (\|\B\|^{2}_{\infty}+\|\nabla \B\|^{2}_{\infty}+\|\nabla\cdot\B\|^{2}_{\infty})T^{4}
				\int_{\Omega_{T}}\xi_{\pm}^{3}|\psi_{\pm}|^{2}\mathrm{~d}x\mathrm{~d}t \nonumber\\
				&\leqslant & C\varepsilon^{2}s^{3}\lambda^{2} 
				\int_{\Omega_{T}}\xi_{\pm}^{3}|\psi_{\pm}|^{2}\mathrm{~d}x\mathrm{~d}t, \label{T4}
			\end{eqnarray}
			for all $s\geq\frac{s_{1}}{\varepsilon}T^{2}(\|\B\|_{\infty}+\|\nabla \B\|_{\infty}+\|\nabla\cdot\B\|_{\infty})$.\\
			By considering \eqref{T1}, \eqref{T2}, \eqref{T3} and \eqref{T4}, we obtain 
			\begin{eqnarray}
				I_{\pm}^{4}=\sum_{i=1}^{4}I_{\pm}^{i4}\geq -C\varepsilon^{3}s^{3}\lambda^{2} 
				\int_{\Omega_{T}}\xi_{\pm}^{3}|\psi_{\pm}|^{2}\mathrm{~d}x\mathrm{~d}t, \label{T}
			\end{eqnarray}
			for all  $s\geq\frac{s_{1}}{\varepsilon}(T^{2}(\|\B\|_{\infty}+\|\nabla \B\|_{\infty}+\|\nabla\cdot\B\|_{\infty}+\|\partial_{t}\B\|_{\infty}^{1/2}+\|\B\|_{\infty}^{1/2})+T\|\B\|_{\infty}^{1/2})$ and $\lambda\geq 1$.\\
			\textcolor{blue}{\textbf{Step 3. First Conclusion.}}\\
			Taking into account \eqref{c8}, \eqref{c12}, \eqref{c16} and \eqref{T}, for any $\lambda\geq \lambda_{1}$ and \\$s\geq \frac{s_{1}}{\varepsilon}(T^{2}+T)(1+\|\B\|_{\infty}+\|\nabla\cdot \B\|_{\infty}+\|\nabla \B\|_{\infty}+\|\B\|_{\infty}^{1/2}+\|\nabla\cdot\B\|_{\infty}^{1/2}+\|\partial_{t}\B\|_{\infty}^{1/2}) $, we obtain
			\begin{eqnarray}
				&&(L_{1}\psi_{\pm},L_{2}\psi_{\pm})_{L^{2}(\Omega_{T})}=\sum_{j=1}^{4}I_{\pm}^{j} \nonumber\\
				&& \quad\geq  c\,\varepsilon^{2}s^{3}\lambda^{4}\int_{\Omega_{T}}\xi_{\pm}^{3}|\psi_{\pm}|^{2}\mathrm{~d}x\mathrm{~d}t-C\varepsilon^{2}s^{3}\lambda^{4}\int_{\omega^{'}\times (0,T)}\xi_{\pm}^{3}|\psi_{\pm}|^{2}\mathrm{~d}x\mathrm{~d}t \nonumber\\
				&& \quad\quad\mp\varepsilon^{2}s^{3}\lambda^{3}\int_{\Gamma_{T}} |\nabla\eta|^{2}\partial_{\n}\eta\xi^{3}|\psi|^{2}\mathrm{~d}\sigma\mathrm{~d}t\;\;\nonumber\\
				&& \quad\quad- 2\varepsilon^{2}s\lambda^{2}\int_{\Gamma_{T}}|\nabla\eta|^{2}\xi\psi\partial_{\n}\psi_{\pm}\mathrm{~d}\sigma\mathrm{~d}t-C\varepsilon^{2}s^{2}\lambda^{4}\int_{\Omega_{T}} \xi_{\pm}^{2}|\psi_{\pm}|^{2}\mathrm{~d}x\mathrm{~d}t \nonumber\\
				&& \quad\quad + c\,\varepsilon^{2}s\lambda^{2}\int_{\Omega_{T}}|\nabla\psi_{\pm}|^{2}\xi_{\pm}\mathrm{~d}x\mathrm{~d}t-C\varepsilon^{2}s\lambda^{2}\int_{\omega^{'}\times (0,T)}|\nabla\psi_{\pm}|^{2}\xi_{\pm}\mathrm{~d}x\mathrm{~d}t \nonumber\\
				&& \quad\quad \mp \varepsilon^{2}s\lambda\int_{\Gamma_{T}}\partial_{\n}\eta \xi|\partial_{\n}\psi_{\pm}|^{2}\mathrm{~d}\sigma\mathrm{~d}t \pm              \varepsilon^{2}s\lambda\int_{\Gamma_{T}}\xi\partial_{\n}\eta|\nabla_{\Gamma}\psi|^{2}\mathrm{~d}\sigma\mathrm{~d}t\nonumber \nonumber\\
				&&\quad\quad +\varepsilon\int_{\Gamma_{T}}\partial_{t}\psi\partial_{\n}\psi_{\pm}\mathrm{~d}\sigma\mathrm{~d}t +\varepsilon\int_{\Gamma_{T}}(\B\cdot\nabla_{\Gamma}\psi)\partial_{\n}\psi_{\pm}\mathrm{~d}\sigma\mathrm{~d}t\;\; \nonumber \\
				&&\quad\quad - C\varepsilon^{2} s^{3}\lambda^{2}\int_{\Omega_{T}} \xi_{\pm}^{3}|\psi_{\pm}|^{2}\mathrm{~d}x\mathrm{~d}t \pm\varepsilon s^{2}\lambda\int_{\Gamma_{T}}  \partial_{t}\alpha\xi\partial_{\n}\eta|\psi|^{2}\mathrm{~d}\sigma\mathrm{~d}t\;\;\nonumber \\
				&&\quad\quad -C\varepsilon^{2}s^{3}\lambda^{3}\int_{\Omega_{T}}\xi_{\pm}^{3}|\psi_{\pm}|^{2}\mathrm{~d}x\mathrm{~d}t. \label{conc1}
			\end{eqnarray}
			Using $\xi_{\pm}\geq\frac{4}{T^{2}}$ and let us increase $s_{1}$ and $\lambda_{1}$, the fifth, twelfth and fourteenth terms in the right hand side of \eqref{conc1} can be absorbed by the first term in the right hand side of \eqref{conc1}, then
			\begin{eqnarray}
				&&(L_{1}\psi_{\pm},L_{2}\psi_{\pm})_{L^{2}(\Omega_{T})}=\sum_{j=1}^{4}I_{\pm}^{j} \nonumber\\
				&& \quad\geq  c\,\varepsilon^{2}s^{3}\lambda^{4}\int_{\Omega_{T}}\xi_{\pm}^{3}|\psi_{\pm}|^{2}\mathrm{~d}x\mathrm{~d}t-C\varepsilon^{2}s^{3}\lambda^{4}\int_{\omega^{'}\times (0,T)}\xi_{\pm}^{3}|\psi_{\pm}|^{2}\mathrm{~d}x\mathrm{~d}t \nonumber\\
				&& \quad\quad\mp\varepsilon^{2}s^{3}\lambda^{3}\int_{\Gamma_{T}} |\nabla\eta|^{2}\partial_{\nu}\eta\xi^{3}\psi^{2}\mathrm{~d}\sigma\mathrm{~d}t- 2\varepsilon^{2}s\lambda^{2}\int_{\Gamma_{T}}|\nabla\eta|^{2}\xi\psi\partial_{\nu}\psi_{\pm}\mathrm{~d}\sigma\mathrm{~d}t \nonumber\\
				&& \quad\quad + c\,\varepsilon^{2}s\lambda^{2}\int_{\Omega_{T}}|\nabla\psi_{\pm}|^{2}\xi_{\pm}\mathrm{~d}x\mathrm{~d}t-C\varepsilon^{2}s\lambda^{2}\int_{\omega^{'}\times (0,T)}|\nabla\psi_{\pm}|^{2}\xi_{\pm}\mathrm{~d}x\mathrm{~d}t \nonumber\\
				&& \quad\quad \mp \varepsilon^{2}s\lambda\int_{\Gamma_{T}}\partial_{\nu}\eta \xi|\partial_{\nu}\psi_{\pm}|^{2}\mathrm{~d}\sigma\mathrm{~d}t \pm                  \varepsilon^{2}s\lambda\int_{\Gamma_{T}}\xi\partial_{\nu}\eta|\nabla_{\Gamma}\psi|^{2}\mathrm{~d}\sigma\mathrm{~d}t\nonumber \nonumber\\
				&&\quad\quad +\varepsilon\int_{\Gamma_{T}}\partial_{t}\psi\partial_{\nu}\psi_{\pm}\mathrm{~d}\sigma\mathrm{~d}t +\varepsilon\int_{\Gamma_{T}}(\B\cdot\nabla_{\Gamma}\psi)\partial_{\nu}\psi_{\pm}\mathrm{~d}\sigma\mathrm{~d}t \nonumber \\
				&&\quad\quad  \pm\varepsilon s^{2}\lambda\int_{\Gamma_{T}}  \partial_{t}\alpha\xi\partial_{\nu}\eta|\psi|^{2}\mathrm{~d}\sigma\mathrm{~d}t. \label{conc2}
			\end{eqnarray}
			From \eqref{c2} and \eqref{conc2}, we obtain
			\begin{align}
				&	\|L_{1}\psi_{\pm}\|^{2}_{L^{2}(\Omega_{T})}+\|L_{2}\psi_{\pm}\|^{2}_{L^{2}(\Omega_{T})}+c\,\varepsilon^{2}s^{3}\int_{\Omega_{T}}\xi_{\pm}^{3}|\psi_{\pm}|^{2}\mathrm{~d}x\mathrm{~d}t+c\,\varepsilon^{2}s\lambda^{2}\int_{\Omega_{T}}|\nabla\psi_{\pm}|^{2}\xi_{\pm}\mathrm{~d}x\mathrm{~d}t \nonumber\\
				& \quad\quad + 2\;\textbf{I}_{\pm} \nonumber\\
				&\leqslant C\left( \|L_{3}\psi_{\pm}\|^{2}_{L^{2}(\Omega_{T})}+\varepsilon^{2}s^{3}\lambda^{4}\int_{\omega^{'}\times (0,T)}\xi_{\pm}^{3}|\psi_{\pm}|^{2}\mathrm{~d}x\mathrm{~d}t+\varepsilon^{2}s\lambda^{2}\int_{\omega^{'}\times (0,T)}|\nabla\psi_{\pm}|^{2}\xi_{\pm}\mathrm{~d}x\mathrm{~d}t  \right), \nonumber\\ \label{c37}
			\end{align}
			where $\textbf{I}_{\pm}$ the sum of all integrals on the boundary that are in the right hand side of \eqref{conc2}.\\
			From \eqref{c0}, we obtain 
			\begin{eqnarray}
				& &	\|L_{3}\psi_{\pm}\|^{2}_{L^{3}(\Omega_{T})} \leqslant  C\left(\varepsilon^{2}s^{2}\lambda^{2}\int_{\Omega_{T}}\xi_{\pm}^{2} |\psi_{\pm}|^{2}\mathrm{~d}x\mathrm{~d}t
				+\varepsilon^{2}s^{2}\lambda^{4}\iint_{\Omega_{T}}\xi_{\pm}^{2} |\psi_{\pm}|^{2}\mathrm{~d}x\mathrm{~d}t \right. \nonumber\\
				& & \left.
				+\|\nabla\cdot\B\|_{\infty}^{2}\int_{\Omega_{T}} |\psi_{\pm}|^{2}\mathrm{~d}x\mathrm{~d}t\right).\nonumber\\
				&\leqslant & C\left(\varepsilon^{2}s^{2}\lambda^{4}\int_{\Omega_{T}}\xi_{\pm}^{2} |\psi_{\pm}|^{2}\mathrm{~d}x\mathrm{~d}t +\|\nabla\cdot\B\|_{\infty}^{2}\int_{\Omega_{T}} |\psi_{\pm}|^{2}\mathrm{~d}x\mathrm{~d}t\right).
				\label{c38}
			\end{eqnarray}
			The first term in the right hand side of \eqref{c38} is aborbed by the third term in the left hand side of \eqref{c37} for $s\geq s_{1}T^{2}$, so is the second one, by taking $s\geq \frac{s_{1}}{\varepsilon^{\frac{2}{3}}}T^{2}\|\nabla\cdot\B\|_{\infty}^{\frac{2}{3}}$.\\
			Finally, we obtain 
			\begin{align}
				& \|L_{1}\psi_{\pm}\|^{2}_{L^{2}(\Omega_{T})}+\|L_{2}\psi_{\pm}\|^{2}_{L^{2}(\Omega_{T})}+c\,\varepsilon^{2}s^{3}\lambda^{4}\int_{\Omega_{T}}\xi_{\pm}^{3}|\psi_{\pm}|^{2}\mathrm{~d}x\mathrm{~d}t\nonumber \\
				&+c\,\varepsilon^{2}s\lambda^{2}\int_{\Omega_{T}}|\nabla\psi_{\pm}|^{2}\xi_{\pm}\mathrm{~d}x\mathrm{~d}t + 2\;\textbf{I}_{\pm} \nonumber\\
				&\leqslant C\left( \varepsilon^{2}s^{3}\lambda^{4}\int_{\omega^{'}\times (0,T)}\xi_{\pm}^{3}|\psi_{\pm}|^{2}\mathrm{~d}x\mathrm{~d}t+\varepsilon^{2}s\lambda^{2}\int_{\omega^{'}\times (0,T)}|\nabla\psi_{\pm}|^{2}\xi_{\pm}\mathrm{~d}x\mathrm{~d}t \right), \label{c21}
			\end{align}
			for the following estimates of the parameters $s$, $\lambda$ and $\varepsilon$
			\begin{eqnarray}
				&&\lambda\geq \lambda_{1},\quad s \geq  \frac{s_{1}}{\varepsilon}(T^{2}+T)\mathcal{B}_{T},\quad \varepsilon\in (0,1) \label{parameters}\\
				&&\mathcal{B}_{T}=1+\|\B\|_{\infty}+\|\nabla\cdot \B\|_{\infty}+\|\nabla \B\|_{\infty}+\|\B\|_{\infty}^{1/2}+\|\nabla\cdot\B\|_{\infty}^{1/2}+\|\nabla\cdot\B\|_{\infty}^{2/3}+\|\partial_{t}\B\|_{\infty}^{1/2}. \nonumber\\
				\nonumber 
			\end{eqnarray}
			\textcolor{blue}{\textbf{Step 4. Elimination of the integral of $|\nabla\psi_{\pm}|^{2}$ on the right-hand side of \eqref{c21}.}}
			We start by adding integral of $|\Delta \psi_{\pm}|^{2}$ to the left-hand side of \eqref{c21}, so that we can eliminate the last term in the right-hand side of \eqref{c21}.\\
			Using \eqref{c-1}, $\xi_{\pm}\geq \frac{4}{T^{2}}$, $s\geq s_{1}T^{2}$ and $|\partial_{t}\alpha_{\pm}|\leqslant T\xi_{\pm}^{2}$, we
			obtain
			\begin{eqnarray}
				\varepsilon^{2}s^{-1}\int_{\Omega_{T}}\xi^{-1}_{\pm}|\Delta\psi_{\pm}|^{2}\mathrm{~d}x\mathrm{~d}t
				&\leqslant & C\left(\varepsilon^{2}s^{3}\lambda^{4}\int_{\Omega_{T}}\xi_{\pm}^{3}|\psi_{\pm}|^{2}\mathrm{~d}x\mathrm{~d}t+s T^{2}\int_{\Omega_{T}}\xi^{3}_{\pm}|\psi_{\pm}|^{2}\mathrm{~d}x\mathrm{~d}t\right. \nonumber\\
				&+& \left. s\lambda^{2}T^{4}\|\B\|_{\infty}^{2}\int_{\Omega_{T}}\xi_{\pm}^{3}|\psi_{\pm}|^{2}\mathrm{~d}x\mathrm{~d}t+ \|L_{2}\psi_{\pm}\|^{2}_{L^{2}(\Omega_{T})}\right).\nonumber
			\end{eqnarray}
			Hence
			\begin{eqnarray}
				\varepsilon^{2}s^{-1}\int_{\Omega_{T}}\xi^{-1}_{\pm}|\Delta\psi_{\pm}|^{2}\mathrm{~d}x\mathrm{~d}t
				&\leqslant & C\left(\varepsilon^{2}s^{3}\lambda^{4}\int_{\Omega_{T}}\xi_{\pm}^{3}|\psi_{\pm}|^{2}\mathrm{~d}x\mathrm{~d}t+\|L_{2}\psi_{\pm}\|^{2}_{L^{2}(\Omega_{T})}\right),\nonumber\\
			\end{eqnarray}
			for  all $s\geq \frac{s_{1}}{\varepsilon}T^{2}\left(1+\|\B\|_{\infty}\right)$. Consequently, we deduce from \eqref{c21} that
			\begin{align}
				& \varepsilon^{2}s^{-1}\int_{\Omega_{T}}\xi^{-1}_{\pm}|\Delta\psi_{\pm}|^{2}\mathrm{~d}x\mathrm{~d}t+\varepsilon^{2}s^{3}\lambda^{4}\int_{\Omega_{T}}\xi_{\pm}^{3}|\psi_{\pm}|^{2}\mathrm{~d}x\mathrm{~d}t\nonumber\\
				&+\varepsilon^{2}s\lambda^{2}\int_{\Omega_{T}}|\nabla\psi_{\pm}|^{2}\xi_{\pm}\mathrm{~d}x\mathrm{~d}t + c\;\textbf{I}_{\pm} \nonumber\\
				&\leqslant C\left( \varepsilon^{2}s^{3}\lambda^{4}\int_{\omega^{'}\times (0,T)}\xi_{\pm}^{3}|\psi_{\pm}|^{2}\mathrm{~d}x\mathrm{~d}t+\varepsilon^{2}s\lambda^{2}\int_{\omega^{'}\times (0,T)}|\nabla\psi_{\pm}|^{2}\xi_{\pm}\mathrm{~d}x\mathrm{~d}t \right), 
				\label{c23}
			\end{align}
			for $\lambda_{1}$ and $s_{1}$ large enough in \eqref{parameters}.\\
			\par 
			We are now ready to eliminate the last term in the right-hand side of \eqref{c21}, let's introduce a cut-off function, between $\omega^{'}$ and $\omega$. More precisely, let $\theta\in \mathcal{C}^{2}(\omega)$ be a positive function such that $\theta=1$ in $\omega^{'}$. We then have
			\begin{eqnarray*}
				C\varepsilon^{2}s\lambda^{2}\int_{\omega^{'}\times(0,T)}\xi_{\pm} |\nabla\psi_{\pm}|^{2}\mathrm{~d}x\mathrm{~d}t
				&\leqslant & C\varepsilon^{2}s\lambda^{2}\int_{\omega\times(0,T)}\theta\xi_{\pm} |\nabla\psi_{\pm}|^{2}\mathrm{~d}x\mathrm{~d}t.
			\end{eqnarray*}
			Since $\theta=0$ on $\partial\omega$, integrating by parts, we obtain 
			\begin{eqnarray*}
				C\varepsilon^{2}s\lambda^{2}\int_{\omega\times(0,T)}\theta\xi_{\pm} |\nabla\psi_{\pm}|^{2}\mathrm{~d}x\mathrm{~d}t &=& -C\varepsilon^{2}s\lambda^{2}\int_{\omega\times(0,T)} \nabla\cdot(\theta\xi_{\pm} \nabla\psi_{\pm})\psi_{\pm}\mathrm{~d}x\mathrm{~d}t\\
				&=& -C\varepsilon^{2}s\lambda^{2}\int_{\omega\times(0,T)} \theta\xi_{\pm} \Delta\psi_{\pm}\psi_{\pm}\mathrm{~d}x\mathrm{~d}t
				\\& &\mp C\varepsilon^{2}s\lambda^{3}\int_{\omega\times(0,T)} \theta\xi_{\pm}(\nabla\eta\cdot \nabla\psi_{\pm})\psi_{\pm}\mathrm{~d}x\mathrm{~d}t\\
				&&- C\varepsilon^{2}s\lambda^{2}\int_{\omega\times(0,T)} \xi_{\pm}(\nabla\theta\cdot\nabla\psi_{\pm})\psi_{\pm}\mathrm{~d}x\mathrm{~d}t.
			\end{eqnarray*}
			Applying Hölder's inequality to each integral, we deduce that
			\begin{eqnarray*}
				& & C\varepsilon^{2}s\lambda^{2}\int_{\omega^{'}\times(0,T)}\xi_{\pm} |\nabla\psi_{\pm}|^{2}\mathrm{~d}x\mathrm{~d}t \leqslant C \varepsilon^{2}s^{3}\lambda^{4}\int_{\omega\times (0,T)}\xi_{\pm}^{3}|\psi_{\pm}|^{2}\mathrm{~d}x\mathrm{~d}t\\
				&&+ \frac{1}{2}\left(\varepsilon^{2}s^{-1}\int_{\omega\times (0,T)}\xi^{-1}_{\pm}|\Delta\psi_{\pm}|^{2}\mathrm{~d}x\mathrm{~d}t+ \varepsilon^{2}s\lambda^{2}\int_{\omega\times (0,T)}|\nabla\psi_{\pm}|^{2}\xi_{\pm}\mathrm{~d}x\mathrm{~d}t
				\right),
			\end{eqnarray*}
			where we have used that $\lambda\geq 1$, $\xi_{\pm}\geq \frac{4}{T^{2}}$, $s\geq s_{1}T^{2}$.\\
			From this estimate and \eqref{c23} we deduce that 
			\begin{align}
				& \varepsilon^{2}s^{-1}\int_{\Omega_{T}}\xi^{-1}_{\pm}|\Delta\psi_{\pm}|^{2}\mathrm{~d}x\mathrm{~d}t+\varepsilon^{2}s^{3}\lambda^{4}\int_{\Omega_{T}}\xi_{\pm}^{3}|\psi_{\pm}|^{2}\mathrm{~d}x\mathrm{~d}t\nonumber\\ &+\varepsilon^{2}s\lambda^{2}\int_{\Omega_{T}}|\nabla\psi_{\pm}|^{2}\xi_{\pm}\mathrm{~d}x\mathrm{~d}t + c\;\textbf{I}_{\pm} \leqslant C \varepsilon^{2}s^{3}\lambda^{4}\int_{\omega\times (0,T)}\xi_{\pm}^{3}|\psi_{\pm}|^{2}\mathrm{~d}x\mathrm{~d}t, \nonumber
			\end{align} 
			for $\lambda_{1}$ and $ s_{1}$ large enough in \eqref{parameters}.\\ 
			In particular, we have
			\begin{align}
				& \varepsilon^{2}s^{3}\lambda^{4}\int_{\Omega_{T}}\xi_{\pm}^{3}|\psi_{\pm}|^{2}\mathrm{~d}x\mathrm{~d}t +\varepsilon^{2}s\lambda^{2}\int_{\Omega_{T}}|\nabla\psi_{\pm}|^{2}\xi_{\pm}\mathrm{~d}x\mathrm{~d}t + c\;\textbf{I}_{\pm} \nonumber\\
				&\leqslant C \varepsilon^{2}s^{3}\lambda^{4}\int_{\omega\times (0,T)}\xi_{\pm}^{3}|\psi_{\pm}|^{2}\mathrm{~d}x\mathrm{~d}t. \label{c24}
			\end{align} 
			\textcolor{blue}{\textbf{Step 5. Return to the original function.}}
			By summing \eqref{c24} for $i=+,-$, in particular we obtain 
			\begin{align}
				&\varepsilon^{2}s^{3}\lambda^{4}\int_{\Omega_{T}}\xi_{+}^{3}|\psi_{+}|^{2}\mathrm{~d}x\mathrm{~d}t +\varepsilon^{2}s\lambda^{2}\int_{\Omega_{T}}|\nabla\psi_{+}|^{2}\xi_{+}\mathrm{~d}x\mathrm{~d}t + c\,\left(\textbf{I}_{+}+\textbf{I}_{-}\right) \nonumber\\
				&\leqslant C\varepsilon^{2}s^{3}\lambda^{4}\int_{\omega\times (0,T)}\left(\xi_{+}^{3}|\psi_{+}|^{2}+\xi_{-}^{3}|\psi_{-}|^{2}\right)\mathrm{~d}x\mathrm{~d}t. 
				\label{c25}
			\end{align}
			From the definitions of $\xi_{\pm}$ and $\alpha_{\pm}$, we have 
			\begin{equation*}
				\xi_{-}\leqslant \xi_{+},\quad \text{and}\quad \alpha_{+}\leqslant \alpha_{-},\;\;\text{in}\;\Omega_{T},
			\end{equation*}
			then, the estimate \eqref{c25} becomes 
			\begin{align}
				&\varepsilon^{2}s^{3}\lambda^{4}\int_{\Omega_{T}}\xi_{+}^{3}|\psi_{+}|^{2}\mathrm{~d}x\mathrm{~d}t +\varepsilon^{2}s\lambda^{2}\int_{\Omega_{T}}|\nabla\psi_{+}|^{2}\xi_{+}\mathrm{~d}x\mathrm{~d}t + c\,\left(\textbf{I}_{+}+\textbf{I}_{-}\right) \nonumber\\
				&\leqslant C\varepsilon^{2}s^{3}\lambda^{4}\int_{\omega\times (0,T)}\xi_{+}^{3}|\psi_{+}|^{2}\mathrm{~d}x\mathrm{~d}t.\label{c26}
			\end{align}
			\par 
			Concerning the integrals on the boundary, using
			\begin{eqnarray*}
				\partial_{\n}\psi_{\pm} = \pm s\lambda\xi_{\pm}\psi_{\pm}\partial_{\n}\eta =\pm s\lambda\xi\psi\partial_{\n}\eta,
			\end{eqnarray*}
			we obtain 
			\begin{eqnarray*}
				\partial_{\nu}\psi_{+}+\partial_{\nu}\psi_{-}=0,\quad\text{and}\quad |\partial_{\nu}\psi_{+}|^{2}=|\partial_{\nu}\psi_{-}|^{2}.
			\end{eqnarray*}
			Thus, we easily get that
			$$\textbf{I}_{+}+\textbf{I}_{-}=0.$$
			Therfore, \eqref{c26} becomes 
			\begin{align}
				s^{3}\lambda^{4}\int_{\Omega_{T}}\xi_{+}^{3}|\psi_{+}|^{2}\mathrm{~d}x\mathrm{~d}t +s\lambda^{2}\int_{\Omega_{T}}|\nabla\psi_{+}|^{2}\xi_{+}\mathrm{~d}x\mathrm{~d}t  \leqslant Cs^{3}\lambda^{4}\int_{\omega\times (0,T)}\xi_{+}^{3}|\psi_{+}|^{2}\mathrm{~d}x\mathrm{~d}t.\label{cc27}
			\end{align} 
			\par
			Now, we  turn back to our original function $\varphi$. From \eqref{c-4} and \eqref{cc27}, we deduce that 
			\begin{align}
				&s^{3}\lambda^{4}\int_{\Omega_{T}}\exp(-2s\alpha_{+})\xi_{+}^{3}|\varphi|^{2}\mathrm{~d}x\mathrm{~d}t +s\lambda^{2}\int_{\Omega_{T}}|\nabla\psi_{+}|^{2}\xi_{+}\mathrm{~d}x\mathrm{~d}t \nonumber\\
				&\quad\quad \leqslant Cs^{3}\lambda^{4}\int_{\omega\times (0,T)}\exp(-2s\alpha_{+})\xi_{+}^{3}|\varphi|^{2}\mathrm{~d}x\mathrm{~d}t. 
				\label{cc28}
			\end{align}
			For $\nabla\varphi$, we use the identity 
			\begin{eqnarray*}
				\exp(-s\alpha_{+})\nabla\varphi=\nabla\psi_{+} - s\lambda\nabla\eta\xi_{+}\psi_{+}
			\end{eqnarray*}
			and applying the triangular inequality to this last identity, we find 
			\begin{eqnarray}
				& &	s\lambda^{2}\int_{\Omega_{T}}\exp(-2s\alpha_{+})\xi_{+}|\nabla\varphi|^{2}\mathrm{~d}x\mathrm{~d}t \nonumber\\
				&&\leqslant  C\left(s\lambda^{2}\int_{\Omega_{T}}\xi_{+}|\nabla\psi_{+}|^{2}\mathrm{~d}x\mathrm{~d}t+s^{3}\lambda^{4}\int_{\Omega_{T}}\exp(-2s\alpha_{+})\xi_{+}^{3}|\varphi|^{2}\mathrm{~d}x\mathrm{~d}t\right). \label{cc29}
			\end{eqnarray}
			Consequently, the integral of $|\nabla\psi_{+}|^{2}$ in \eqref{cc28} can be changed by the integral of
			$|\nabla\varphi|^{2}$ in to the left-hand side of \eqref{cc29}:
			\begin{align}
				&s^{3}\lambda^{4}\int_{\Omega_{T}}\exp(-2s\alpha_{+})\xi_{+}^{3}|\varphi|^{2}\mathrm{~d}x\mathrm{~d}t +s\lambda^{2}\int_{\Omega_{T}}|\nabla\varphi|^{2}\xi_{+}\mathrm{~d}x\mathrm{~d}t \nonumber\\
				&\quad\quad \leqslant Cs^{3}\lambda^{4}\int_{\omega\times (0,T)}\exp(-2s\alpha_{+})\xi_{+}^{3}|\varphi|^{2}\mathrm{~d}x\mathrm{~d}t,
				\label{cc30}
			\end{align}
			for $\lambda_{1}$ and $s_{1}$ large enough in \eqref{parameters}.
		\end{proof}
		Under the same conditions of Proposition \ref{P4}, we have the following corollary:
		\begin{corollary} 
			Let $\B\in W^{1,\infty}(\Omega\times (0,+\infty))^{d}$ be a tangential transport, and assume the same conditions in Proposition  \ref{P4}. Then, for all $0<\kappa<1$, there are two constants $C$ independent of $\varepsilon$ and  $C_{1}>0$ independent of $\varepsilon$ and $T$ such that for all $t\in [0,\kappa\,T]$, we have 
			\begin{equation}
				\int_{\Omega}|\varphi(x,t)|^{2}\mathrm{~d} x \leqslant C\exp\left(\frac{C_{1}}{\varepsilon}\left(1+\frac{1}{T}\right)\right)\int_{\omega\times (0,T)}|\varphi(x,t)|^{2}\mathrm{~d} x \mathrm{~d}t.\label{io1}
			\end{equation}
		\end{corollary}
		\begin{proof}
			By the Carleman estimate \eqref{Carleman}, we obtain
			\begin{equation}
				\begin{aligned}
					\int_{\Omega\times(0, T)} \exp(-2 s \alpha_{+})\xi_{+}^3|\varphi|^2 \mathrm{d} x \mathrm{~d}t \leq C  \int_{\omega \times(0, T)} \exp(-2 s \alpha_{+}) \xi_{+}^3|\varphi|^2 \mathrm{~d} x \mathrm{~d}t,
				\end{aligned}
			\end{equation}
			where $\lambda=\lambda_{1}$ and $s = \frac{s_{1}}{\varepsilon}(T^{2}+T)\mathcal{B}_{\infty}$; see definition of $\mathcal{B}_{T}$ in \eqref{constante s}. Note that  $\mathcal{B}_{\infty}$ is well defined and independent of $T$, since $\B\in W^{1,\infty}(\Omega\times (0,\infty))^{d}$.\\
			Taking lower and upper estimates with respect to $x$ of the weight functions, we get 
			\begin{equation}
				\begin{aligned}
					\int_{\Omega\times(0, T)}\underline{f}(t) |\varphi(x,t)|^2 \mathrm{d} x \mathrm{~d} t \leq C\int_{\omega \times(0, T)}\overline{f}(t)|\varphi(x,t)|^2 \mathrm{~d} x \mathrm{~d}t,
				\end{aligned}
			\end{equation}
			where 
			$$\begin{cases*} \displaystyle
				\underline{f}(t):=\exp \left(-2 s \max _{x \in \overline{\Omega}} \alpha_{+}(x,t)\right) \min _{x \in \overline{\Omega}} \xi^{3}_{+}(x,t),\\
				\displaystyle \overline{f}(t):=\exp \left(-2 s \min _{x \in \overline{\Omega}} \alpha_{+}(x,t)\right) \max _{x \in \overline{\Omega}} \xi^{3}_{+}(x,t).
			\end{cases*}$$
			For $\lambda_{1}$ and $s_{1}$ large enoughh, it is easy to check that the function $g$ admits a maximum on $[0,T]$ at $t=\frac{T}{2}$ and
			$f$ admits a minimum on $\left[\kappa T,\left(\frac{1+\kappa}{2}\right)T\right]$ at $t=\kappa_{0}T$ where $\kappa_{0}\in \left[\kappa,\frac{1+\kappa}{2}\right]$. Hence 
			\begin{eqnarray}
				& & \int_{\Omega\times \left(\kappa T,\left(\frac{1+\kappa}{2}\right)T\right)} |\varphi(x,t)|^2 \mathrm{d} x \mathrm{~d} t \leq  C\frac{\overline{f}\left(\frac{T}{2}\right)}{\underline{f}(\kappa_{0}T)}\int_{\omega \times(0, T)}|\varphi(x,t)|^2 \mathrm{~d} x \mathrm{~d}t \nonumber\\
				& & \quad\quad\leq  \exp\left(\frac{C_{1}}{\varepsilon}\left(1+\frac{1}{T}\right)\right)\int_{\omega \times(0, T)}|\varphi(x,t)|^2 \mathrm{~d} x \mathrm{~d}t. \label{e2}
			\end{eqnarray}
			for somme $C_{1}$ depends only on $\Omega$, $\omega$, $\kappa$ and $\mathcal{B}_{\infty}$.\\
			Using the dissipation estimate \eqref{A3}
			for the solutions of \eqref{s2}, we get for all $ 0\leqslant t\leqslant s\leqslant T$,
			\begin{eqnarray*}
				\int_{\Omega}|\varphi(x,t)|^{2}\mathrm{~d}x&\leqslant& \exp\left(C_{\B}T \right)
				\int_{\Omega}|\varphi(x,s)|^{2}\mathrm{~d}x,
			\end{eqnarray*}
			where $C_{\B}=\|\nabla\cdot\B\|_{L^{\infty}(\Omega_{T})}$. By integrating this inequality on $\left(\kappa T,\left(\frac{1+\kappa}{2}\right)T \right)$, we obtain, for all $t\in [0,\kappa T]$,
			\begin{eqnarray}
				\int_{\Omega}|\varphi(x,t)|^{2}dx\leqslant \frac{2\exp\left(C_{\B}T \right)}{(1-\kappa)T}
				\int_{\Omega\times \left(\kappa T,\left(\frac{1+\kappa}{2}\right)T\right)}|\varphi(x,s)|^{2}\mathrm{~d} x \mathrm{~d}s.\label{e3}
			\end{eqnarray}
			From \eqref{e2} and \eqref{e3}, we obtain \eqref{io1}.
		\end{proof}
		\section{Proof of main results} \label{Section 7}
		\subsection{Proof of Theorem \ref{m1}} \label{Section 71}
		According to the estimates \eqref{io1} and \eqref{disspation 3} we have:\\
		For all $0<\kappa<1$, there are two constants $C$ independent of $\varepsilon$ and  $C_{1}>0$ independent of $\varepsilon$ and $T$ such that for all $t\in [0,\kappa T]$, we have 
		\begin{equation}
			\int_{\Omega}|\varphi(x,t)|^{2}\mathrm{~d} x \leqslant C\exp\left(\frac{C_{1}}{\varepsilon}\left(1+\frac{1}{T}\right)\right)\int_{\omega\times (0,T)}|\varphi(x,t)|^{2}\mathrm{~d} x \mathrm{~d}t \label{e4}
		\end{equation}
		and there are two constants $C^{'}$ independent of $\varepsilon$ and  $C_{0}>0$ independent of $\varepsilon$ and $T$ such that
		for any integer $m$ such that $1\leqslant m\leqslant\frac{T}{T_{0}}$, there exists $C>0$  independent of $\varepsilon$ such that for all $t\in [mT_{0},T]$, we have 
		\begin{equation} 
			\|\varphi(\cdot, 0)\|^{2}_{L^{2}(\Omega)}\leqslant C^{'}\left(\exp\left(\frac{-mC_{0}}{\varepsilon}\right)\|\varphi(\cdot, t)\|^{2}_{L^{2}(\Omega)}+\|\varphi\|^{2}_{L^{2}(0,T;L^{2}(\omega))}\right). \label{e5}
		\end{equation}
		For $0<\kappa<1$ fixed, taking $m:=\left[\frac{C_{1}}{C_{0}} \right] +1$ where $\left[\frac{C_{1}}{C_{0}} \right]$ denotes the integer part of $\frac{C_{1}}{C_{0}}$ and $\rho_{0}:=\max\left(\frac{m}{\kappa},\frac{1}{T_{0}\left(m\frac{C_{0}}{C_{1}}-1\right)}\right)$.\\
		Let $T\geq \rho_{0} T_{0}$. Then $mT_{0}\leq \kappa T$, so applying \eqref{e4} and \eqref{e5}, we get
		\begin{equation} \label{e6}
			\|\varphi(\cdot, 0)\|^{2}_{L^{2}(\Omega)}\leqslant C^{'}\left(C\exp\left(\frac{C_{1}\left(1+1/T\right)-mC_{0}}{\varepsilon}\right)+1\right)\|\varphi\|^{2}_{L^{2}(0,T;L^{2}(\omega))}.
		\end{equation} 
		On the other hand, the choice of $\rho_{0}$ implies that  $C_{1}\left(1+1/T\right)-mC_{0}\leq 0$, for all $T\geq \rho_{0} T_{0}$. Finally combining \eqref{e6} and \eqref{cost of control}, we obtain \eqref{cost1} for $\varepsilon>0$ small enough and $T\geq \rho_{0} T_{0}$.
		\subsection{Proof of Theorem \ref{m2}} \label{Section 72}
		Let $x_{0}\in\Omega$ such that the condition \ref{C3m2} of the Theorem \ref{m2} is fulfilled. From the continuity of $x\mapsto\varPhi(t,T,x)$ uniform in $t$, there exists $r_{0}>0$ such that
		$$\varPhi(t,T,x)\in \Omega\setminus\overline{\omega},\;\;\forall t\in [0,T],\;\forall x\in \overline{B}(x_{0},4r_{0}).$$
		Consider $\varphi_{T}\in D(B(x_{0},r_{0}))$ and $\varphi$ be the strong solution of \eqref{s2} with the data $\varphi_{T}$, and $\theta$ the function defined in Lemma \ref{lemma} with this choice of $x_{0}$, $r=r_{0}$, $t_{1}=0$, $t_{2}=T$. Let $\vartheta_{1}$ and $\vartheta_{2}$ be regular functions such that
		\begin{equation*}
			\left\{
			\begin{aligned}
				\vartheta_{1}(x,t) &=0, & & \forall (x,t)\in \mathcal{D}_{2r_{0}}(x_{0},0,T), \\
				\vartheta_{1}(x,t) &=1, & & \forall (x,t)\notin \mathcal{D}_{3r_{0}}(x_{0},0,T)	
			\end{aligned}
			\right.\quad\text{and}\quad 	
			\left\{
			\begin{aligned}
				\vartheta_{2}(x,t) &=0, & & \forall (x,t)\in \mathcal{D}_{3r_{0}}(x_{0},0,T), \\
				\vartheta_{2}(x,t) &=1, & & \forall (x,t)\notin \mathcal{D}_{4r_{0}}(x_{0},0,T).	
			\end{aligned}
			\right.
		\end{equation*}
		For reasons of simplicity, we will divide the proof into three steps.\\
		\textcolor{blue}{\textbf{Step 1:}} We will show that there are $C_{1},C_{2}>0$ independent of $\varepsilon$ such that
		\begin{equation} \label{7}
			\int_{\Omega_{T}}|\vartheta_{1}\varphi|^{2}\mathrm{~d}x\mathrm{~d}t+
			\int_{\Omega_{T}}|\nabla(\vartheta_{1}\varphi)|^{2}\mathrm{~d}x\mathrm{~d}t\leqslant C_{1}\exp\left(-\frac{C_{2}}{\varepsilon}\right)\int_{\Omega}|\psi(x,T)|^{2}\mathrm{~d}x,\\
		\end{equation}
		where $\psi=\exp\left(\frac{\theta}{\varepsilon}\right)\varphi$ and $T, \varepsilon$ are small enough.\\
		For all $(x,t)\in\text{supp}(\vartheta_{1})$, we have $\theta(x,t)\geq c_{0}r_{0}^{2}$, then 
		\begin{eqnarray} \label{17}
			\int_{\Omega_{T}}|\vartheta_{1}\varphi|^{2}\mathrm{~d}x\mathrm{~d}t&\leqslant &  \exp\left(-\frac{2c_{0}r_{0}^{2}}{\varepsilon}\right)\int_{\Omega_{T}}|\vartheta_{1}v|^{2}\mathrm{~d}x\mathrm{~d}t \nonumber\\
			&\leqslant & \|\vartheta_{1}\|^{2}_{\infty}\exp\left(-\frac{2c_{0}r_{0}^{2}}{\varepsilon}\right)\int_{\Omega_{T}}|\psi|^{2}\mathrm{~d}x\mathrm{~d}t.
		\end{eqnarray}
		From $\nabla(\vartheta_{1}\varphi)=\exp\left(-\frac{\theta}{\varepsilon}\right)\left(\psi\left(\nabla\vartheta_{1}-\vartheta_{1}\frac{\nabla\theta}{\varepsilon}\right)+\vartheta_{1}\nabla \psi\right)$, we obtain
		\begin{eqnarray} \label{18}
			\int_{\Omega_{T}}|\nabla(\vartheta_{1}\varphi)|^{2}\mathrm{~d}x\mathrm{~d}t&\leqslant &  4\exp\left(-\frac{2c_{0}r_{0}^{2}}{\varepsilon}\right)\left(\|\nabla\vartheta_{1}\|^{2}_{\infty}+\frac{\|\vartheta_{1}\nabla\theta\|^{2}_{\infty}}{\varepsilon^{2}}\right)\int_{\Omega_{T}}|\psi|^{2}\mathrm{~d}x\mathrm{~d}t \nonumber\\
			&& + 2\exp\left(-\frac{2c_{0}r_{0}^{2}}{\varepsilon}\right) \|\vartheta_{1}\|^{2}_{\infty}\int_{\Omega_{T}}|\nabla \psi|^{2}\mathrm{~d}x\mathrm{~d}t. 
		\end{eqnarray}
		By application of Agmon inequality \eqref{A1}, we get 
		\begin{equation} \label{19}
			\int_{\Omega_{T}}|\psi|^{2}\mathrm{~d}x\mathrm{~d}t\leqslant T\exp\left(\frac{C}{\varepsilon}T\right)\int_{\Omega}|\psi(x,T)|^{2}\mathrm{~d}x
		\end{equation}
		and 
		\begin{equation} \label{20}
			\int_{\Omega_{T}}|\nabla \psi|^{2}\mathrm{~d}x\mathrm{~d}t\leqslant \exp\left(\frac{C}{\varepsilon}T\right)\int_{\Omega}|\psi(x,T)|^{2}\mathrm{~d}x.
		\end{equation}
		Taking $0<T<\frac{2c_{0}r_{0}^{2}}{C}$ and $\varepsilon$ small enough, from \eqref{17}, \eqref{18}, \eqref{19} and \eqref{20}, we obtain \eqref{7}.\\
		\textcolor{blue}{\textbf{Step 2:}} We will prove that there are
		$C_{1},C_{2}>0$ independent of $\varepsilon$ such that
		\begin{equation} \label{10}
			\int_{0}^{T}\int_{\omega}|\varphi(x,t)|^{2}\mathrm{~d}x\mathrm{~d}t\leqslant C_{1}\exp\left(-\frac{C_{2}}{\varepsilon}\right)\int_{\Omega}|\varphi_{T}(x)|^{2}\mathrm{~d}x,
		\end{equation}
		for $T$ and $\varepsilon$ are small enough.\\
		Let $\phi=\vartheta_{2}\varphi$, then $\varPhi$ is the solution of the system $\mathcal{S}(\varnothing,0,T,f,0,\varepsilon,\B)$, 
		where $f=(\partial_t \vartheta_{2}+\varepsilon\Delta\vartheta_{2}-\B(x,t)\cdot\nabla\vartheta_{2})\varphi+2\varepsilon\nabla\vartheta_{2}\cdot\nabla\varphi$.\\
		By the estimate \eqref{dissw}, we obtain  
		\begin{equation} \label{diss11}
			\int_{\Omega}|\phi(x,t)|^{2}\mathrm{~d}x\leqslant C\exp\left(\frac{CT}{\varepsilon}\right)\int_{0}^{T}\int_{\Omega}|f|^{2}\mathrm{~d}x\mathrm{~d}t,\;\;\forall t\in [0,T],
		\end{equation}
		for $\varepsilon$ is small enough and $C>0$ independent of $\varepsilon$.\\
		Since $\vartheta_{1}=1$ on the supports of the functions  $\partial_t \vartheta_{2}$ and $\nabla\vartheta_{2}$, we obtain from \eqref{7} the existence of the constants $C_{1},C_{2}>0$ independent of $\varepsilon$ such that
		\begin{eqnarray} \label{diss10}
			\int_{\Omega_{T}}|f|^{2}\mathrm{~d}x\mathrm{~d}t &\leqslant& C_{1}\exp\left(-\frac{C_{2}}{\varepsilon}\right)\int_{\Omega}|\psi(x,T)|^{2}\mathrm{~d}x \nonumber\\
			&=& C_{1}\exp\left(-\frac{C_{2}}{\varepsilon}\right)\int_{\Omega}|\varphi_{T}(x)|^{2}\mathrm{~d}x,
		\end{eqnarray}
		because $\theta(\cdot,T)=0$ on the support of $\varphi_{T}$.\\
		Using \eqref{diss11} and \eqref{diss10}, there are constants $C_{1},C_{2}>0$ independent of $\varepsilon$ such that
		\begin{equation*}
			\int_{0}^{T}\int_{\Omega}|\phi(x,t)|^{2}\mathrm{~d}x\mathrm{~d} t\leqslant C_{1}\exp\left(-\frac{C_{2}}{\varepsilon}\right)\int_{\Omega}|\varphi_{T}(x)|^{2}\mathrm{~d}x,
		\end{equation*}
		for $T>0$ and $\varepsilon$ are small enough, since $\vartheta_{2}=1$ on $\omega\times (0,T)$, we deduce that \eqref{10} is true.\\\\
		\textcolor{blue}{\textbf{Step 3:}} Finally, since $\varphi$ is the weak solution of the adjoint system \eqref{s2}, from \eqref{weak solution}, we obtain, for all $v\in D(]0,T[)$
		\begin{eqnarray*}
			\int_{0}^{T}\left(\int_{\Omega}\varphi(x,t)\d x\right) v^{'}(t)\d t=-\int_{0}^{T}\mathfrak{a}_{w}(T-t,\varphi(t),v(t))\d t=0,
		\end{eqnarray*}
		since $v$ independent of $x$. Hence $t\longmapsto \displaystyle\int_{\Omega}\varphi(x,t)\d x$ is weakly differentiable and $\displaystyle\frac{\d}{\d t} \int_{\Omega}\varphi(x,t)\d x=0$, thus
		\begin{equation} 
			\int_{\Omega}\varphi(x,0)\mathrm{~d}x=\int_{\Omega}\varphi(x,T)\mathrm{~d}x. \label{e1}
		\end{equation}
		Choosing the initial data $\varphi_{T}\in D(B(x_{0},r_{0}))$ such that $\displaystyle\int_{\Omega}\varphi_{T}(x)\mathrm{~d}x\neq 0$.\\ By Hölder's inequality and \eqref{e1}, we have 
		\begin{equation} \label{11}
			\int_{\Omega}|\varphi(x,0)|^{2}\mathrm{~d}x\geq\frac{1}{|\Omega|}\left|	\int_{\Omega}\varphi_{T}(x)\mathrm{~d}x\right|^{2}.
		\end{equation}
		Finally, combining \eqref{10}, \eqref{11} and \eqref{cost of control}, we obtain \eqref{cost2}.
\section{Open problems}
		Here we present some relevant problems that remain open:
		\begin{itemize}
			\item Let $\mathbf{B}\in W^{1,\infty}(\Omega_T)$. Under the assumption:
			  $$\exists T_{0}\in (0,T),\;\exists r_{0}>0,\;\mbox{such that}\; (T,T_{0},r_{0},\B,\Omega)\; \mbox{satisfies}\; \eqref{Flushing Condition}\; \mbox{for}\; \mathbb{R}^{d}\setminus\overline{\Omega}.$$
			 Show that
			\begin{eqnarray*}
				\mathcal{K}(\varepsilon,T,\Omega,\omega)\leq C_{1}\exp\left(-\frac{C_{2}}{\varepsilon}\right),
			\end{eqnarray*}
			where $\mathcal{K}$ the cost of the null controllability of \eqref{s1}. The difficulty in this case is to establish a Carleman estimate without the condition $\B(x,t)\cdot\n(x)=0$, as this condition requires the trajectories of \eqref{OD} still in $\overline{\Omega}$.\\
			\item Study the null-controllability cost of $\eqref{s1}_{1}$ with  speed $\B\in L^{\infty}(\Omega_{T})$ and any kind of boundary conditions. As far as we know, this problem is open even for Dirichlet boundary conditions.  
		\end{itemize}


\begin{thebibliography}{plain}
			
			\bibitem{arendt2014maximal} Arendt, W., Dier, D., Laasri, H., \& Ouhabaz, E. M.:  Maximal regularity for evolution equations governed by non-autonomous forms. Advances in Differential Equations. 19, 1043-1066 (2014)
			\bibitem{arendt2011dirichlet} Arendt, W., \& ter Elst, A. F.: The Dirichlet-to-Neumann operator on rough domains. Journal of Differential Equations, 251(8), 2100-2124 (2011)
			\bibitem{barcena2021cost} Bárcena-Petisco, J. A.:  Cost of null controllability for parabolic equations with vanishing diffusivity and a transport term. ESAIM: Control, Optimisation and Calculus of Variations, 27, 106 (2021)
			\bibitem{brezis1983analyse} Brezis, H.: Analyse fonctionnelle. Théorie et applications. Masson, Paris (1983)
			\bibitem{carreno2016cost} Carre{\~n}o, N., \& Guzm{\'a}n, P.: On the cost of null controllability of a fourth-order parabolic equation. Journal of Differential Equations, 261(11), 6485-6520 (2016)
			\bibitem{carreno2015non}  Carre{\~n}o, N., \& Guerrero, S.: On the non-uniform null controllability of a linear KdV equation. Asymptotic Analysis, 94(1-2), 33-69 (2015)
			
			
			\bibitem{coron2005singular}  Coron, J-M., \&  Guerrero, S.: Singular optimal control: a linear 1-D parabolic--hyperbolic example. Asymptotic Analysis, 44(3-4), 237-257 (2005)
			
			\bibitem{Lions} Dautray, R., \& Lions, J. L.:  Mathematical analysis and numerical methods for science and technology: volume 1 physical origins and classical methods. Springer Science \& Business Media (2012)
			\bibitem{ettahri:hal-04131920} Et-Tahri, F., Bárcena-Petisco, J. A., Boutaayamou, I., \& Maniar, L.:  Asymptotic behavior of null-controllablity cost for parabolic equations with vanishing diffusivity and a transport term. https://hal.science/hal-04131920/file/Cost\%20of\%20null\%20controllability2.pdf (2023)
			\bibitem{fernandez2006exact} Fernández-Cara, E., González-Burgos, M., Guerrero, S., \& Puel, J. P.:  Exact controllability to the trajectories of the heat equation with Fourier boundary conditions: the semilinear case. ESAIM: Control, Optimisation and Calculus of Variations, 12(3), 466-483 (2006)
			\bibitem{fursikov1996controllability} Fursikov, A., \& Imanuvilov, O. Y.:  Controllability of Evolution Equations, Lecture Notes. 34, Seoul National University. Korea (1996)
			\bibitem{glass2010complex} Glass, O.: A complex-analytic approach to the problem of uniform controllability of a transport equation in the vanishing viscosity limit. Journal of Functional Analysis, 258(3), 852-868 (2010)
			\bibitem{grisvard1985elliptic} Grisvard, P.: Elliptic Problems in Nonsmooth Domains. Monographs and Studies in Mathematics, vol. 24. Pitman, Boston 49--52 (1985)
			\bibitem{guerrero2007singular} Guerrero, S., \& Lebeau, G.: Singular optimal control for a transport-diffusion equation. Communications in Partial Differential Equations, 32(12), 1813-1836 (2007)
			\bibitem{kato2013perturbation} Kato, T.:  Perturbation theory for linear operators (Vol. 132). Springer Science \& Business Media,New York (2013)
			\bibitem{laurent2021uniform} Laurent, C., \& Léautaud, M.: On uniform observability of gradient flows in the vanishing viscosity limit. Journal de l’École polytechnique—Mathématiques, 8, 439-506 (2021)
			\bibitem{laurent2022uniform} Laurent, C., \& Léautaud, M.: On uniform controllability of 1D transport equations in the vanishing viscosity limit. Comptes Rendus. Mathématique, 361(G1), 265-312 (2023)
			\bibitem{lions1988controlabilite} Lions, J. L.:  Contr{\^o}labilit{\'e} exacte, stabilisation et perturbations de systemes distribu{\'e}s. Tome 1. Contr{\^o}labilit{\'e} exacte. Rech. Math. Appl, 8 (1988)
			\bibitem{lissy2012link} Lissy, P.: A link between the cost of fast controls for the 1-D heat equation and the uniform controllability of a 1-D transport-diffusion equation. Comptes Rendus Mathematique, 350(11-12), 591-595 (2012)
			\bibitem{lissy2014application} Lissy, P.: An application of a conjecture due to Ervedoza and Zuazua concerning the observability of the heat equation in small time to a conjecture due to Coron and Guerrero concerning the uniform controllability of a convection–diffusion equation in the vanishing viscosity limit. Systems \& Control Letters, 69, 98-102 (2014)
			\bibitem{lissy2015explicit} Lissy, P.: Explicit lower bounds for the cost of fast controls for some 1-D parabolic or dispersive equations, and a new lower bound concerning the uniform controllability of the 1-D transport–diffusion equation. Journal of Differential Equations, 259(10), 5331-5352 (2015)
			\bibitem{lopez2000null} L{\'o}pez, A., Zhang, X., \& Zuazua, E.: Null controllability of the heat equation as singular limit of the exact controllability of dissipative wave equations. Journal de mathématiques pures et appliquées, 79(8), 741-808 (2000)
			\bibitem{phung2002null} Phung, K. D.: Null controllability of the heat equation as singular limit of the exact controllability of dissipative wave equation under the Bardos-Lebeau-Rauch geometric control condition. Computers \& Mathematics with Applications, 44(10-11), 1289-1296 (2002)
			\bibitem{russell1978controllability}  Russell, D. L.:  Controllability and stabilizability theory for linear partial differential equations: recent progress and open questions. Siam Review, 20(4), 639-739 (1978)
			\bibitem{showalter2013monotone} Showalter, R. E.:  Monotone operators in Banach space and nonlinear partial differential equations (Vol. 49). American Mathematical Soc..(2013)
			
		\end{thebibliography}
	\end{document}